\documentclass[10pt]{article}
\usepackage[utf8]{inputenc}
\usepackage[T1]{fontenc}
\usepackage[english]{babel}

\usepackage{amsmath, amssymb}

\usepackage{geometry}
\usepackage{authblk}

\usepackage{amsmath, amsfonts, amssymb}
\usepackage{bm} % boldmath symbols

\usepackage[left]{lineno}
\usepackage[algo2e, ruled, lined, longend]{algorithm2e}

\usepackage[dvipsnames]{xcolor} % must be loaded before tikz
\usepackage{listofitems}

\usepackage{hyperref}
\usepackage{csquotes}

\usepackage[backend=biber,
            style=numeric-comp,
            sorting=nty,
            giveninits=true,
            maxbibnames=99,
            date=year,      % no months in references
        hyperref]{biblatex}

\AtEveryBibitem{
  \clearfield{urlyear}
  \clearfield{urlmonth}
  \clearfield{month}
  \clearfield{day}
  \clearfield{eprintclass}
  \clearlist{language}
}
\hypersetup{
    hidelinks,
    colorlinks
}

\usepackage{subcaption}
\usepackage{siunitx}
\usepackage{booktabs}

\usepackage{tikz}
\usetikzlibrary{shapes.geometric,shapes.arrows}
\usetikzlibrary{positioning,fit}

\tikzset{
  nn node/.style={thick,circle,draw=NavyBlue,minimum size=20,inner sep=0.5,outer sep=0.6},
  nn node in/.style={nn node,Green!20!black,draw=Green,fill=Green!25},
  nn node hidden/.style={nn node,NavyBlue!20!black,draw=NavyBlue,fill=NavyBlue!20},
  nn node convol/.style={nn node,orange!20!black,draw=orange,fill=orange!20},
  nn node out/.style={nn node,Red!20!black,draw=Red,fill=Red!20},
  nn connect/.style={thick,Gray}, %,line cap=round
  nn connect arrow/.style={-{Latex[length=4,width=3.5]},thick,NavyBlue,shorten <=0.5,shorten >=1},
  nn node 1/.style={nn node in}, % node styles, numbered for easy mapping with \nstyle
  nn node 2/.style={nn node hidden},
  nn node 3/.style={nn node out},
  roundedbox/.style={rectangle,draw,thick,fill=white,outer sep=0pt,align=center,minimum height=3.5em,rounded corners}
}

\tikzset{
  unet txt/.style={text=Gray!60!black, font=\large},
  unet img/.style={rectangle, draw, semithick, inner sep=0pt, outer sep=0pt},
  unet arrow/.style={unet txt, font=\Large, single arrow, draw, semithick, inner sep=3pt, outer sep=0pt, minimum height=40pt},
  unet in/.style={draw=Green, fill=Green!25},
  unet out/.style={draw=Red, fill=Red!25},
  unet enc/.style={draw=NavyBlue, fill=NavyBlue!20},
  unet trans/.style={draw=Cyan, fill=Cyan!20},
  unet updown/.style={draw=orange, fill=orange!20},
  unet dec/.style={draw=Plum, fill=Plum!20},
  unet img in/.style={unet img, unet in},
  unet img out/.style={unet img, unet out},
  unet img enc/.style={unet img, unet enc},
  unet img dec/.style={unet img, unet dec},
  unet img concat/.style={unet img, draw=none, inner sep=0pt},
  unet arrow in/.style={unet arrow, unet trans},
  unet arrow out/.style={unet arrow, unet out},
  unet arrow enc/.style={unet arrow, unet trans, thick, inner ysep=10pt, minimum width=50pt},
  unet arrow updown/.style={unet arrow, unet updown, minimum height=20pt},
  unet arrow bot/.style={unet arrow in},
  unet arrow concat/.style={unet arrow enc},
  unet skip/.style={->, draw=Gray!60!black, line width=3pt},
}

\newlength{\clength}
\newlength{\hzero}
\newlength{\hone}
\newlength{\htwo}
\newlength{\hthree}
\newlength{\hfour}

\DeclareMathOperator*{\argmin}{arg\,min}
\DeclareMathOperator{\Span}{span}

\def\FGMRES{CNN-FGMRES } %if I want to give a special name of FGMRES preconditioned by CNN
\def\Fc{{\mathcal{F}}}
\def\Pc{{\mathcal{P}}}
\def\Uc{{\mathcal{U}}}
\def\K{L} 
\newcommand{\bz}{\mathbf{z}}
\newcommand{\bv}{\mathbf{v}}
\newcommand{\bw}{\mathbf{w}}
\newcommand{\bu}{\mathbf{u}}
\newcommand{\bb}{\mathbf{b}}
\newcommand{\bc}{\mathbf{c}}
\newcommand{\bx}{\mathbf{x}}
\newcommand{\by}{\mathbf{y}}
\newcommand{\bt}{\mathbf{t}}
\newcommand{\br}{\mathbf{r}}
\newcommand{\be}{\mathbf{e}}

\newcommand{\etab}{\eta_b}
\def\tol{\varepsilon}
\def\its{its}
\def\ET{ET}
\def\maxit{maxIts}

\newcommand{\ignore}[1]{}
\newcommand{\NN}{\mathcal{N}_\theta}

\newcommand{\loss}{\mathcal{L}_\theta}

\newcommand{\cost}{J}

\newcommand{\simiid}{\stackrel{\textnormal{iid}}{\sim}}
\newcommand{\bydef}{:=}

\providecommand{\vec}{}
\renewcommand{\vec}[1]{\mathbf{#1}}

\newcommand{\R}{\mathbb{R}}
\newcommand{\C}{\mathbb{C}}

\newcommand{\vx}{\vec{x}}

\newcommand{\inner}{\textnormal{inner}}

\newcommand{\pml}{\textnormal{pml}}

\newcommand{\Opml}{\Omega_{\pml}}
\newcommand{\Oinner}{\Omega_{\inner}}

\providecommand{\L}{}
\renewcommand{\L}{L}
\newcommand{\Lpml}{\ell}
\newcommand{\Linner}{\L_{\inner}}

\def\its{Its}
\def\ET{ET}
\def\maxit{maxIts}

\begin{document}

\title{Neural operator preconditioning from mixed dataset for the Helmholtz equations: Application to transcranial ultrasound}

\author[1,*]{Yanfei Xiang}

\affil[1]{IRMA, Université de Strasbourg, MACARON, Inria, France.}
\affil[*]{Corresponding author (\url{yfxiang0amber@gmail.com} and \url{yanfei.xiang@math.unistra.fr}).}

\date{}

\maketitle

\begin{abstract}
This work develops a neural operator preconditioned subspace method for sequences of linear systems arising from the discretization of the two-dimensional Helmholtz equation in transcranial ultrasound applications. The problem involves strongly heterogeneous, patient-dependent velocity fields that induce severe wave distortion and pose significant challenges for standard iterative solvers. Building on neural network preconditioning framework of Giraud et al. (HAL RR-9593, 2025) and the idealized skull dataset used for the learned optimizer of Stanziola et al. (JCP 441, 2021), neural operator preconditioners are trained on six mixed velocity-source datasets combining randomized source configurations and idealized skull-based velocity fields with random noise. The proposed mixed-dataset strategy aims to improve both computational efficiency and generalization across varying configurations.
The neural operator is trained on a coarse grid using a physics-informed loss based on the relative residual of the discrete Helmholtz equation and is incorporated as a nonlinear preconditioner within flexible GMRES (FGMRES). Numerical experiments demonstrate that the resulting hybrid method efficiently solves practical transcranial ultrasound problems on grids 64 times larger than those used during training, whereas both classical GMRES and the learned optimizer fail to converge within comparable computational budgets. Moreover, the proposed method achieves arbitrary solution accuracies and exhibits strong out-of-distribution generalization across diverse source and velocity configurations.
These results demonstrate that appropriately designed mixed training datasets substantially improve the effectiveness of learned neural operator preconditioners. More broadly, this work highlights the importance of dataset design in scientific machine learning and provides a practical framework for integrating matrix-free neural operator preconditioning with Krylov subspace methods for solving practical large-scale Helmholtz problems.
\end{abstract}


\section{Introduction}

This work focuses on developing neural operator preconditioners trained on mixed datasets to accelerate the solution of sequences of linear systems arising from the discretization of the two-dimensional parametric Helmholtz equation in transcranial ultrasound therapy~\cite{skullDataset2018}.
In this setting, the Helmholtz equation incorporates a heterogeneous velocity field derived from human head CT data, resulting in complex wave interactions within the adult human skull.
Since skull geometries and acoustic properties vary significantly across patients, the corresponding Helmholtz operators can exhibit considerable variability, making the development of robust and efficient hybrid numerical solvers, potentially enhanced by machine learning, essential for practical clinical applications.
This challenging transcranial ultrasound problem was previously investigated by Stanziola et al.~\cite{Stanziola2021_HelmholtzEquationSolver}, who proposed a learned optimizer based on a recurrent neural network embedded within a Richardson iteration. For the representative transcranial ultrasound example considered in their study, the learned optimizer outperformed classical GMRES~\cite{sasc:86} in computational efficiency while exhibiting encouraging generalization capabilities. However, its attainable accuracy remained limited, and it failed to converge on the additional transcranial ultrasound examples considered in the present work (see Table~\ref{tab_results_test_dirac1_idea_skull} and Figure~\ref{fig:nn-fgmres_gen_sos} in Section~\ref{subsec:skull_exp}). 
Rather than replacing classical iterative solvers, our approach in~\cite{Giraud2025_NeuralNetworkPreconditioning} employs neural networks as nonlinear preconditioners within flexible GMRES (FGMRES)~\cite{saad:93}, combining the approximation and generalization capabilities of scientific machine learning~\cite{SciML25}, particularly neural operator learning~\cite{neural_operator}, with the robustness and reliability of Krylov subspace methods. The resulting hybrid methods converge reliably to arbitrary prescribed accuracies without requiring repeated network training and exhibit strong generalization across a broad range of Helmholtz problems.
Nevertheless, the previous framework was trained exclusively on randomly generated datasets, and its computational efficiency deteriorates when applied to the realistic challenging application-oriented transcranial ultrasound simulations. This observation motivates the present work. We investigate whether carefully designed mixed training datasets can improve the effectiveness of neural operator preconditioning while preserving its robustness and generalization capabilities.
Specifically, we train the same neural operator architecture using six different mixed datasets for the Helmholtz equation. These datasets combine application-oriented samples that mimic realistic transcranial ultrasound configurations with randomized source and velocity fields to maintain sufficient diversity during training. Their relative performance is systematically evaluated to identify dataset compositions that provide a favorable balance between specialization for the target application and generalization to previously unseen configurations.
The present work can therefore be viewed as a practical extension of our previous neural operator preconditioning framework~\cite{Giraud2025_NeuralNetworkPreconditioning}. Whereas our earlier study relied exclusively on randomly generated training data, the central contribution of the present work is the design and systematic investigation of mixed training datasets that improve efficiency for challenging transcranial ultrasound simulations while preserving strong out-of-distribution performance on other Helmholtz problems.
More broadly, the present work demonstrates that careful dataset design is a key ingredient for translating scientific machine learning methods from idealized benchmarks to challenging real-world applications.

In parallel, recent advances in scientific machine learning, particularly in neural operator learning and diffusion models, have also highlighted the importance of exposing models to diverse training data distributions to improve robustness and out-of-distribution generalization across varying physical parameters, geometries, and discretizations~\cite{li2021fourier, DeepONet, brandstetter2022message, kovachki2023neural,pathak2022fourcastnetglobaldatadrivenhighresolution, PSAROS2022meta_learning, XianqiLi2025_DataConsistency, XianqiLi2026_BrainTumor}. 
Some related developments in neural network preconditioning for classical numerical methods are discussed in~\cite{dimola2025numerical, Kopanikov2025, trifonov2025efficient, li2025neural, Azulay2023_MultigridAugmentedDeepLearning}.
The field of scientific machine learning is evolving rapidly, and the references cited above are intended to highlight representative developments closely related to the present work rather than to provide a comprehensive survey.

The remainder of this paper is organized as follows. Section~\ref{sec:backgroud} introduces the governing Helmholtz equation arising in transcranial ultrasound therapy together with its numerical solution using nonlinear preconditioned subspace method. We present the mathematical formulation, including the boundary conditions, discretization, and iterative solution strategies, and briefly review recent developments in neural network solvers and neural operator preconditioners for the Helmholtz equation, highlighting the limitations that motivate the present work.
Section~\ref{sec:nn-components} presents the proposed convolutional neural network (CNN)-based neural operator preconditioner and the construction of the six mixed training datasets. The resulting preconditioner is designed to accelerate FGMRES for solving sequences of heterogeneous Helmholtz equations arising in transcranial ultrasound simulations while maintaining strong generalization to other Helmholtz applications.
Section~\ref{sec:experiments} evaluates the proposed approach on a variety of transcranial ultrasound examples with previously unseen physical parameters. We systematically compare the different mixed training datasets and demonstrate that the resulting preconditioners substantially accelerate FGMRES while exhibiting strong out-of-distribution generalization across a broad range of Helmholtz problems.
Finally, Section~\ref{sec:conclusion} summarizes the main findings and discusses their implications for scientific machine learning, neural operator preconditioning, and large-scale PDE simulations.

Key notations used throughout this paper are summarized below. Unless otherwise specified, $\|\cdot\|$ denotes the Euclidean norm for both vectors and matrices. The superscript~$^H$ denotes the conjugate transpose.
The symbols~$\R$ and~$\C$ denote the sets of real and complex numbers, respectively.
For convenience, some \texttt{Python}-style notation is adopted in the presentation of the algorithm.
To improve readability, we use the following notation throughout the paper.
Scalars are denoted by lowercase letters. Bold lowercase letters, e.g., $\bx$, denote vectors (or columns of a matrix), while uppercase letters, e.g., $A$, denote matrices. Functions, operators, and mappings are represented by calligraphic letters, e.g., $\Fc.$
Unless otherwise stated, the subscript~$_j$ denotes the quantity obtained at iteration~$j$, while the positive integer subscript~$_m$ denotes the maximum number of iterations within a Krylov cycle.
The inner product between two vectors~$\bx$ and~$\by$ is denoted by $\langle \bx, \by\rangle$, and the linear span of the vectors~$\bx_1,\ldots,\bx_j$ is denoted by $\operatorname{span}\{\bx_1,\ldots,\bx_j\}.$

%%*************************************************
\section{Background}~\label{sec:backgroud}

%%*************************************************
\subsection{Governing equation: The Helmholtz equation}\label{sec:helm}

\begin{figure}[!ht]\centering%
\resizebox{.4\textwidth}{!}{\begin{tikzpicture}[x=8pt,y=8pt]%
\node[rectangle, draw=black, line width=1.5pt, fill=gray!20, fit={(0,0)(20,20)}, inner sep=0pt] {};
\node[rectangle, draw=black, line width=1pt, fill=white, fit={(2,2)(18,18)}, inner sep=0pt] {$\Oinner$};
\node at (10, 19) {$\Opml$};
\draw[dashed, line width=.5pt] (-4,20) -- (2,20);
\draw[dashed, line width=.5pt] (-2,18) -- (2,18);
\draw[dashed, line width=.5pt] (-4,0) -- (2,0);
\draw[<->, line width=1pt] (-1,18) -- node[midway, left]{$\ell$} (-1,20) ;
\draw[<->, line width=1pt] (-3,0) -- node[midway, left]{$2L$} (-3,20); 
\end{tikzpicture}}
\caption{Illustration of the PML boundary condition on two-dimensional domain~(\cite[Figure~1]{Giraud2025_NeuralNetworkPreconditioning}).}
\label{fig:pml}
\end{figure}
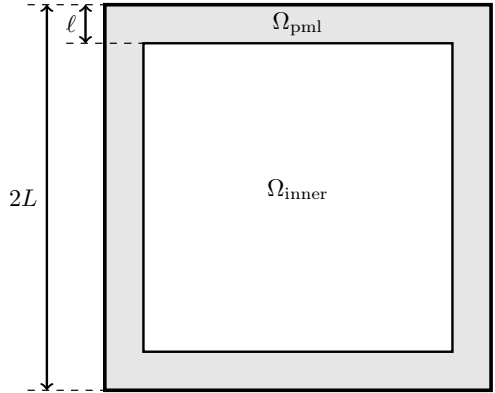

We consider the two-dimensional Helmholtz equation, subject to the so-called Sommerfeld radiation condition,
\begin{equation}\label{eq:helmholtz}
    \begin{cases}
        \nabla^2 u + k^2 u = f \; \mbox{ in } \R^2,\\[\medskipamount]
        \lim_{\|\vx\|_2 \to \infty} \|\vx\|_2^{1/2} 
        \left(
            \dfrac{\partial u}{\partial \|\vx\|_2}(\vx) - \jmath 
            k(\vx)
            u(\vx)
        \right)
        = 0,
    \end{cases}
\end{equation}
where $\jmath \in \C$ denotes the imaginary unit such that $\jmath^2=-1$, and $u \colon \R^2 \to \C$, $f \colon \R^2 \to \C$, $k \colon \R^2 \to \R^+$ are scalar fields, referred to as the solution field, the source field, and the wavenumber field, respectively.
The latter is defined by $k \colon \vx \mapsto \omega / c(\vx)$, where $\omega \in \R^+$ is the angular frequency of the source, and $c \colon \R^2 \to \R^+$ is the velocity (speed of sound) field.
In what follows, we restrict ourselves to the  case $\omega = 1$.
Numerically, the unbounded domain $\R^2$ is truncated to a computational, square domain $\Omega \bydef [-\L, \L]^2 \subset \R^2$ and
the derivatives in Equation~\eqref{eq:helmholtz} are discretized using Fourier differentiation~\cite{Stanziola2021_HelmholtzEquationSolver, Xiang2022_SolutionLargeLinear, Shpakovych2023_NeuralNetworkPreconditioning} with $N$ points in each direction.
Consequently, there is a total of $n \bydef N^2$ degrees of freedom.
The boundary conditions are (approximately) enforced using perfectly matched layers (PMLs)~\cite{Berenger1994_PerfectlyMatchedLayer} in the outer region $\Opml \bydef \Omega \setminus \Oinner$, where $\Oinner \bydef [-\Linner, \Linner]^2$, with $\Linner \bydef \L - \Lpml$ and $\Lpml < \L$, as illustrated in Figure~\ref{fig:pml}.

Within the PML framework, the absorbing layer for truncating the domain can be characterized by a damping function
\begin{equation}\label{eq:pml}
    \sigma(\bx) = 
    \begin{cases}
        0, \;  \mbox{ for } \bx \in \Oinner \subset \Omega,\\[\medskipamount]
        \sigma_{\text{max}} \dfrac{(\bx - \Linner)^2}{l^2}, \;  \mbox{ for } \bx \in \Opml \subset \Omega,
    \end{cases}
\end{equation}
where $\sigma_{\text{max}} \in \R^+$ is a prescribed parameter denoting the maximum damping coefficient.
With this construction, the solution remains unchanged in the inner domain $\Oinner \subset \Omega,$ while it decays exponentially in the PML region $\Opml \subset \Omega.$ As a result, outgoing waves are effectively absorbed without producing spurious reflections at the interface between the physical domain and the absorbing layer.

The discretization of the Helmholtz problem~\eqref{eq:helmholtz} with varying physical parameters (i.e., varying source filed $f$ and velocity filed $c$) yields a groups of linear systems with varying left and right hand sides 
\begin{equation}~\label{eq:discrePDEs}
    A(\bc)^{(\ell)} \bu = \bb^{(\ell)}, \text{ with } A :=
	\nabla^2 + k^2 \in \C^{n\times n}, \ \ell =1,2,....\K,
\end{equation}
where, associated with the $\ell$-{th} family from each discrete Helmholtz equation
involving the discrete velocity field $\bc \in \R^{n}$ and the associated discretized operator (matrix) $A(\bc) \in \C^{n\times n}$, as well as the discretized source term $\bb \in \C^{n}$ and solution $\bu \in \C^{n}.$

%******************************************************
\subsection{Subspace solver with nonlinear preconditioner}\label{sec:nla}
The solution of large linear systems of equations is commonly addressed using iterative methods.
These approaches are attractive due to their modest memory requirements and their ability to stop the iterations when the quality of the solution is similar to the possible uncertainty in the matrix entries or right-hand sides, which may arise from discretization errors in PDE or from uncertainties in the input data. In this context, backward error analysis~\cite{wilkinson1963rounding, Rigal1967} provides powerful techniques for designing meaningful stopping criteria.
Among iterative techniques, subspace methods, such as Krylov subspace methods, have become standard tools over the past decades. For unsymmetric problems, the best known method is GMRES~\cite{sasc:86}, which minimizes the residual norm over an increasing Krylov subspace. This technique relies on the so-called Arnoldi process, which incrementally builds an orthonormal basis of the nested subspaces. At iteration $k$, the Arnoldi algorithm computes a set of orthonormal vectors $V_k = [\bv_1, \dots, \bv_k ]$ that satisfy the so-called Arnoldi relation, which is written in matrix form:
$$
A V_k = V_{k+1}\bar{H}_k \quad \mbox{\rm with } V_{k+1}^H V_{k+1} = I_{k+1}
$$
where $\bar{H}_k \in \mathbb{C}^{(k+1)\times k}$ is upper Hessenberg.
Assume the approximated solution at iteration~$k$ is denoted as $\bx_k,$ the minimum residual norm iterate $\bx_k \in \Span\{\bv_1, \dots, \bv_k \}$, that is $\bx_k = V_k \by_k$ with  $\by_k =\argmin_{ \by \in \C^k} (\| \bar{H}_k \by - \beta \be_1 \| )$ where $\beta =\| \bb \|$.
In practice, a preconditioner $M$ is used to speed up the convergence of GMRES, where $M$ is expected to approximate $A^{-1}$ somehow. For GMRES it is recommended to use the preconditioner on the right, i.e., GMRES solves $AM \bt = \bb$, with $M \bt=\bx$. Thus, at each iteration, GMRES still minimizes the residual norm of the residual associated with the original linear system. The Arnoldi relation is
$$
A M V_k = V_{k+1}\bar{H}_k \quad \mbox{\rm with } V_{k+1}^H V_{k+1} = I_{k+1}.
$$
In~\cite{saad:93}, Saad introduced the idea of a flexible preconditioner, which allows to have a different preconditioning matrix $M_i$ at each iteration, so that the generalized Arnoldi relation becomes
\begin{equation}\label{eq:generalizedArnoldi}
    A Z_k = V_{k+1}\bar{H}_k \quad \mbox{\rm with } V_{k+1}^H V_{k+1} = I_{k+1}, 
\end{equation}
where
% $Z_k = [M_1 v_1, \dots, M_k v_k]$.
$Z_k = [\bz_1, \dots, \bz_k]$ with $\bz_i = M_i \bv_i~(i=1,\dots,k)$.
This idea can easily be extended to the situation where the preconditioner is no longer linear, that is, $\bz_i = {\cal M}_i(\bv_i)$, where each ${\cal M}_i$ may be a nonlinear operator. 
In the latter case, Equation~\eqref{eq:generalizedArnoldi} still holds and uniquely defines each iterate as long as $Z_k$ remains full rank.
In this work, we adopt ${\cal M}_i() = \NN()$, where $\NN$ denotes the trained neural network inference (see Section~\ref{sec:nn-components}). The nonlinearity has two origins; first, the activation functions in the neurons are nonlinear, and second, the preconditioner is computed in 32-bit arithmetic while the rest of the computation is done in 64-bit arithmetic. This truncation due to casting 64 bits to 32 bits is also nonlinear.
A sketch of the CNN-preconditioned FGMRES algorithm is provided in Algorithm~\ref{alg:nn-fgmres}.

{\DontPrintSemicolon\LinesNumbered%
\begin{algorithm2e}[ht]%
  $\br_0 = \bb - A \bx_0$, $\beta = \|\br_0\|$, $\bv_1 = \br_0 / \beta$ \;
  \For{$j=1,\ldots,m$}{
    $\bz_j = \NN(\bv_j, \text{**kwargs})$ {(Here, **kwargs means ``keyword arguments" in Python)} \; % **kwargs means ``keyword arguments" in Python is a way to pass a variable number of named arguments to a function
    $\bw = A \bz_j$ \; 
    \For{$i=1,\ldots,j$}{
      $h_{i,j} = \langle \bw , \bv_i \rangle$\;
      $\bw \leftarrow \bw - h_{i,j} \bv_i$ \;
    }
    $h_{j+1, j} = \|\bw\|$ \;
    $\bv_{j+1} = \bw / h_{j+1, j}$ \;
  }
  $\bar{H}_m = [h_{i,j}]_{1\leq j \leq m, 1 \leq i \leq j+1}$ \;
  $Z_m = [\bz_1, \dots, \bz_m]$ \;
  $\by_m = \argmin_{\by \in \C^m} \| \bar{H}_m \by - \beta \be_1 \|$ \;
  $\bx_m = \bx_0 + Z_m \by_m$ \;
  \caption{Sketch of the CNN-preconditioned FGMRES algorithm~\cite[Algorithm~1]{Giraud2025_NeuralNetworkPreconditioning}.} %
  \label{alg:nn-fgmres}
\end{algorithm2e}%
}

%******************************************************
\subsection{Related work on neural network solvers and neural operator preconditioners for the Helmholtz equation}\label{sec:pre_nns_pre}

The present work is motivated by a sequence of previous studies on machine learning approaches for solving Helmholtz equations arising in transcranial ultrasound simulations.
An early contribution in this direction is the learned optimizer proposed by Stanziola et al.~\cite{Stanziola2021_HelmholtzEquationSolver}. Their approach employed a modified U-Net architecture~\cite{Ronneberger2015_UNetConvolutionalNetworks} within a Richardson iteration framework, in which a recurrent neural network generated a corrective update at each iteration. The model was trained using a fixed source configuration and velocity fields derived from idealized skull geometries.
While this work demonstrated the potential of neural-network-based solvers for Helmholtz problems, it also revealed several important limitations. In particular, the attainable accuracy was limited (approximately $10^{-4}$), and the method successfully solved only one practical transcranial ultrasound example while failing on the remaining cases considered in the present work (see Table~\ref{tab_results_test_dirac1_idea_skull} and Figure~\ref{fig:nn-fgmres_gen_sos} in Section~\ref{subsec:skull_exp}). 
Moreover, the training procedure exhibited instabilities, including producing divergent training losses and sensitivity to initialization (see~\cite[Figure~5.2]{Xiang2022_SolutionLargeLinear}), which occasionally required training attempts to obtain a satisfactory model.
To address these issues, an optimal step-size strategy was introduced in~\cite[Chapter~5]{Xiang2022_SolutionLargeLinear}, leading to a minimum-residual Richardson iteration. By explicitly minimizing the residual at each iteration, this formulation substantially improved both the attainable accuracy and the stability of the training procedure. Beyond using neural networks as correction operators, this work also introduced the idea of employing them as nonlinear preconditioners within two types of Krylov subspace methods.
Subsequent studies~\cite{Shpakovych2023_NeuralNetworkPreconditioning,Giraud2025_NeuralNetworkPreconditioning} extended this neural operator preconditioning framework to more general Helmholtz problems and other classes of PDEs. 
Unlike purely data-driven solvers, these hybrid approaches combine the approximation capability and generalization of neural operators with the robustness of Krylov subspace methods, enabling convergence to arbitrary prescribed tolerances while retaining strong performance across a broad range of benchmark problems.
Although these methods demonstrated excellent robustness and generalization across diverse Helmholtz problems, their computational efficiency remained limited for challenging application-specific transcranial ultrasound simulations. The objective of the present work is therefore to investigate whether application-oriented mixed training datasets can improve the efficiency of neural operator preconditioning while preserving its robustness and out-of-distribution generalization.

%******************************************************
\section{Components of nonlinear neural operator preconditioning}\label{sec:nn-components}

Assume the CNN-based preconditioning operator is denoted as $\NN$ with trainable parameters~${\theta}$~(i.e., the trainable weights and biases of the CNNs). We focus on learning an operator $\NN$
\begin{equation}~\label{eq:operator}
	\NN: \Pc \longrightarrow \Uc, 
\end{equation}
with coefficient input $\left\{\bb^{(\ell)}, \text{**kwarg}\right\} \in \Pc, \text{ **kwarg} = \left\{\bc^{(\ell)}, \sigma \right\}, \ \ell =1,\ldots, \K$ (here $\K$ denotes the size of the training dataset, similar to the number of the linear systems shown in Equation~\eqref{eq:discrePDEs}) and output solution $\bu^{(\ell)}_{\theta} \in \Uc$ returned from CNNs with $\theta.$
Our target is to find the learning operator such that 
\begin{equation}~\label{eq:operator_approximation}
	\NN(\bb^{(\ell)}, \text{**kwarg}) \sim \bu^{(\ell)}_{\theta}, \mbox{ with } \text{**kwarg} = \left\{\bc^{(\ell)}, \sigma \right\}, \ i =1,\ldots, \K,
\end{equation}
which means the learned operator $\NN$ approximates to the inverse of the coefficient matrix $A^{-1}$ since we have $\bu = A^{-1} (\bb, \text{**kwarg})$ from Equation~\eqref{eq:discrePDEs}.
Because of this, the learned operator $\NN$ can be used as a \text{nonlinear} \text{CNN-based} preconditioner in the FGMRES method described in Algorithm~\ref{alg:nn-fgmres}, where the preconditioned Krylov basis is computed as: 
\begin{equation}~\label{eq:nns_pre}
    \bz_j = \NN(\bv_j, \bc, \sigma)
\end{equation}
for the application in the parametric Helmholtz equation~\eqref{eq:helmholtz}. Except for use the Krylov basis as the input of the trained inference, refer to~\cite[Chapter~5]{Xiang2022_SolutionLargeLinear} for more details about other possibilities with similar numerical performance.

To distinguish from the original FGMRES~\cite{saad:93}, we denote FGMRES with the \text{CNN-based} preconditioning step shown in Equation~\eqref{eq:nns_pre} for Algorithm~\ref{alg:nn-fgmres} as the \FGMRES algorithm.
Compared to the original FGMRES, these \text{CNN-based preconditioning} steps are performed in 32-bit, consistent with the precision used during the training processes. The other remaining operations in the \FGMRES algorithm are conducted in 64-bit. 
Thus, this \FGMRES algorithm is realized under mixed-precision calculations.

\subsection{Network architecture}\label{sec:network-architecture}

As introduced in Section~\ref{sec:nla}, our approach uses a neural network as a nonlinear preconditioner.
To do so, we design a convolutional neural network based on the U-Net architecture~\cite{Ronneberger2015_UNetConvolutionalNetworks}, which takes the discretized source field, velocity field and the PML damping function $\sigma$ as inputs.
To make them suitable to the convolutional nature of the U-Net, these 2D fields are recast as image-like tensors.
Specifically, recalling that $n=N^2$, the complex-valued source term $\bb \in \C^n$ may be recast as a $2 \times N \times N$ tensor composed of two $N \times N$ channels representing the real and imaginary parts of the field.
So does the complex $\sigma \in \C^n$ function.
Similarly, the real-valued velocity field $\bc \in \R^n$ may be recast as a $1 \times N \times N$ tensor.
These five inputs are concatenated along channel dimension to get a $5 \times N \times N$ input tensor.
The output of the network is a $2 \times N \times N$ tensor, representing a 2D, complex-valued field.
As we shall see in Section~\ref{sec:training}, the network will be trained to provide approximations (predictions) $\hat{\bu}_\theta$ of discrete solution fields $\bu$ of the linear system $A(\bc) \bu = \bb$ corresponding to a discretized version of the Helmholtz problem \eqref{eq:helmholtz} with seetings described in Table~\ref{tab:helm_model_parameters}.

\begin{figure}[!ht]\centering%
\includegraphics[width=0.80\textwidth]{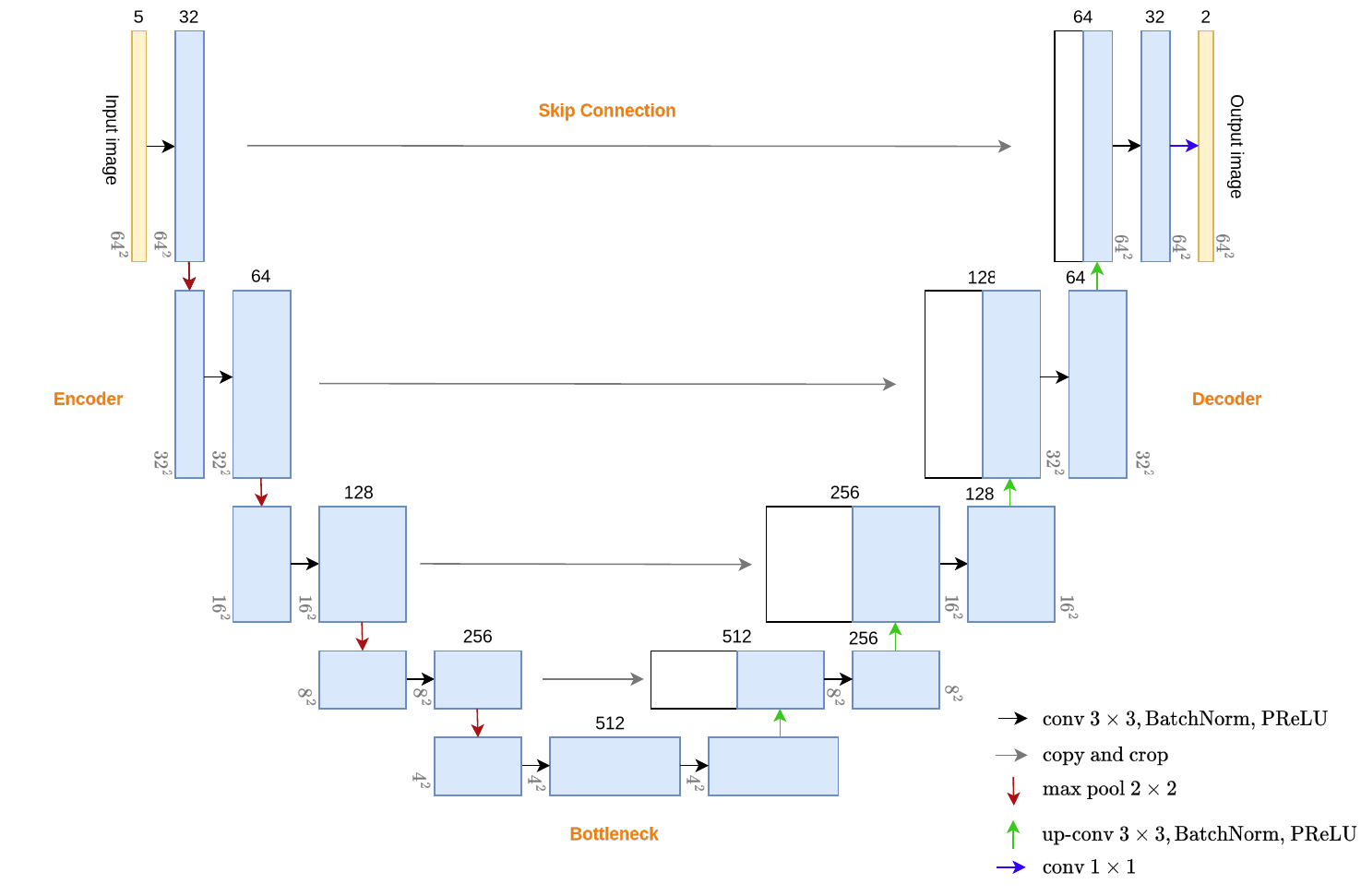}
\caption{Overview of our U-Net architecture with 4 depth.}
\label{fig:Unet}
\end{figure}

The primary hyperparameter of the U-Net architecture is its depth. For the Helmholtz problems considered in Table~\ref{tab:helm_model_parameters}, a depth of four is sufficient to construct an effective neural preconditioner.
The resulting U-Net architecture is illustrated in Figure~\ref{fig:Unet}. It consists of an encoder–decoder structure, where the encoder progressively downsamples the input and the decoder reconstructs the output through successive upsampling operations. The encoder and decoder are linked through a bottleneck layer and a set of skip connections. These skip connections play a crucial role in preserving fine-scale information while combining it with the high-level features extracted by the encoder, thereby enabling the network to capture both global context and local details.
The corresponding hyperparameter settings are summarized in Table~\ref{tab:model_hyperparameters}.
With the Helmholtz problem parameters and network hyperparameters specified in Tables~\ref{tab:helm_model_parameters} and~\ref{tab:model_hyperparameters}, our U-Net architecture contains approximately \num{832}~K trainable parameters, which is independent of $n$ owing to its convolutional nature.

%******************************************************
%******************************************************

%%%**************************************************************
\begin{table}[!h]
	\centering
	\begin{minipage}{0.99\textwidth}
		\centering
		\caption{Helmholtz equation with PML in 2D and Fast Fourier Transforms (FFT) parameters}
		\label{tab:helm_model_parameters}
		\begin{tabular}{lr}
			\toprule
			\textbf{Parameter} & \textbf{Value} \\
			\midrule
			Limits of the Domain & $[[-20, 20], [-20, 20]]$ \\
			Discretisation points $n$ in Domain & $[64, 64]$ \\
			Axes for which FFT is applied & $[-2, -1] $ \\
            Configurations for PML & $l=4, \ \sigma_{\text{max}} = 2$ \\
			%Percentage Size of PML within the Domain & 10\% \\
			\bottomrule
		\end{tabular}
	\end{minipage} \\
%	\end{table}
	\vskip .5cm
%\begin{table}[h]
	\hskip .10cm
	\begin{minipage}{0.45\textwidth}
		\centering
		\caption{U-Net Model Hyper-parameters}
		\label{tab:model_hyperparameters}
		\begin{tabular}{lr}
			\toprule
			\textbf{Parameter} & \textbf{Value} \\
			\midrule
			Input Channels & 5 \\
			Output Channels & 2 \\
			Activation Function & PReLU \\
			Depth & 4 \\
			Channels per Layer & 32 \\
			Channels per State & 32 \\
			Convolving kernel size, Stride, Padding & 8, 2, 3 \\
			Batch normalization & True \\
			Bias in concolutional layers & True \\
			\bottomrule
		\end{tabular}
	\end{minipage} 
%\end{table}
\hskip 1.2cm
%\begin{table}[!h]
	\begin{minipage}{0.45\textwidth}
		\centering
		\caption{Training Hyper-parameters}
		\label{tab:training_hyperparameters}
		\begin{tabular}{lr}
			\toprule
			\textbf{Parameter} & \textbf{Value} \\
			\midrule
            Maximal Epoch & 500 \\
			Optimizer & Adam \\
			Batch Size & 32 \\
			%	Buffer Size & 100 \\
			Gradient Clipping value & 1 \\
			Learning Rate (LR) & $10^{-3}$ \\
			Minimal LR & $10^{-5}$ \\
			Seed for Random Number & 42 \\
			\bottomrule
		\end{tabular}
	\end{minipage}
\end{table}
%******************************************************

The abstract, functional representation of the neural network preconditioning with trainable parameters~$\theta$ is thus a nonlinear operator,
\begin{equation}
  \label{eq:nn_operator}
  \begin{aligned}
    \NN \colon \C^n \times \R^n \times \C^n & \to \C^n,\\
    (\bb, \bc, \sigma) & \mapsto \NN(\bb, \bc, \sigma) = \NN(\bc, \sigma)(\bb),
  \end{aligned}
\end{equation}
where the conversions to and from image-like tensors are implicit.

\subsection{Mixed training dataset}\label{sec:dataset-training}
\begin{figure}[!ht]%
    \centering%
    \includegraphics[width=\textwidth]{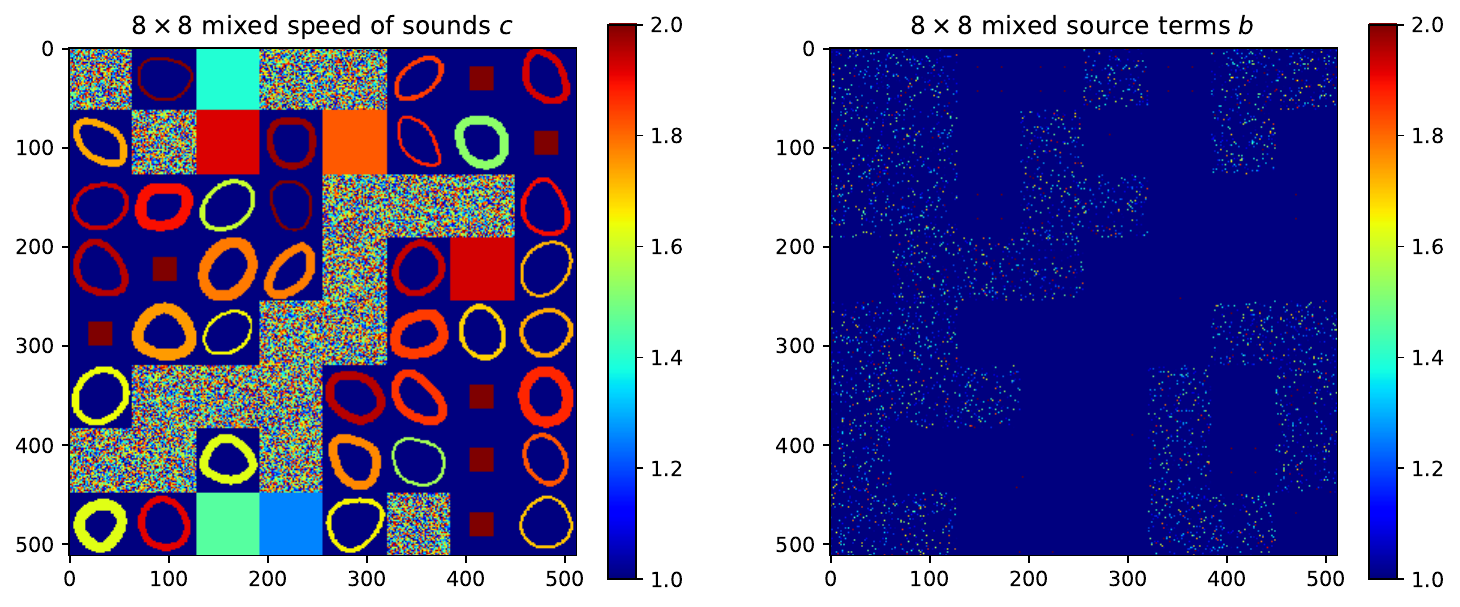}%
    \caption{Visualize some samples from the different types of the training dataset $(\bc, \bb)$ illustrated at blocks~\ref{list_data_c}-\ref{list_data_b} (colorbar of the mixed source terms $\bb$ shown in the right plot has been scaled for better visualisation).}
    \label{fig:visu_train_dataset}
\end{figure}

In this work, mixed training datasets are employed to train the neural operator preconditioners. Each training dataset is composed of pairs of velocity fields~$\bc$ and source fields~$\bb$, denoted by $(\bc,\bb)$. By defining the distributions of~$\bc$ and~$\bb$ independently, different combinations of structured, application-oriented samples and randomized samples can be constructed.
We first define a fully mixed dataset by considering the following distributions for~$\bc$ and~$\bb$:
\begin{enumerate}
    \item For the velocity field $\bc \in \R^n:$ $\bc \in \left\{{\cal U}([1,2]), \ \text{\texttt{TruncIdealSkulls}}, \ \text{\texttt{randConstant}}, \ \text{\texttt{square}} \right\},$ the four configurations are sampled with probabilities $[0.4, 0.4, 0.1, 0.1],$ respectively:~\label{list_data_c} 
            \begin{enumerate}
            \item $\bc \simiid {\cal U}([1, 2])$: each entry of $\bc$ is independently sampled from the uniform distribution on $[1, 2]$;
            \item $\bc \simiid \text{\texttt{TruncIdealSkulls}}$: $\bc$ is generated from the centrally truncated regions of the training dataset of idealized skull models introduced in~\cite{Stanziola2021_HelmholtzEquationSolver}. The truncated region preserves the main geometric features of the idealized skull structures;
            \item $\bc \simiid \text{\texttt{randConstant}}$: all entries of $\bc$ are assigned the same value, sampled uniformly at random from $[1, 2]$;
            \item $\bc \simiid \text{\texttt{square}}$: a square region centered in the computational domain is assigned the value 2, while all remaining entries are set to 1.
        \end{enumerate}
    \item For the source field $\bb \in \C^n:$ $\bb \in \left\{{\cal N}(0,1), \ \text{\texttt{dirac}}, \ \text{\texttt{dirac1}}, \ \text{\texttt{dirac4}} \right\},$ the four configurations are sampled with probabilities $[0.6, 0.3, 0.05, 0.05],$ respectively:~\label{list_data_b}
        \begin{enumerate}
            \item $\bb \simiid {\cal N}(0,1)$: each entry of $\bb$ is independently sampled from a standard normal distribution with zero mean and unit variance;
            \item $\bb \simiid \text{\texttt{dirac}}$: a single entry of $\bb$ is assigned the value 60 at a randomly selected grid location, while all remaining entries are set to zero;~\label{list_data_b_b}
            \item $\bb \simiid \text{\texttt{dirac1}}$: a single entry of $\bb$ is assigned the value 60 at the center of the computational domain, while all remaining entries are set to zero;
            \item $\bb \simiid \text{\texttt{dirac4}}$: four entries of $\bb$ are assigned the value $30$ at four fixed grid locations within the computational domain, while all remaining entries are set to zero.~\label{list_data_b_d}
        \end{enumerate}
\end{enumerate}
Refer to Figure~\ref{fig:visu_train_dataset} for a visualisation of representative samples from the different types of training data pairs~$(\bc, \bb)$ described in blocks~\ref{list_data_c}-\ref{list_data_b}.
We note that the amplitudes 60 or 30 was applied to the \texttt{dirac}-type source configurations 
so that the resulting source field has a magnitude comparable to that of the standard normal setting ${\cal N}(0,1)$ over the computational domain. This is important in the training process, which normalization helps maintain the effectiveness of the neural operator preconditioning across different source distributions.
Such localized \texttt{dirac}-type sources are frequently encountered in practical wave propagation problems, including transcranial ultrasound simulations and spherical wave propagation, where the source term is concentrated at a small number of points or within a localized region and vanishes elsewhere.

According to the velocity and source distributions defined in blocks~\ref{list_data_c}-\ref{list_data_b}, five additional mixed datasets are then obtained by selecting different subsets or combinations of these distributions. The final six training datasets considered in this work are summarized in Table~\ref{tab:cases_mixed_dataset}.
Each dataset is denoted by \texttt{D0}-\texttt{D5}, where \texttt{D2} corresponds to the mixed setting described in blocks~\ref{list_data_c}-\ref{list_data_b}.  
As shown in Table~\ref{tab:cases_mixed_dataset}, datasets \texttt{D0} and \texttt{D1} share the same configuration of $\bc$ but differ in the source term $\bb.$ Similarly, \texttt{D2} and \texttt{D3} share the same configuration of $\bc$ while differing in $\bb.$ 
In contrast, \texttt{D4} and \texttt{D5} use the same configuration of $\bb$ but vary in $\bc.$
In the following, we describe how neural networks are trained using these six mixed datasets to construct effective neural operator preconditioners, with a primary focus on applications to transcranial ultrasound simulation.

\begin{table}[!ht]
    \centering\small
    \begin{tabular}{lll}
        \toprule    
        \bfseries \#Dataset & \bfseries \begin{tabular}{c} Source term $\bb$ \\ (with probability) \end{tabular} & \bfseries  \begin{tabular}{c} Velocity field (or Speed of sound) $\bc$ \\ (with probability) \end{tabular}  \\
        \midrule % \texttt{D5}
        \texttt{D0} & \begin{tabular}{c} $\left\{{\cal N}(0,1), \ \text{\texttt{dirac}}, \ \text{\texttt{dirac1}}, \ \text{\texttt{dirac4}} \right\}$ \\ ($[0.6, 0.3, 0.05, 0.05]$) \end{tabular} & \begin{tabular}{c} $\left\{{\cal U}([1,2]), \ \text{\texttt{TruncIdealSkulls}} \right\}$ \\ ($[0.5, 0.5]$) \end{tabular} \\
        \midrule % \texttt{D2}
        \texttt{D1} & \begin{tabular}{c} ${\cal N}(0,1)$ \end{tabular} & \begin{tabular}{c} $\left\{{\cal U}([1,2]), \ \text{\texttt{TruncIdealSkulls}} \right\}$ \\ ($[0.5, 0.5]$) \end{tabular} \\
        \midrule % \texttt{D0}
        \texttt{D2} & \begin{tabular}{c} $\left\{{\cal N}(0,1), \ \text{\texttt{dirac}}, \ \text{\texttt{dirac1}}, \ \text{\texttt{dirac4}} \right\}$ \\ ($[0.6, 0.3, 0.05, 0.05]$) \end{tabular} & \begin{tabular}{c} $\left\{{\cal U}([1,2]), \ \text{\texttt{TruncIdealSkulls}}, \ \text{\texttt{randConstant}}, \ \text{\texttt{square}} \right\}$ \\ ($[0.4, 0.4, 0.1, 0.1]$) \end{tabular} \\
        \midrule % \texttt{D1}
        \texttt{D3} & \begin{tabular}{c} ${\cal N}(0,1)$ \end{tabular} & \begin{tabular}{c} $\left\{{\cal U}([1,2]), \ \text{\texttt{TruncIdealSkulls}}, \ \text{\texttt{randConstant}}, \ \text{\texttt{square}} \right\}$ \\ ($[0.4, 0.4, 0.1, 0.1]$) \end{tabular} \\
        \midrule % \texttt{D3}
        \texttt{D4} & \begin{tabular}{c} ${\cal N}(0,1)$ \end{tabular} & \begin{tabular}{c} \text{\texttt{TruncIdealSkulls}} \end{tabular} \\
        \midrule % \texttt{D4}
        \texttt{D5} & \begin{tabular}{c} ${\cal N}(0,1)$ \end{tabular} & \begin{tabular}{c} ${\cal U}([1,2])$ \end{tabular} \\
        \bottomrule
    \end{tabular}
    \caption{Six types of mixed training dataset pairs $(\bb, \bc)$ used for training neural operator preconditioning.} % Default mixed training dataset \texttt{D0}, and other five cases of the training dataset.
    \label{tab:cases_mixed_dataset}
\end{table}

\subsection{Training neural operator preconditioner}\label{sec:training}

\begin{figure}[!ht]\centering%
\includegraphics[width=0.80\textwidth]{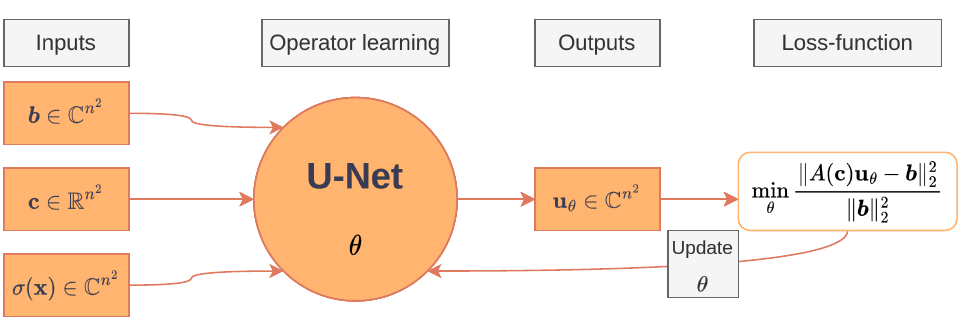}
\caption{Schematic representation of the training framework.}
\label{fig:nn_general_schamatic}
\end{figure}

The training process aims at finding the parameters $\theta^* \in \Theta$ of the neural network that minimize a certain cost functional $\cost \colon \Theta \to \R$, where $\Theta$ denotes the parameter space, that is,
\begin{equation}
  \label{eq:nn_optim}
  \theta^* \in \argmin_{\theta \in \Theta} \cost(\theta).
\end{equation}
For our specific problem, the cost function takes the form
\begin{equation}
  \label{eq:nn_cost_fn}
  \cost(\theta) = \mathbb{E}_{\bb, \bc, \sigma} [\loss(\bb, \bc, \sigma)],
\end{equation}
where $\bb,$ $\bc$ and $\sigma$ are the inputs of the network, $\loss$ is a loss function and $\mathbb{E}_{\bb, \bc, \sigma}$ denotes the expectation operator with respect to $\bb,$ $\bc$ and $\sigma.$
In this work, we consider the loss function defined by
\begin{equation}
  \label{eq:nn_loss_fn}
  \loss(\bb, \bc, \sigma)
  =
  \dfrac{\|\bb - A(\bc)\NN(\bb,\bc, \sigma)\|^2}{\|\bb\|^2}
  =
  \dfrac{\|\bb - A(\bc)\hat{\bu}_\theta \|^2}{\|\bb\|^2}
  =
  \dfrac{\|\bb - \hat{\bb}_\theta \|^2}{\|\bb\|^2},
\end{equation}
where $\hat{\bu}_\theta \bydef \NN(\bb, \bc, \sigma) = \NN(\bc, \sigma)(\bb)$ denotes the prediction of the neural network, and $\hat{\bb}_\theta = A(\bc)\hat{\bu}_\theta$.
The neural network is thus trained to be a good approximation of $A^{-1}$, i.e., roughly speaking, so that $\NN(\bc, \sigma) \approx A^{-1}(\bc)$ for any $\bc$ with fixed $\sigma$.
The training process can here be interpreted as unsupervised, in the sense that the prediction $\hat{\bu}_\theta$ are not (directly) compared to reference values $\bu^*$ in the loss function.
We remark that the residual of the linear system for a given prediction $\hat{\bu}_\theta = \NN(\bc, \sigma)(\bb)$, namely $\bb - A(\bc)\hat{\bu}_\theta$, which is at the core of the loss definition, corresponds to an algebraic (discretized) version of the residual of the PDE, including the boundary conditions.
As such, the loss may be interpreted as being physics-informed.
In terms of numerical linear algebra, we further note that for any given $\bc \in \R^n$, $\loss(\bb,\bc, \sigma) = \eta_b(\hat{\bu}_\theta)^2$, where $\eta_b(\tilde{\bu})$ denotes the backward error with respect to $\bb$ of $\tilde{\bu}$ as an approximate solution to the linear system $A(\bc) \bu = \bb$~\cite{higham2002accuracy, Rigal1967}.
An alternative interpretation that can be drawn from the last representation of the loss function \eqref{eq:nn_loss_fn} is that the neural network is trained so that $A(\bc)\NN(\bc, \sigma)$ be a good approximation of the identity operator, similar to an autoencoder~\cite{Bishop2024_Autoencoders}.
Thus, the neural network $\NN$ can then be interpreted as an encoder, trained such that it is decoded by $A$.
The training process is summarized in Figure~\ref{fig:nn_general_schamatic}.

%******************************************************
%******************************************************

\subsection{Some training details and remarks}
\label{sec:training_remarks}
%%***********************************************************
To train the CNN-based neural operator preconditioning described in Equation~\eqref{eq:nn_operator} or Figure~\ref{fig:nn_general_schamatic}, 10,000 random samples are generated for each training dataset listed in Table~\ref{tab:cases_mixed_dataset}, including distributions for $b$ and~$c.$
The training is conducted using the Adam optimizer, %~\cite{Kingma2017_AdamMethodStochastic}, 
with the maximum number of epochs~(denoted as \texttt{max\_epoch}) set to 500, that is $\text{\texttt{max\_epoch}}=500$. The learning rate is initialized at 0.001 and gradually decreases to a minimum value of $10^{-5}.$
Refer to Table~\ref{tab:helm_model_parameters}-\ref{tab:training_hyperparameters} 
for more details about the Helmholtz operator parameters, the model, and training hyper-parameters of the U-Net model designed for the Helmholtz operators on a two-dimensional domain.
This U-Net preconditioning model, configured with a $\text{\texttt{depth}}$ of 4, is trained under the single precision~(that is \texttt{float32} (32-bit) and \texttt{complex64}) on 
4~NVIDIA~V100~GPUs (32.5~GB) with \text{2x~16-core~Intel~Skylake~CPUs} (Model name: Intel(R) Xeon(R) Gold 6130 CPU @ 2.10GHz) located on the CAIUS cluster\footnote{le Cluster de cAlcul Intensif à l'Université de Strasbourg (CAIUS): \url{https://hpc.pages.unistra.fr}} of the University of Strasbourg Computing Centre (CCUS). 
The network architecture and the Helmholtz operator training utilities are implemented using the \texttt{PyTorch} and \texttt{PyTorch Lightning} libraries.
The training time required for this U-Net preconditioning model under these different training datasets for the discrete Helmholtz equations is summarized in Table~\ref{tab:train_time}.

%%%**************************************************************
\begin{table}[!h]
	\centering
	\begin{minipage}{.8\textwidth}
		\centering
		\caption{Training time for learning the Helmholtz operator with the six different datasets \texttt{D0}-\texttt{D5} described in Table~\ref{tab:cases_mixed_dataset}}
		\label{tab:train_time}
		\begin{tabular}{lrrr}
			\toprule
			\textbf{Training Dataset}  & \textbf{Training Time} & \textbf{Saved Epoch} & \textbf{Validation Loss} \\  
			\midrule
			\texttt{D0}   & 36.95 mins & 437 & 0.52514 \\
			\texttt{D1}   & 31.61 mins & 494 & 0.27424 \\
			\texttt{D2}   & 31.58 mins & 492 & 0.60667 \\
			\texttt{D3}  & 36.85 mins & 498 & 0.29508 \\
			\texttt{D4}   & 36.64 mins & 498 & 0.3199 \\
			\texttt{D5}  & 31.62 mins & 491 & 0.18456 \\
			\bottomrule
		\end{tabular}
	\end{minipage}
\end{table}
%******************************************************

\section{Numerical experiments}~\label{sec:experiments}
%******************************************************
This numerical section focus on illustrating different numerical features of these six trained U-Net neural operator preconditioners, including testing similar but unseen scenarios and checking network generalizability from various applications.
Most importantly, the dataset from the transcranial ultrasound therapy~\cite{skullDataset2018}.
Specifically, in these testing phases, numerical experiments are carried out with the discrete parametric Helmholtz equations~\eqref{eq:discrePDEs} that the neural operators have never seen before.
Without special notes, we display results of GMRES and \text{\FGMRES} (that is FGMRES preconditioned by the trained U-Net preconditioner with 4~\texttt{depth}) without restart.
A classical stopping criterion for the numerical linear algebra solvers is based on backward error analysis and consists of stopping the iteration when
\begin{equation}\label{eq:BEb}
	{\etab} =		\frac{\|A{\mathbf{u}} - {\boldsymbol{b}}\|} {\| {\boldsymbol{b}} \|} \le \tol, \text{ with } \tol = 10^{-12} \text{ by default}.
\end{equation}
The performance of the involved algorithms is evaluated in terms of the number of iterations, denoted as~$\its$, as well as the execution time, denoted as~$\ET$, required to converge. 
The maximum dimension of the Krylov search space $m$ shown in Algorithm~\ref{alg:nn-fgmres} is set to be 512, that is $m=512.$ The number of restart is $l=1 \text{ or } 10$ (here $l=1$ refers to the case without restart). Thus, the whole maximum iteration (denoted as $\maxit$) is $l \times 512$. 
We stop the algorithm when satisfying formula~\eqref{eq:BEb} or when $\maxit = l \times 512.$ 
Notation $\its_{max}, \its_{avg}, \its_{min}$ refers to the maximal, average, and minimum number of consumed $\its$ when testing multiple examples (that is $\K>1$).
Symbol~$(\its_{max}, \its_{avg}, \its_{min})_{E}^{*}$ denotes the algorithm diverges with $E$ linear systems among the whole $\K$ testing systems ($0 < E \leq \K$), or these $E$ linear systems are unable to reach the targeted accuracy within $\maxit$.
% %

Recall that this 2D~U-Net inference was trained on GPUs, whereas the classical subspace solvers have been generally implemented on CPUs.
In the following subsections, the GMRES and \text{\FGMRES} solvers are implemented in a \texttt{Python} prototype, supporting both CPU and GPU backends for computations with \texttt{numpy} and \texttt{cupy} Python libraries.
Without special notes, those two involved solvers in each of the subsections are running on the same 4~V100 GPUs devices used in the training process for a fair comparison.

\subsection{Neural operator preconditioning on the test dataset} 
\begin{figure}[!ht]%
    \centering%
    \begin{subfigure}{.49\textwidth}
        \includegraphics[width=\textwidth]{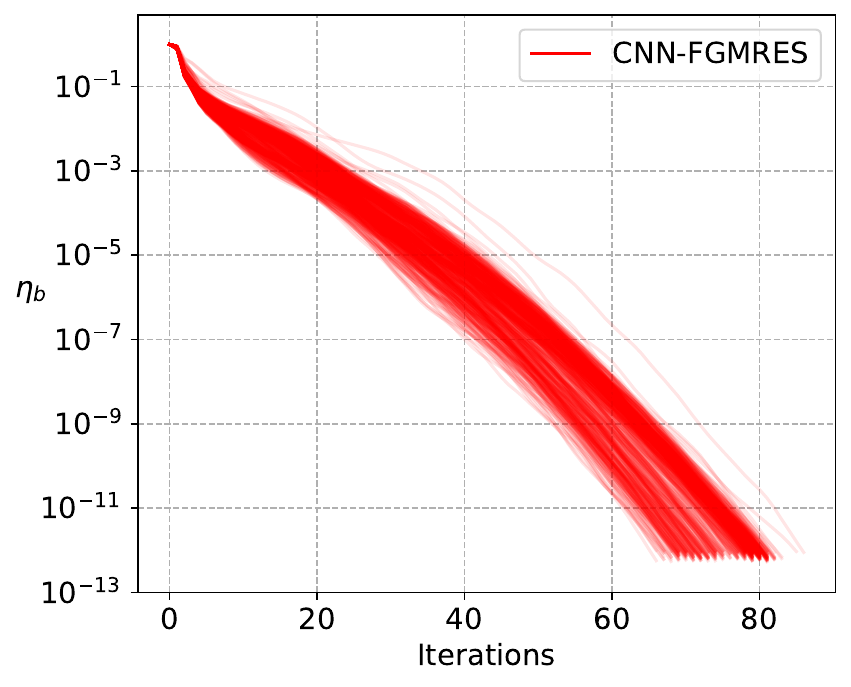}
        \caption{\FGMRES with $\K=1000,$ dataset \texttt{D0}}
    \end{subfigure}\hfill
    \begin{subfigure}{.49\textwidth}%
    \includegraphics[width=\textwidth]{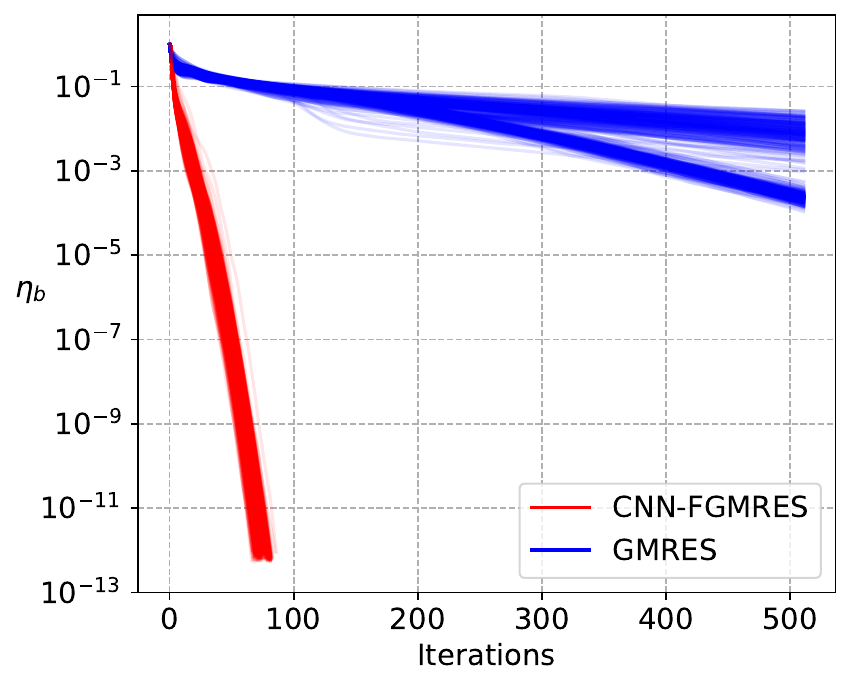}%
    \caption{GMRES, \FGMRES with $\K=500,$ dataset \texttt{D0}}%
    \end{subfigure}
    \caption{Results of the CNN-preconditioning trained with dataset \texttt{D0} in Table~\ref{tab:cases_mixed_dataset}: Convergence in terms of $\eta_b=10^{-12}$ of \FGMRES and GMRES with $\maxit=512$ without restrat for $\K$ examples from the test dataset with mixed data.
    }
    \label{fig:nn-fgmres_test_idea_skull}
\end{figure}

\begin{figure}[!ht]
    \centering
    \begin{subfigure}{.49\textwidth}
        \includegraphics[width=\textwidth]{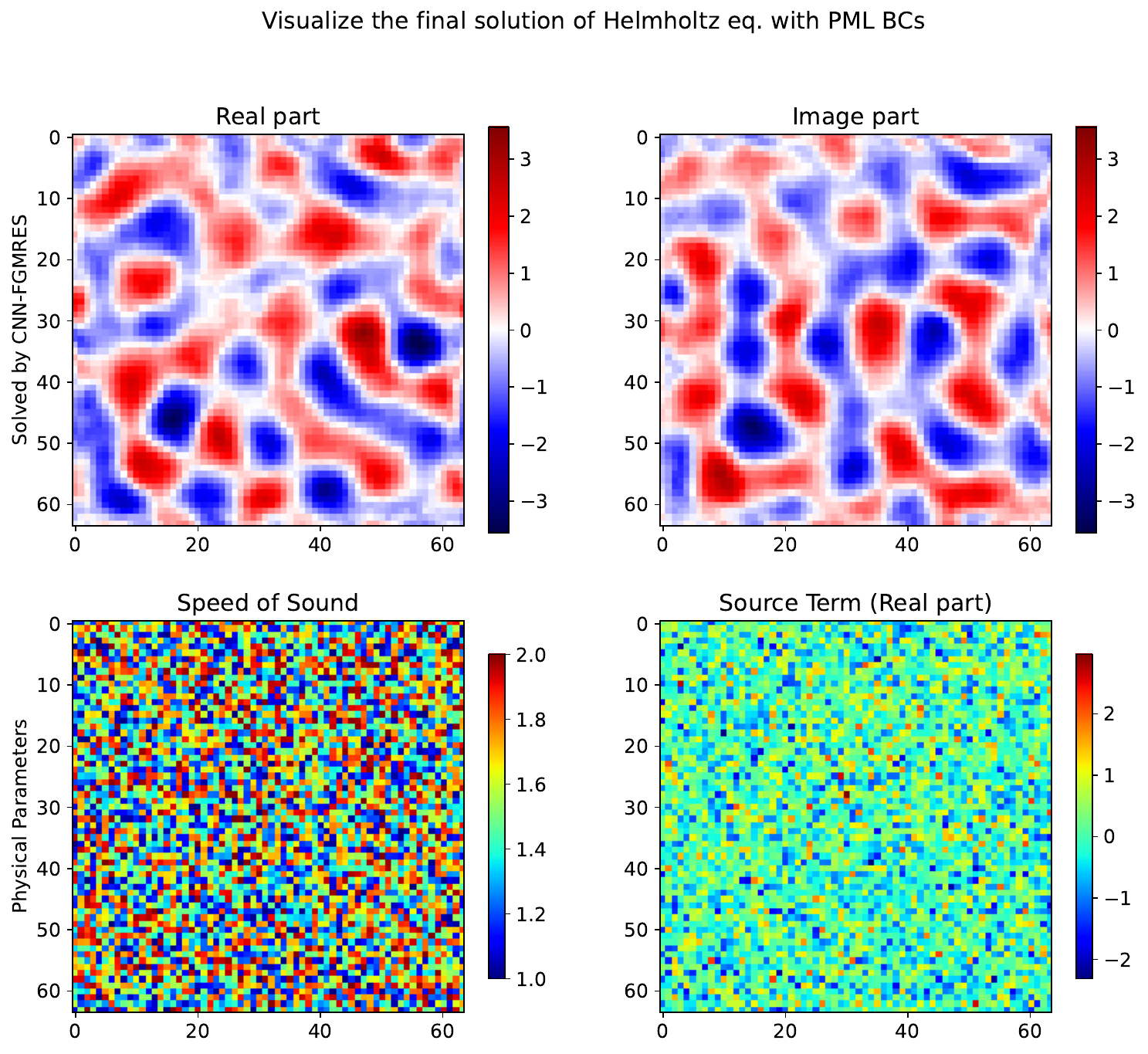} % 0, 3, 9, 11, 12, 14, 17, 18, 19
        \caption{$\bc \simiid {\cal U}([1,2]) + \bb \simiid {\cal N}(0,1)$}
    \end{subfigure}\hfill
    \begin{subfigure}{.49\textwidth}
        \includegraphics[width=\textwidth]{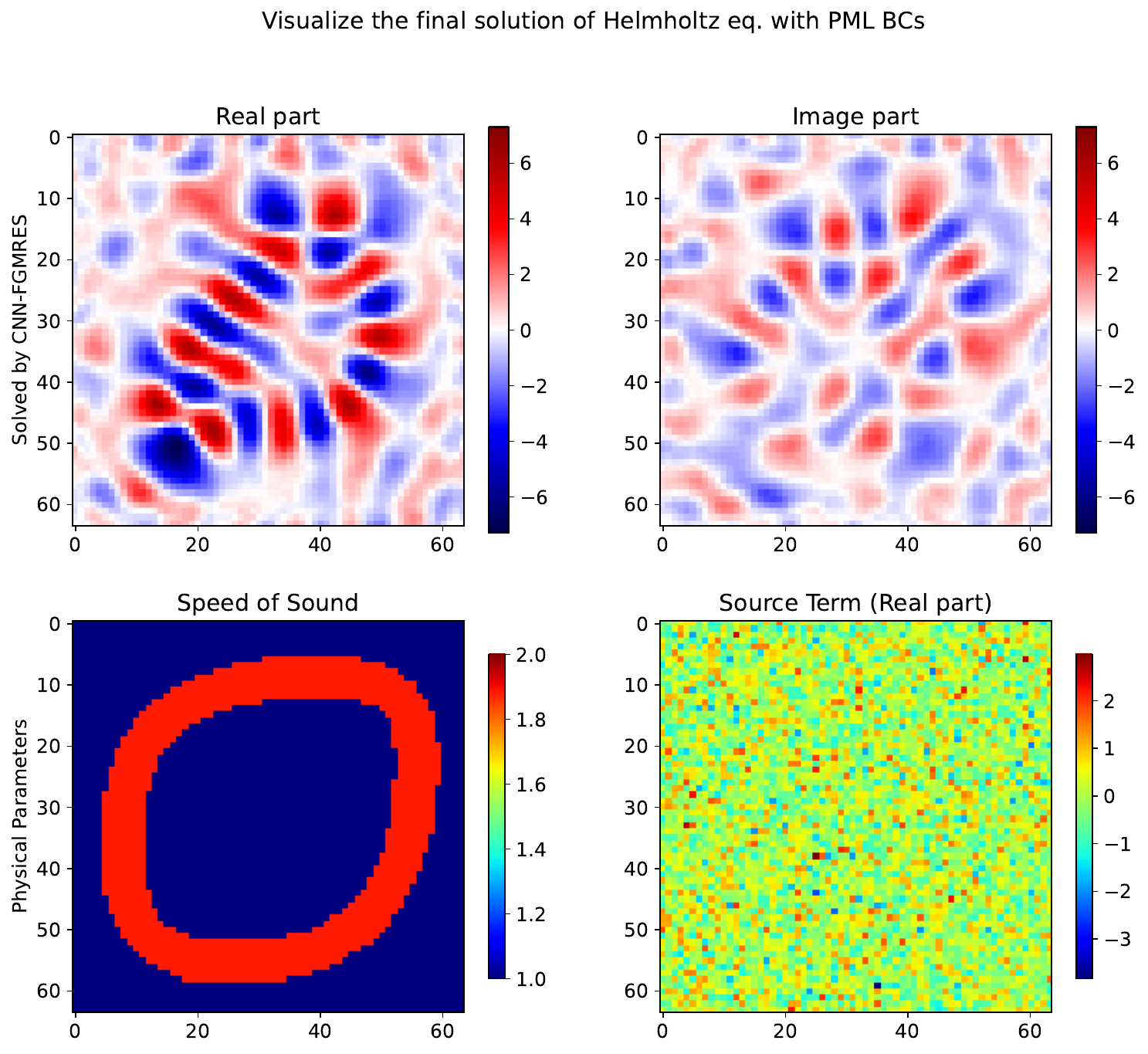} % 1, 7, 8, 15
        \caption{$\bc \simiid \text{\texttt{TruncIdealSkulls}} + \bb \simiid {\cal N}(0,1)$}
    \end{subfigure}
    \medskip
    \begin{subfigure}{.49\textwidth}
        \includegraphics[width=\textwidth]{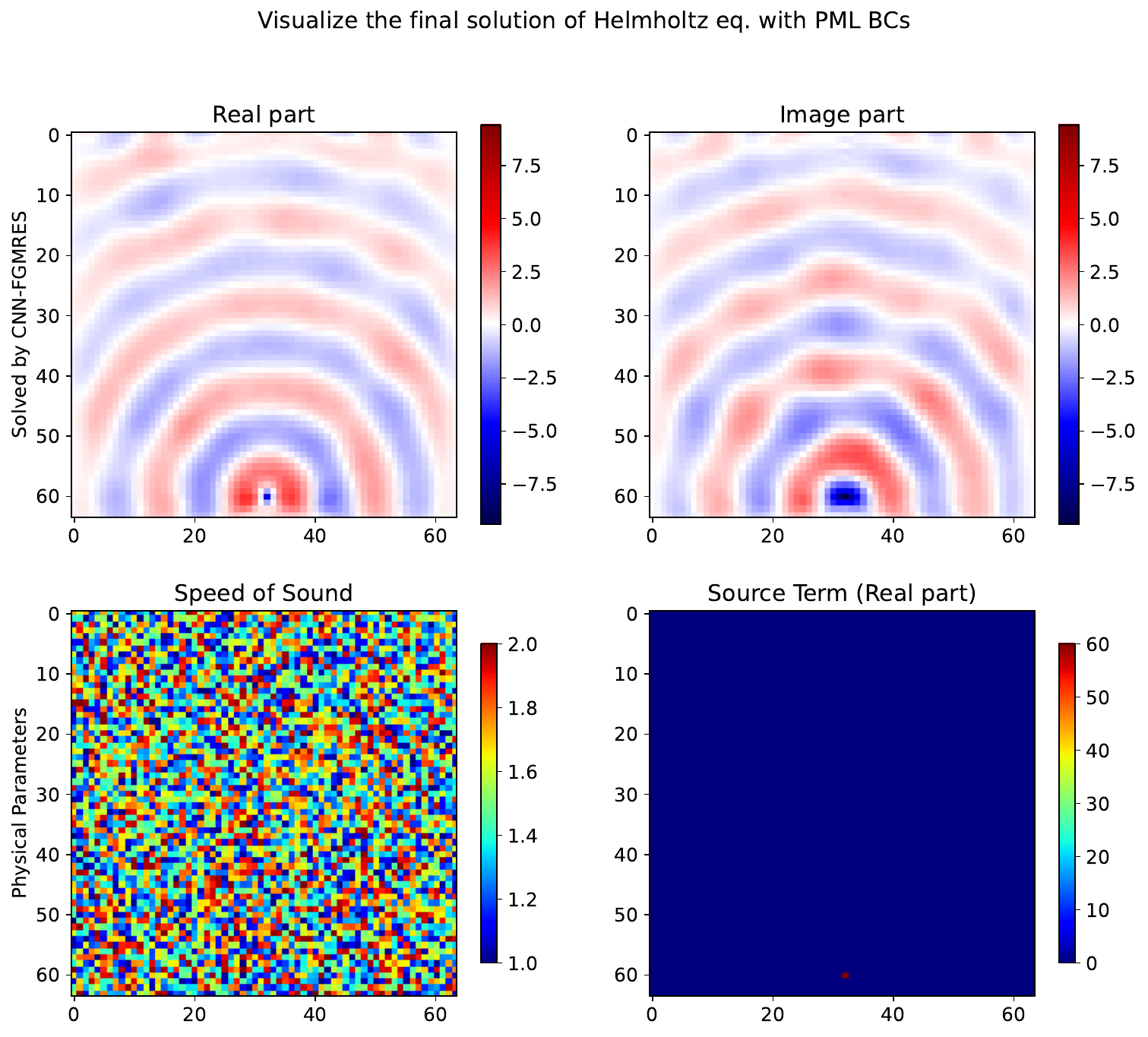}
        \caption{$\bc \simiid {\cal U}([1,2]) + \bb \simiid \text{\texttt{dirac1}}$} % 2 or 16
    \end{subfigure}\hfill
    \begin{subfigure}{.49\textwidth}
        \includegraphics[width=\textwidth]{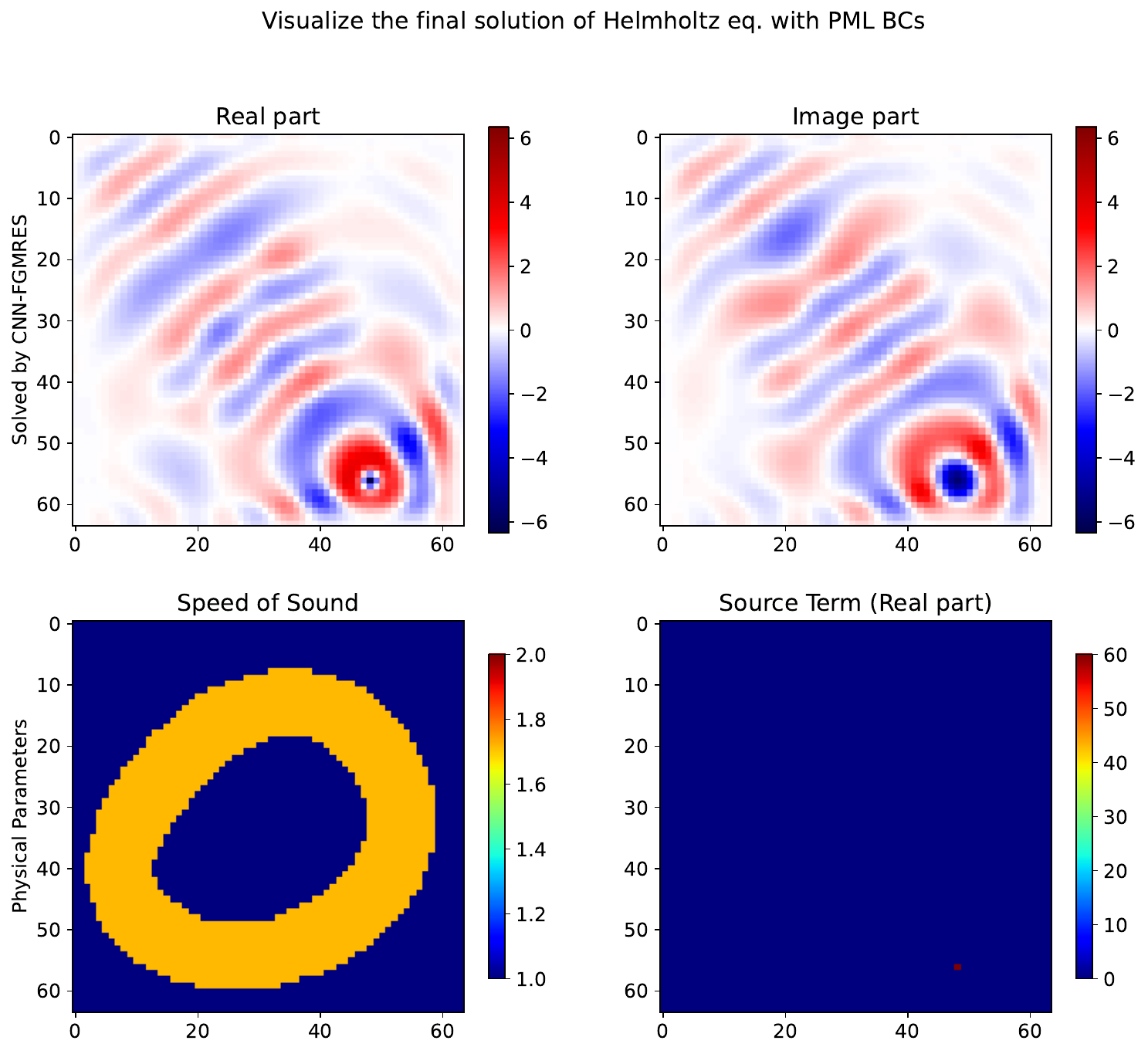} % 4, 5
        \caption{$\bc \simiid \text{\texttt{TruncIdealSkulls}} + \bb \simiid \text{\texttt{dirac}}$}
    \end{subfigure}
    \caption{Visualize the complex wavefield in real-imag part (upper two plots of each subplot (a) to (d)) and its corresponding velocity $\bc$ and source term $\bb$ (lower two plots of each subplot) of four examples solved by \FGMRES (the CNN was trained with dataset \texttt{D0}), and $\eta_b = 10^{-12},$ the four examples are from the test dataset with mixed data $(\bc, \bb).$}
    \label{fig:nn-fgmres_test_4_idea_skull}
\end{figure}

%******************************************************
We first assess the effectiveness of the six trained neural operator preconditioning inferences on problems similar to those used during the training, that is, linear systems arising from the discretization on a \text{$64 \times 64$}~grid, with physical parameters $\bc$ and $\bb$ drawn from the randomly mixed dataset \texttt{D0}-\texttt{D5} described in Table~\ref{tab:cases_mixed_dataset}.
The corresponding numerical results for GMRES and \FGMRES applied to the $\K$ linear systems are reported in Table~\ref{tab_results_test_idea_skull}. For the testing cases with $\K=500,$ we observe that GMRES (mostly) fails to reach the prescribed accuracy within $\maxit.$ In contrast, \FGMRES converges with significantly reduced computational cost, thanks to the trained neural operator preconditioner. This effectiveness and robustness are maintained when increasing the number of test problems to $\K=1000,$ thereby covering a broader range of physical parameter variations. This conclusion holds consistently across all trained models.
These results are further illustrated in Figure~\ref{fig:nn-fgmres_test_idea_skull}, which presents the convergence history of GMRES and \FGMRES (with the \text{CNN-preconditioning} trained on dataset \texttt{D0}) with $\K$ systems.
Figure~\ref{fig:nn-fgmres_test_4_idea_skull} provides visualizations of the obtained complex wavefield solutions under different configurations of the physical parameters.

%*** Begain Table 6 *****************************************
%%***********************************************************
\begin{table}[!htbp]
	\center{
		\begin{tabular}
			{@{}lllrr@{}}
			\toprule
			\raisebox{-1.50ex}[0cm][0cm]{\begin{tabular}{c} \# \\ {\small  Domain/Dataset}\end{tabular}}&
			\raisebox{-1.50ex}[0cm][0cm]{Method}&
			\raisebox{-1.50ex}[0cm][0cm]{$\K$}&
			\raisebox{-1.50ex}[0cm][0cm]{$(\its_{max}, \its_{avg}, \its_{min})$}&
			\raisebox{-1.50ex}[0cm][0cm]{$\ET$}\\
            &
			&
			&
			&\\
			\midrule
			\raisebox{-3.000ex}[0cm][0cm]{$64 \times 64$/\texttt{D0}}
            & GMRES & 500 & (512, 512, 512)$^{*}_{500}$ & 50359.62s \\ 
            & \FGMRES  & 500 & (87, 77, 67) & 478.97s \\ 
            & \FGMRES  & 1000 & (84, 77, 67) &  888.59s \\ 
			\midrule
			\raisebox{-3.000ex}[0cm][0cm]{~\qquad/\texttt{D1}}
            & GMRES & 500 & (512, 512, 512)$^{*}_{500}$ & 50370.80s \\ 
            & \FGMRES  & 500 & (51, 35, 30) & 134.40s \\
            & \FGMRES  & 1000 &  (50, 35, 30) & 261.95s \\ 
			\midrule
			\raisebox{-3.000ex}[0cm][0cm]{~\qquad/\texttt{D2}}
            & GMRES & 500 & (512, 511, 500)$^{*}_{493}$ &  51754.15s \\ 
            & \FGMRES  & 500 & (95, 79, 68) & 465.44s \\ 
            & \FGMRES  & 1000 &  (96, 79, 67) & 925.23s \\ 
			\midrule
			\raisebox{-3.000ex}[0cm][0cm]{~\qquad/\texttt{D3}}
            & GMRES & 500 & (512, 512, 512)$^{*}_{500}$ & 50367.13s \\ 
            & \FGMRES  & 500 & (49, 35, 26) &  136.70s \\  
            & \FGMRES  & 1000 & (49, 34, 26) & 237.14s \\ 
			\midrule
			\raisebox{-3.000ex}[0cm][0cm]{~\qquad/\texttt{D4}}
            & GMRES & 500 & (512, 512, 512)$^{*}_{500}$ & 51274.51s \\
            & \FGMRES  & 500 & (52, 38, 31) & 163.01s  \\
            & \FGMRES  & 1000 & (50, 38, 30) & 296.01s \\ 
			\midrule
			\raisebox{-3.000ex}[0cm][0cm]{~\qquad/\texttt{D5}}
            & GMRES  & 500 &  (512, 512, 512)$^{*}_{500}$ & 53071.27s \\  
            & \FGMRES & 500 & (29, 27, 27) & 114.50s \\ 
            & \FGMRES  & 1000 &  (29, 27, 27) &  195.50s \\ 
			\bottomrule
		\end{tabular}}
		\caption{{\small  Numerical results shown in Figure~\ref{fig:nn-fgmres_test_idea_skull} in terms of $(\its_{max}, \its_{avg}, \its_{min})$ and the GPUs~$\ET$ on 4~V100~GPUs for $\K=500 \text{ or } 1000$~examples from the test dataset with $\eta_b = 10^{-12}$ and $\maxit=512$ without restrat.}}\label{tab_results_test_idea_skull}
\end{table}
%%***********************************************************

%******************************************************
%******************************************************

%******************************************************
%******************************************************

On the other hand, the effectiveness of the neural operator preconditioner depends on the composition of the training dataset. For example, the model trained on dataset \texttt{D5}, which is constructed from a single configuration of the physical parameters $(\bc, \bb)$, yields the largest reduction in computational cost relative to GMRES. In contrast, the model trained on dataset \texttt{D2} (or \texttt{D0}), which contains a more diverse mixture of parameter configurations, exhibits relatively lower acceleration. This behavior is expected: training on a single configuration allows the network to specialize to that particular distribution, whereas mixed datasets require the model to represent a broader range of physical configurations. Consequently, mixed datasets generally trade a small amount of efficiency on individual configurations for improved robustness and generalization across diverse problems.
This trend is further illustrated by comparing the results obtained with datasets \texttt{D0} and \texttt{D1}, as well as with datasets \texttt{D2} and \texttt{D3}. In each pair, the velocity distributions $\bc$ are generated using the same mixed strategy, while the source fields $\bb$ differ between a mixed distribution and a single structured configuration (see Table~\ref{tab:cases_mixed_dataset}).

\subsection{Network generalizability in human head~CT scan examples}~\label{sec:gen} 
%%***********************************************************
The following subsections investigate the generalization capabilities of the trained CNN-based preconditioning operators in the primary application of head~CT scan examples. To evaluate their generalizability, we test their performance under various heterogeneous parameter configurations in the parametric Helmholtz equation~\eqref{eq:discrePDEs}, extending beyond the settings used during training.
We first consider the truncated idealized skulls on a $64 \times 64$ grid with variations in terms of the source field in Section~\ref{sec:raw_nn_gen_b}. Then, in Section~\ref{subsec:skull_exp}, we focus on the primary application involving a practical transcranial ultrasound dataset derived from adult human head~CT scans, where both the source and velocity fields exhibit structured heterogeneity on a $512 \times 512$ grid, which is 64~times larger than that used in training.

It is worth noting that, for challenging large-scale problems such as the practical transcranial ultrasound applications considered in Section~\ref{subsec:skull_exp}, many classical algebraic preconditioners—including Jacobi, algebraic multigrid, and incomplete LU (ILU) or incomplete Cholesky factorizations—may become impractical. Their construction and application typically require explicit access to the matrix representation of the discretized operator, which can lead to prohibitive memory requirements and setup costs on fine computational grids. 
In contrast, the proposed neural operator preconditioner, which approximates the inverse of the discretized differential operator, is inherently matrix-free and is applied directly through network inference, thereby naturally avoiding the assembly, storage, and factorization of large sparse matrices.
For further details and a quantitative comparison between neural operator preconditioners and several classical algebraic preconditioning techniques on moderate-scale academic examples, we refer the reader to~\cite[Section~4.2]{Xiang2025_NeuralNetworkPreconditioning, Giraud2025_NeuralNetworkPreconditioning}.

\subsubsection{Test on the idealized skulls with fixed and smoothed source field}\label{sec:raw_nn_gen_b} 
%%***********************************************************
Before evaluating network generalizability of the models trained with the six dataset described in Table~\ref{tab:cases_mixed_dataset} on the more challenging transcranial ultrasound dataset, we first assess their performance on 1000 idealized skull samples (truncated test idealized skulls introduced in~\cite{Stanziola2021_HelmholtzEquationSolver}) with some given source fields. In this setting, the velocity field is given by $\bc \simiid \text{\texttt{TruncIdealSkulls}}.$
For these 1000 realizations of $\bc$, which differ from those used during training, we consider the following two configurations of the source term $\bb$:
\begin{enumerate}
    \item $\bb \simiid \text{\texttt{dirac1}}$: a single nonzero entry of $\bb$ is set to 1 at the fixed grid location $(60, 32)$, while all other entries are zero;

    \item $\bb \simiid \texttt{dirac1\_smooth}$: a smoothed version of $\bb \simiid \text{\texttt{dirac1}}$, defined as a localized Gaussian source $\bb(\bx, \by) = \exp^{{- \left((x-c_x)^2 + (x-c_y)^2\right)} / {2 \sigma^2}},$
    where $\sigma = 2$ and $(c_x,c_y) = (60,32)$ denotes the grid location.
\end{enumerate}
With these physical parameter configurations, the numerical results for solving the $\K$ linear systems arising from Equation~\eqref{eq:discrePDEs} are reported in Table~\ref{tab_results_test_dirac1_idea_skull} in terms of $(\its_{max}, \its_{avg}, \its_{min})$ and the elapsed time $\ET$.

%***Begain Table 7*****************************************
%%***********************************************************
\begin{table}[!htbp]
	\center{
		\begin{tabular}
			{@{}llllrr@{}}
			\toprule
            \raisebox{-1.50ex}[0cm][0cm]{\begin{tabular}{c} \# \\ {\small  Source term}\end{tabular}}&
			\raisebox{-1.50ex}[0cm][0cm]{\begin{tabular}{c} \# \\ {\small  Domain/Dataset}\end{tabular}}&
			\raisebox{-1.50ex}[0cm][0cm]{Method}&
			\raisebox{-1.50ex}[0cm][0cm]{$\K$}&
			\raisebox{-1.50ex}[0cm][0cm]{$(\its_{max}, \its_{avg}, \its_{min})$}&
			\raisebox{-1.50ex}[0cm][0cm]{$\ET$}\\
            &
            &
			&
			&
			&\\
			\midrule
            \raisebox{-20.00ex}[0cm][0cm]{\begin{tabular}{c} $\bb  \simiid \text{\texttt{dirac1}}$ \end{tabular}}
			& \raisebox{-0.00ex}[0cm][0cm]{$64 \times 64$/}
            & GMRES & 20 & (512, 512, 512)$^{*}_{20}$ & 2142.70s  \\  
            \cmidrule{2-6}
            \raisebox{-20.00ex}[0cm][0cm]{}
			& \raisebox{-1.500ex}[0cm][0cm]{~\qquad/\texttt{D0}}
            & \FGMRES  & 20 & (79, 72, 67) &   21.95s \\ 
            & & \FGMRES  & 1000 & (89, 72, 66) & 810.03s \\ 
            \cmidrule{3-6}
			& \raisebox{-1.500ex}[0cm][0cm]{~\qquad/\texttt{D1}}
            & \FGMRES  & 20 & (47, 40, 33) & 17.87s \\ 
            & & \FGMRES  & 1000 & (51, 40, 32) & 306.07s \\ 
            \cmidrule{3-6}
			& \raisebox{-1.500ex}[0cm][0cm]{~\qquad/\texttt{D2}}
            & \FGMRES  & 20 & (83, 74, 70) & 18.69s \\ 
			& & \FGMRES  & 1000 & (93, 76, 67) & 877.26s  \\ 
            \cmidrule{3-6}
			& \raisebox{-1.500ex}[0cm][0cm]{~\qquad/\texttt{D3}}
            & \FGMRES  & 20 & (46, 40, 35) & 13.08s \\ 
            &  & \FGMRES  & 1000 &  (50, 40, 32) & 295.92s \\ 
            \cmidrule{3-6}
			&  \raisebox{-1.500ex}[0cm][0cm]{~\qquad/\texttt{D4 $\star$}}
            & \FGMRES  & 20 & (47, 38, 33) & 7.93s \\ 
            &  & \FGMRES  & 1000 & (51, 38, 30) & 289.49s \\ 
            \cmidrule{3-6}
			&  \raisebox{-1.500ex}[0cm][0cm]{~\qquad/\texttt{D5}}
            & \FGMRES  & 20 & (134, 107, 81) & 39.17s \\ 
            &  & \FGMRES  & 1000 &  (141, 111, 76) & 1696.05s \\ 
            \midrule
            \raisebox{-20.00ex}[0cm][0cm]{\begin{tabular}{c} $\bb  \simiid \texttt{dirac1\_smooth}$ \end{tabular}}
			&  \raisebox{-0.000ex}[0cm][0cm]{$64 \times 64$/}
            & GMRES & 20 & (512, 512, 512)$^{*}_{20}$ & 2135.96s \\ 
            \cmidrule{2-6}
            \raisebox{-20.00ex}[0cm][0cm]{}
			&  \raisebox{-1.500ex}[0cm][0cm]{~\qquad/\texttt{D0}}
            & \FGMRES  & 20 & (81, 73, 68) & 18.10s \\ 
            &  & \FGMRES  & 1000 & (92, 73, 68) &  820.14s \\ 
            \cmidrule{3-6}
			&  \raisebox{-1.500ex}[0cm][0cm]{~\qquad/\texttt{D1}}
            & \FGMRES  & 20 & (48, 40, 34) & 8.14s \\ 
            &  & \FGMRES  & 1000 & (52, 40, 32) & 313.74s \\ 
            \cmidrule{3-6}
            &  \raisebox{-1.500ex}[0cm][0cm]{~\qquad/\texttt{D2}}
            & \FGMRES  & 20 &   (85, 77, 72) & 27.15s \\ 
            &  & \FGMRES  & 1000 & (96, 77, 69) & 893.96s \\ 
            \cmidrule{3-6}
			&  \raisebox{-1.500ex}[0cm][0cm]{~\qquad/\texttt{D3}}
            & \FGMRES  & 20 &  (47, 40, 34) & 9.11s \\
            &  & \FGMRES  & 1000 & (51, 40, 32) & 313.31s \\ 
            \cmidrule{3-6}
			&  \raisebox{-1.500ex}[0cm][0cm]{~\qquad/\texttt{D4 $\star$}}
            & \FGMRES  & 20 &  (45, 36, 33) &  7.35s \\ 
            &  & \FGMRES  & 1000 & (52, 38, 31) & 295.11s \\
            \cmidrule{3-6}
			&  \raisebox{-1.500ex}[0cm][0cm]{~\qquad/\texttt{D5}}
            & \FGMRES  & 20 & (137, 110, 84) & 37.15s \\ 
            &  & \FGMRES  & 1000 & (142, 111, 76) & 1784.98s \\ 
			\bottomrule
		\end{tabular}}
		\caption{{\small  Numerical results shown in Figure~\ref{fig:nn-fgmres_test_dirac1_idea_skull}
        in terms of $(\its_{max}, \its_{avg}, \its_{min})$ and the GPUs~$\ET$ on 4~V100~GPUs for \text{$\K=20, \text{ or } 1000$} examples with varying velocity $\bc \simiid \text{\texttt{TruncIdealSkulls}}$, and the fixed source term \text{$\bb  \simiid \text{\texttt{dirac1}},$} and its polluted one $\bb  \simiid \texttt{dirac1\_smooth}$ with center radial decay Gaussian noise by $\texttt{radius}=5, \texttt{noise\_level}=0.05.$ The $\eta_b = 10^{-12}$ and $\maxit=512$ without restrat.}}\label{tab_results_test_dirac1_idea_skull}
\end{table}
%%***********************************************************

%******************************************************
\begin{figure}[!ht]%
    \centering%
    \begin{subfigure}{.49\textwidth}
        \includegraphics[width=\textwidth]{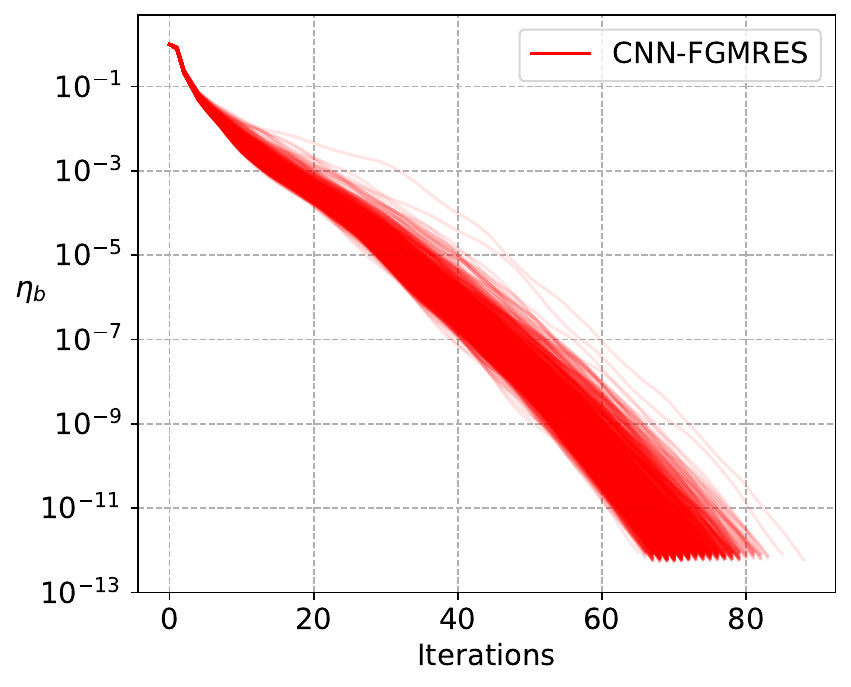}
        \caption{\FGMRES with $\K=1000, \bb \simiid \text{\texttt{dirac1}},$ dataset~\texttt{D0}}
    \end{subfigure}\hfill
    \begin{subfigure}{.49\textwidth}%
    \includegraphics[width=\textwidth]{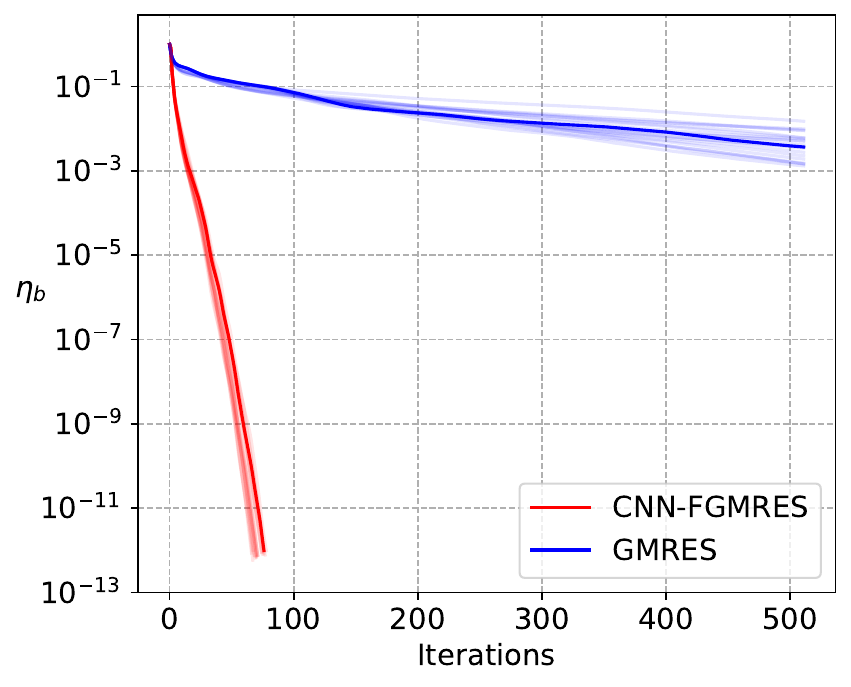}%
    \caption{\FGMRES, GMRES with $\K=20, \bb \simiid \text{\texttt{dirac1}},$ dataset \texttt{D0}}%
    \end{subfigure}
    \vskip .1cm
    \medskip
    \begin{subfigure}{.50\textwidth}
        \includegraphics[width=\textwidth]{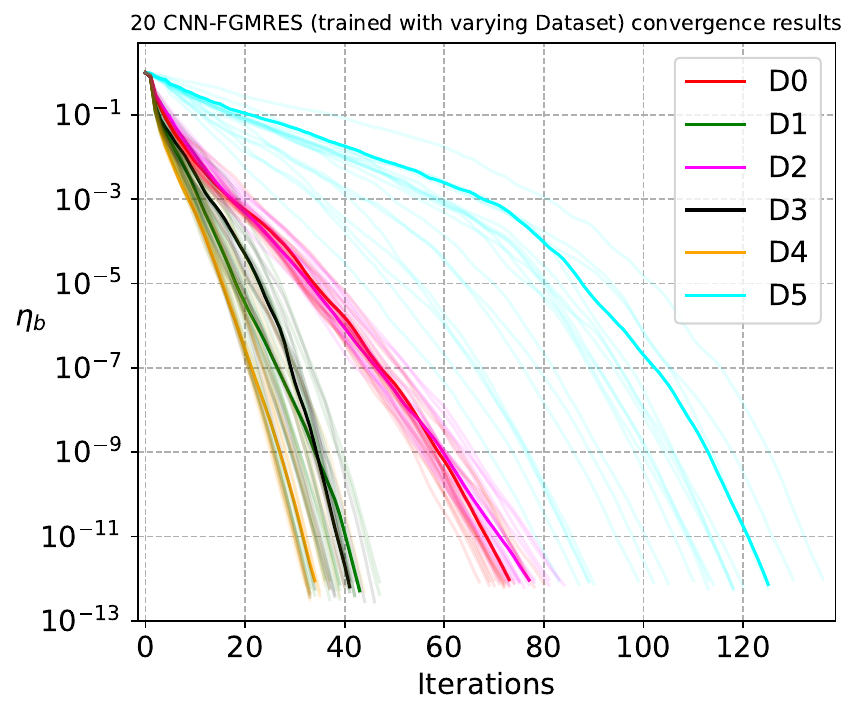}
        \caption{\FGMRES with $\K=20,$ $\bb \simiid \texttt{dirac1\_smooth},$ and CNN was trained with the six different datasets \texttt{D0}-\texttt{D5}} 
    \end{subfigure}
    \caption{Convergence in terms of $\eta_b=10^{-12}$ of \FGMRES and GMRES with $\maxit=512$ without restrat for $\K$ examples with varying velocity $\bc \simiid \text{\texttt{TruncIdealSkulls}}$ from the test dataset and the fixed source term $\bb  \simiid \text{\texttt{dirac1}}$ with fixed grid-point locates at $60 \times 32$ within 2D domain $64 \times 64,$ and its polluted one $\bb  \simiid \texttt{dirac1\_smooth}$ with center radial decay Gaussian noise by $\texttt{radius}=5, \texttt{noise\_level}=0.05.$ 
    In (a)-(b), the presented \FGMRES results are for the CNN-preconditioning trained with dataset \texttt{D0}, and $\bb  \simiid \text{\texttt{dirac1}}.$ In (c), it shows corresponding results from the CNN-preconditioning trained with the six different datasets \texttt{D0}-\texttt{D5} described in Table~\ref{tab:cases_mixed_dataset}, and $\bb  \simiid \texttt{dirac1\_smooth}.$}
    \label{fig:nn-fgmres_test_dirac1_idea_skull}
\end{figure}
%******************************************************

As in the previous experiments, the application of the trained CNN-based operator preconditioners enables \text{\FGMRES} to converge with substantially fewer iterations and lower computational cost than GMRES, which fails to achieve the prescribed accuracy within $\maxit.$ Moreover, the robustness of the proposed approach is maintained when the number of test problems is increased to $\K=1000.$ The corresponding convergence histories for the model trained on dataset~\texttt{D0} are shown in Figure~\ref{fig:nn-fgmres_test_dirac1_idea_skull}~(a)-(b).
Among the six trained models, the operator preconditioner trained on dataset~\texttt{D4} exhibits the best performance. This is expected, since the tested velocity fields $\bc \simiid \text{\texttt{TruncIdealSkulls}}$ share structural similarities with those included in dataset~\texttt{D4}. Nevertheless, it is noteworthy that the operator preconditioner trained on the fully random dataset~\texttt{D5}, which contains no particular geometric structure, still provides a significant acceleration of \FGMRES, albeit with somewhat reduced effectiveness. The differences among the six trained preconditioners can be clearly observed in Figure~\ref{fig:nn-fgmres_test_dirac1_idea_skull}~(c), which presents the convergence histories of \FGMRES for a set of 20 linear systems.
These results indicate that the choice of the training dataset does impact performance on some targed testing cases. This effect becomes even more pronounced in the practical transcranial ultrasound examples presented in the next section.

\begin{figure}[!ht]
    \centering
    \begin{subfigure}{.45\textwidth}  
        \includegraphics[width=\textwidth]{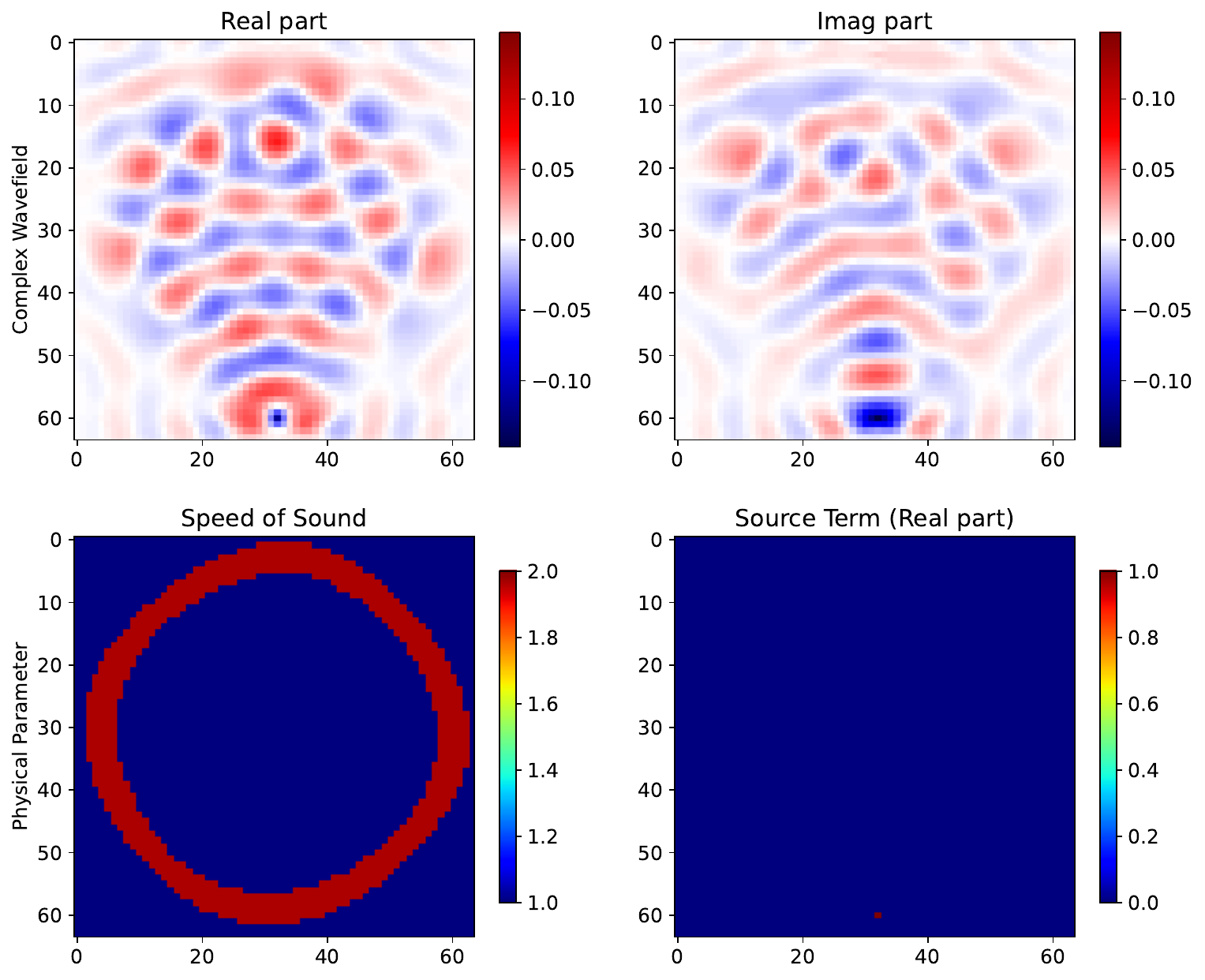}
        \caption{$\bb  \simiid \text{\texttt{dirac1}}$}
    \end{subfigure}\hfill
    \begin{subfigure}{.45\textwidth}
        \includegraphics[width=\textwidth]{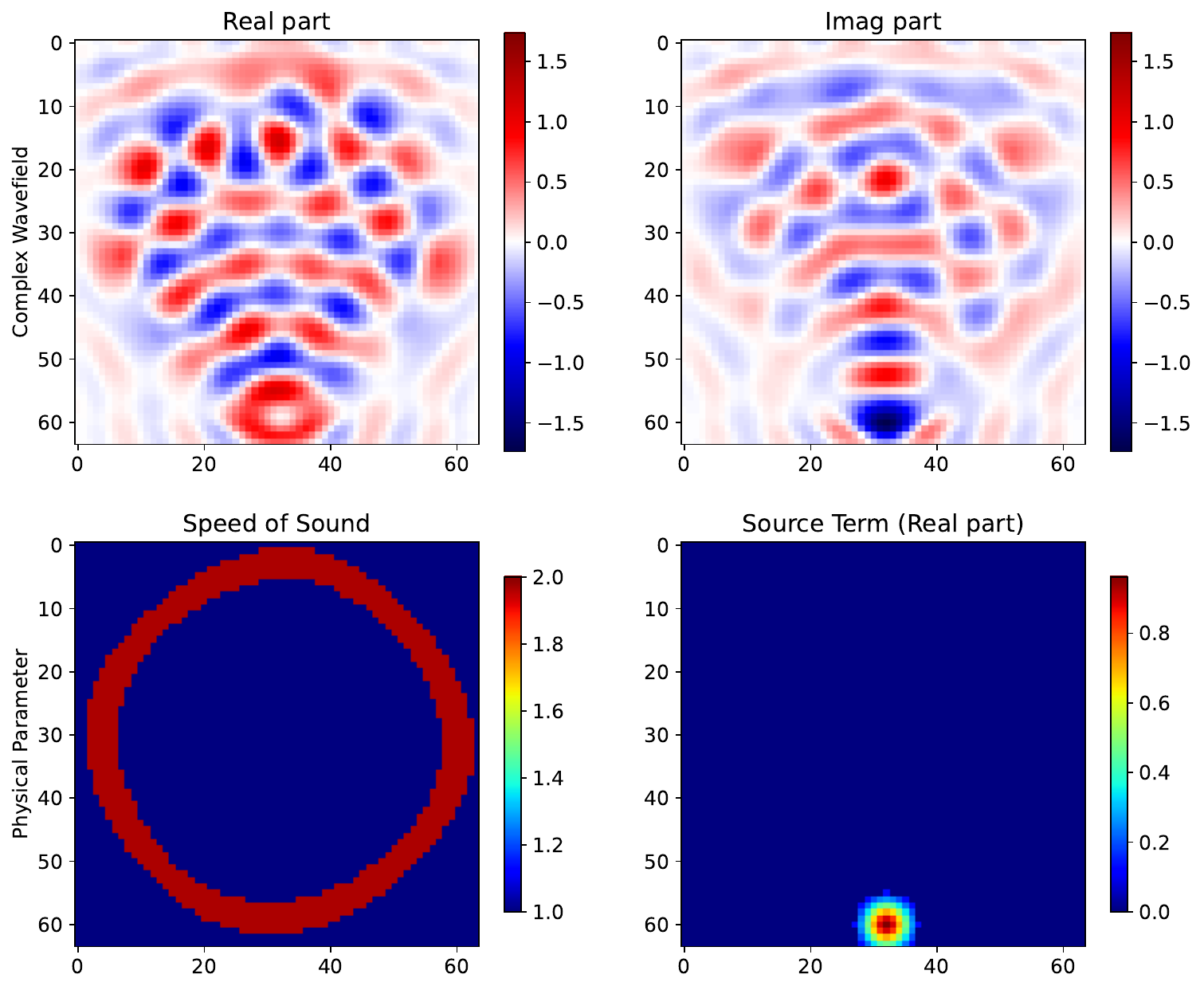}
        \caption{polluted $\bb  \simiid \texttt{dirac1\_smooth}$}
    \end{subfigure}
    \medskip
    \begin{subfigure}{.45\textwidth}
        \includegraphics[width=\textwidth]{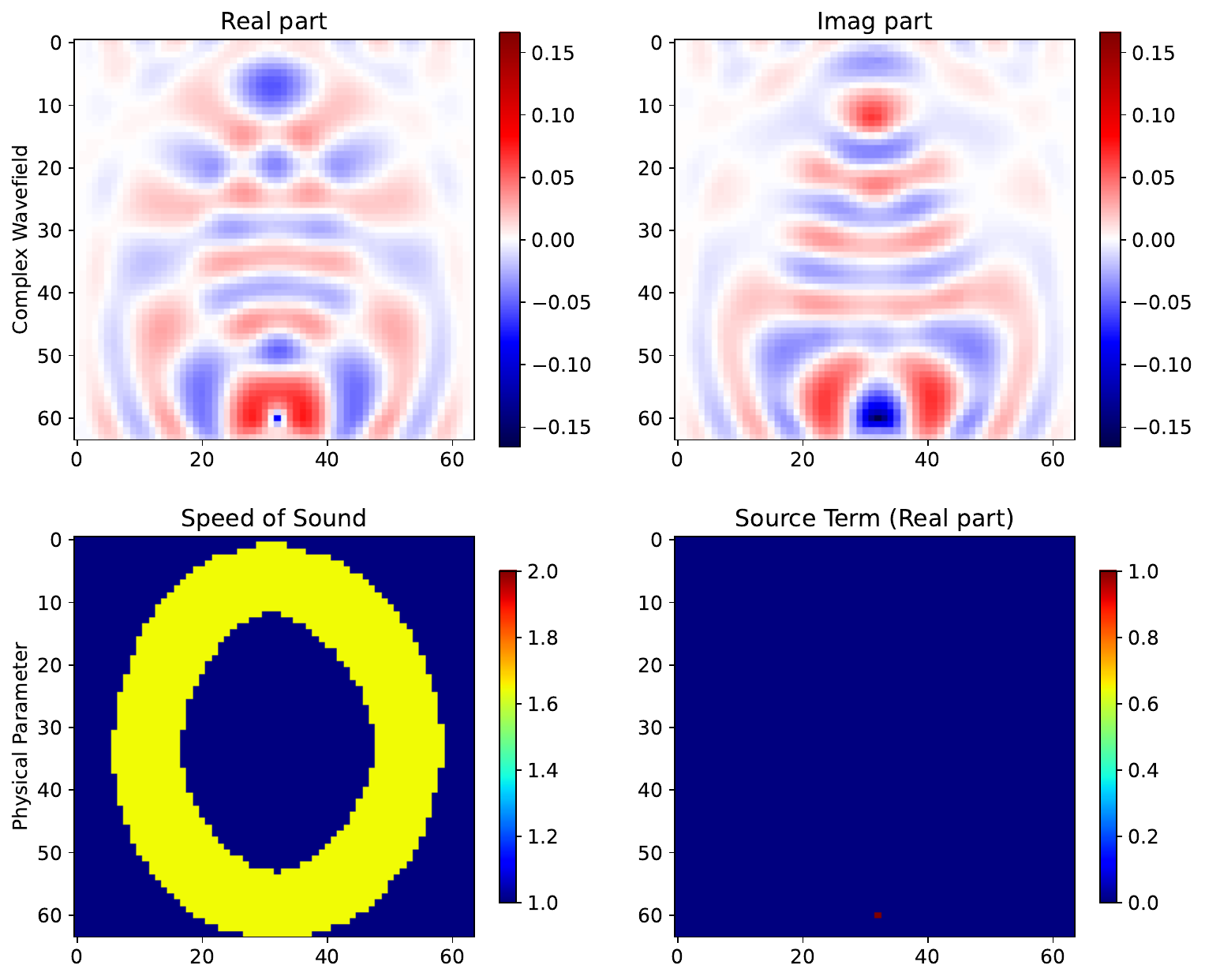}
        \caption{$\bb  \simiid \text{\texttt{dirac1}}$}
    \end{subfigure}\hfill
    \begin{subfigure}{.45\textwidth}
        \includegraphics[width=\textwidth]{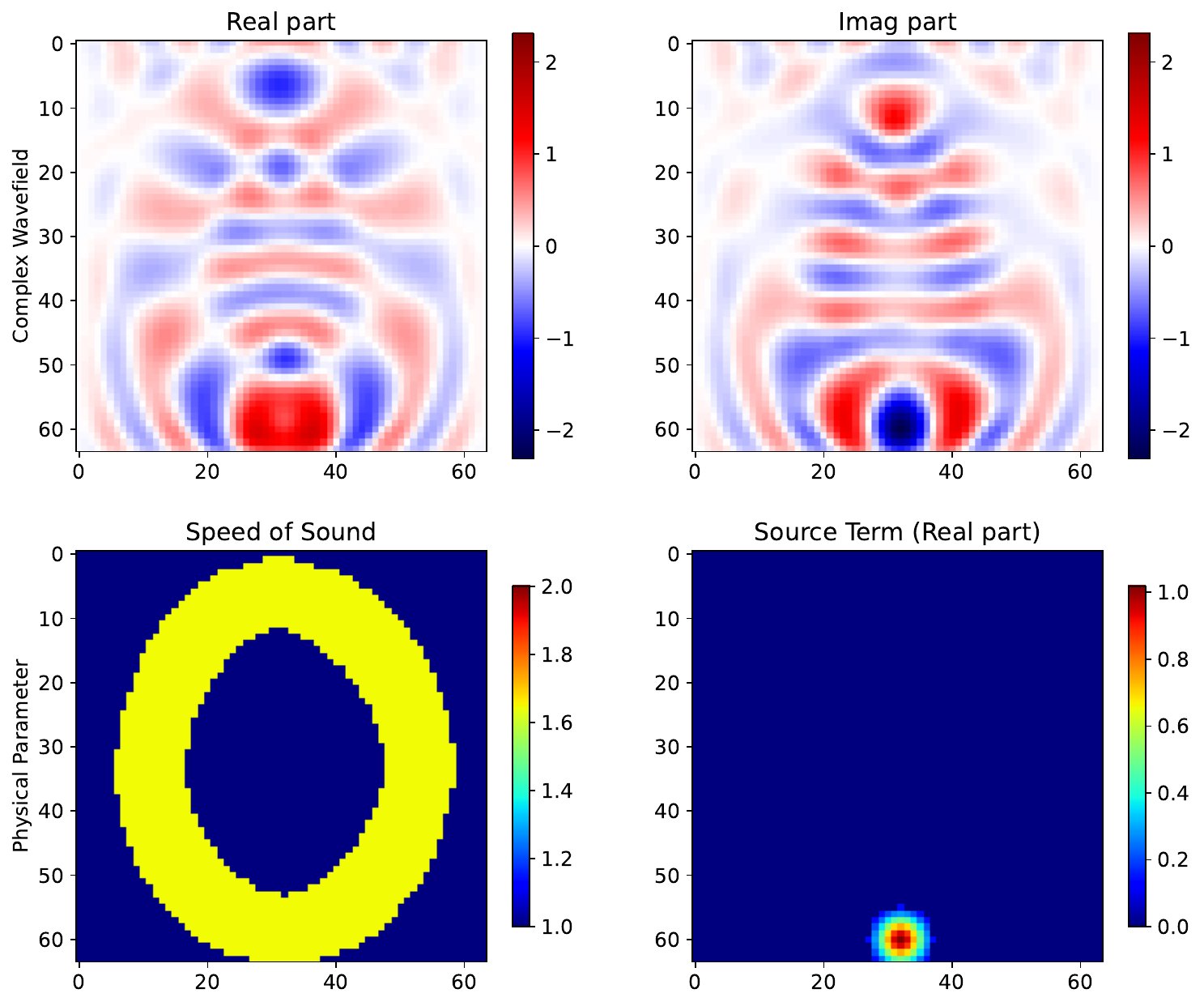}
        \caption{polluted $\bb  \simiid \texttt{dirac1\_smooth}$}
    \end{subfigure}
    \caption{Visualize the complex wavefield in real-imag part (upper two plots of each subplot (a) to (d)) and its corresponding velocity $\bc$ and source term $\bb$ (lower two plots of each subplot)  of two examples solved by \FGMRES (the CNN was trained with dataset \texttt{D0}), and $\eta_b = 10^{-12},$ where the two varying $\bc \simiid \text{\texttt{TruncIdealSkulls}}$ is from the test dataset, and the fixed $\bb  \simiid \text{\texttt{dirac1}}$ with fixed grid-point locates at $60 \times 32$ within 2D domain $64 \times 64,$ and its polluted one $\bb  \simiid \texttt{dirac1\_smooth}$ with center radial decay Gaussian noise by $\texttt{radius}=5, \texttt{noise\_level}=0.05.$}
    \label{fig:nn-fgmres_test_4_dirac1_idea_skull}
\end{figure}

On the other hand, no significant differences are observed between the source field configurations $\bb \simiid \text{\texttt{dirac1}}$ and $\bb \simiid \texttt{dirac1\_smooth}$ in the numerical results reported in Table~\ref{tab_results_test_dirac1_idea_skull}. To further investigate their effects, we visualize the corresponding wavefields in Figure~\ref{fig:nn-fgmres_test_4_dirac1_idea_skull}. The smoothed-source configuration produces less localized irregularity near the source location, thereby providing a clearer visualization of both the velocity distribution and the resulting wavefield patterns.

%******************************************************
%******************************************************

\subsubsection{Application to practical transcranial ultrasound examples}~\label{subsec:skull_exp}
%%***********************************************************
In this section, we evaluate the trained neural operator preconditioners on four practical adult human head~CT scan examples obtained from {qure.ai}~\cite{skullDataset2018}. 
The dataset is publicly available through Kaggle.\footnote{\url{https://www.kaggle.com/datasets/crawford/qureai-headct}}
Each example is constructed from a transverse~CT slice of an adult human skull, yielding a heterogeneous velocity $\bc \in \R^{512^2}$ distributed on a $512 \times 512$ grid. 
This computational grid is 64 times larger than the $64 \times 64$ grid used during training, providing a significantly challenging out-of-distribution testing case. 
Visualizations of the four velocity fields are shown in the lower-left panels of Figure~\ref{fig:nn-fgmres}~(a)-(d). The corresponding source field $\bb \in \R^{512^2}$ is modeled as a focused ultrasound transducer represented by a one-dimensional arc, illustrated by the red curve in the lower-right panels of Figure~\ref{fig:nn-fgmres} for its visualization. The transducer aperture diameter and radius of curvature are both set to 60~mm, and the source frequency is fixed to 490~kHz. 

To solve these challenging large-scale problems, restarted GMRES and \FGMRES are employed with restart parameter $l=10$, resulting in a maximum of $\maxit=5120$ iterations throughout this section.
Furthermore, in addition to GMRES, we also benchmark our approach against the closely related state-of-the-art neural network solver: the learned optimizer proposed by Stanziola {et al.} in~\cite{Stanziola2021_HelmholtzEquationSolver}.
Specifically, we use their released \texttt{IterativeSolver} (i.e., the learned optimizer) together with the pretrained checkpoint provided in the authors' GitHub repository\footnote{\url{https://github.com/ucl-bug/helmnet/tree/main}} and run the learned optimizer solver up to 5120 nonlinear Richardson iterations.
Since the downloaded physical parameter pairs $(\bc, \bb)$ are loaded directly into CPU memory, both the learned optimizer, GMRES and \FGMRES are executed on the same CPU to avoid unnecessary data transfers between CPU and GPU devices. Specifically, all experiments in this section are performed on the CAIUS cluster, %\footnote{Le Cluster de cAlcul Intensif à l'Université de Strasbourg (CAIUS): \url{https://hpc.pages.unistra.fr}}, 
using Intel Skylake compute nodes equipped with 2x 16-core Intel Xeon Gold 6126 processors running at 2.60~GHz and 192~GB of RAM per node.

%***Begain Table 8*****************************************
%%***********************************************************
\begin{table}[!htbp]
	\center{
		\begin{tabular}
			{@{}llllrr@{}}
			\toprule
			\raisebox{-1.50ex}[0cm][0cm]{\begin{tabular}{c} \# \\ {\small  Domain/Dataset}\end{tabular}}&
			\raisebox{-1.50ex}[0cm][0cm]{Method}&
			\raisebox{-1.50ex}[0cm][0cm]{$\K$}&
      \raisebox{-1.50ex}[0cm][0cm]{$\eta_b$}&
			\raisebox{-1.50ex}[0cm][0cm]{$[\its_{1}, \ldots, \its_{\K}]$}&
			\raisebox{-1.50ex}[0cm][0cm]{$\ET$}\\
            &
            &
			&
			&
			&\\
			\midrule
			\raisebox{-0.000ex}[0cm][0cm]{$512 \times 512$/}
            & GMRES & 4 &   $10^{-3}$ & [5120$^{*}$, 4546, 5120$^{*}$, 5120$^{*}$] & 44967.68s {\tiny (12.49 hours)} \\ %16683550 (5120$^{*}_{3}$, 4976, 4546)
            \midrule
			\raisebox{-1.500ex}[0cm][0cm]{~\qquad/\texttt{D0} $\star$} 
            & \FGMRES  & 4 & $10^{-3}$ & [191, 159, 166, 162] & 737.70s {\tiny (12.30 mins)} \\ % (191, 169, 159)
            & \FGMRES  & 4 & $10^{-6}$ & [559, 390, 392, 382] & 2773.70s {\tiny (46.23 mins)} \\ % 16684204 (559, 430, 382)
			\cmidrule{2-6}
			\raisebox{-1.500ex}[0cm][0cm]{~\qquad/\texttt{D1}}
            & \FGMRES  & 4 & $10^{-3}$ &  [367, 173, 190, 207] &  1278.78s {\tiny (21.31 mins)} \\ % 16683554* (367, 234, 173) /
            & \FGMRES  & 4 & $10^{-6}$ &  [1217, 347, 359, 417] & 4405.26s {\tiny (73.42 mins)} \\ % 16683684 (1217, 585, 347)
			\cmidrule{2-6}
			\raisebox{-1.50ex}[0cm][0cm]{~\qquad/\texttt{D2}}
            % & \FGMRES  & 4 & $10^{-2}$ & (113, 98, 87) [94, 87, 113, 101] & 280.26s (run on skylake 24 CPUs) \\
			& \FGMRES  & 4 & $10^{-3}$ & [403, 149, 158, 170] & 1472.62s {\tiny (24.54 mins)} \\ % 16683551 (403, 220, 149)
            & \FGMRES  & 4 & $10^{-6}$ & [1467, 378, 391, 404] & 5265.71s {\tiny (87.76 mins)} \\ % 16683683 (1467, 660, 378)
			\cmidrule{2-6}
			\raisebox{-0.000ex}[0cm][0cm]{~\qquad/\texttt{D3}}
            & \FGMRES  & 4 & $10^{-3}$ &  [395, 198, 216, 240] & 1749.17s {\tiny (29.15 mins)} \\ % 16683553 (395, 262, 198)
			\cmidrule{2-6}
			\raisebox{-0.000ex}[0cm][0cm]{~\qquad/\texttt{D4}}
            & \FGMRES  & 4 & $10^{-3}$ & [627, 671, 441, 430] & 4488.37s {\tiny (74.81 mins)} \\ % 16683555 (671, 542, 430)
			\cmidrule{2-6}
			\raisebox{-0.000ex}[0cm][0cm]{~\qquad/\texttt{D5}}
            & \FGMRES  & 4 & $10^{-3}$ & [5120$^{*}$, 115, 272, 143] & 12834.75s {\tiny (3.57 hours)} \\ % 16683558 (5120$^{*}_{1}$,  1412, 115)
            \midrule
			\raisebox{-0.000ex}[0cm][0cm]{~\qquad/from~\cite{Stanziola2021_HelmholtzEquationSolver}} 
            & Learned optimizer~\cite{Stanziola2021_HelmholtzEquationSolver}  & 4 & $10^{-3}$ &  [3084, 5120$^{*}$, 5120$^{*}$, 5120$^{*}$] & 922.22s {\tiny (15.37 mins)} \\  % 17518918
			\bottomrule
		\end{tabular}}
		\caption{{\small  Numerical results shown in Figure~\ref{fig:nn-fgmres_gen_sos} in terms of $\its$ and the CPUs~$\ET$ on \text{2x~16-core~Intel~Skylake~CPUs} for the four (i.e., $\K=4$) head CT scans problems with $\eta_b = 10^{-3}, 10^{-6}$ and $\maxit=5120$ with restrat equals 10. %(recall that the training time of this 2D U-Net on 4~NVIDIA~V100~GPUs is 45.53 min).
        }}\label{tab_results_skull}
\end{table}
%%***********************************************************

%******************************************************
\begin{figure}[!ht]%
    \centering%
    \begin{subfigure}{.49\textwidth}
        \includegraphics[width=\textwidth]{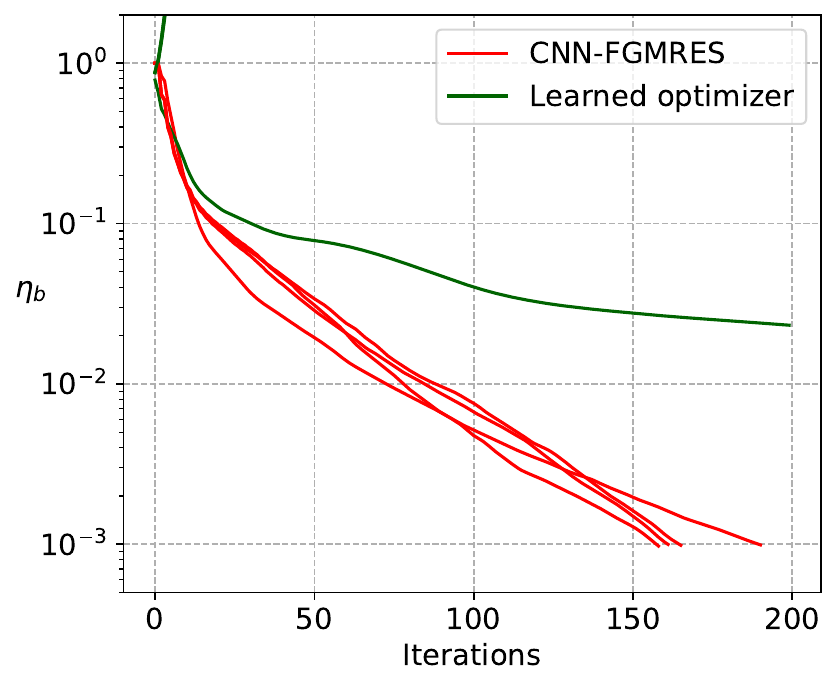}
        \caption{\FGMRES (from dataset \texttt{D0}) and Learned optimizer~\cite{Stanziola2021_HelmholtzEquationSolver} with $\K=4, \ \maxit=200$}
    \end{subfigure}\hfill
    \begin{subfigure}{.49\textwidth}%
        \includegraphics[width=\textwidth]{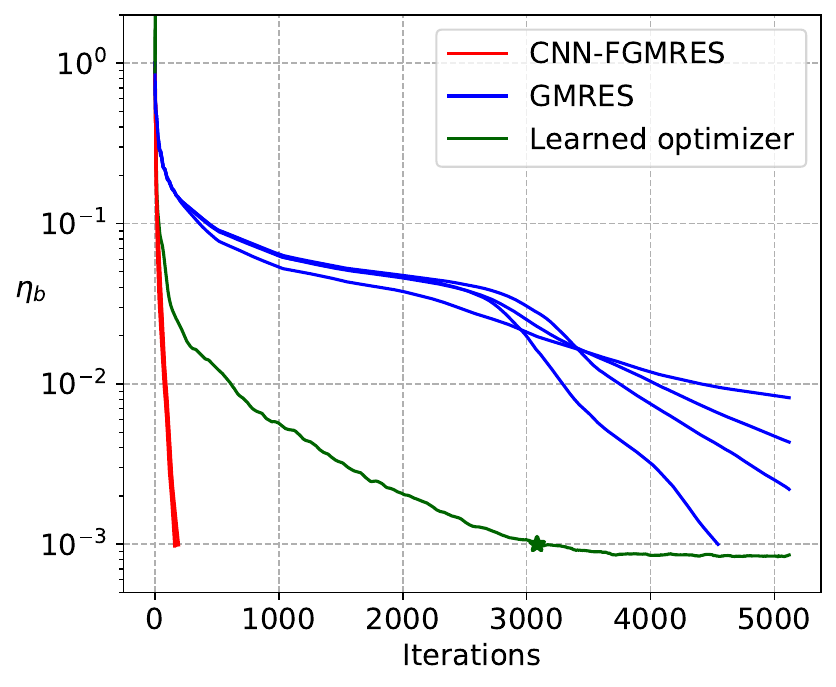}
        \caption{\FGMRES (from dataset \texttt{D0}), GMRES and Learned optimizer~\cite{Stanziola2021_HelmholtzEquationSolver} with $\K=4, \ \maxit=5120$}
    \end{subfigure}
    \vskip .1cm
    \medskip
    \begin{subfigure}{.49\textwidth}
        \includegraphics[width=\textwidth]{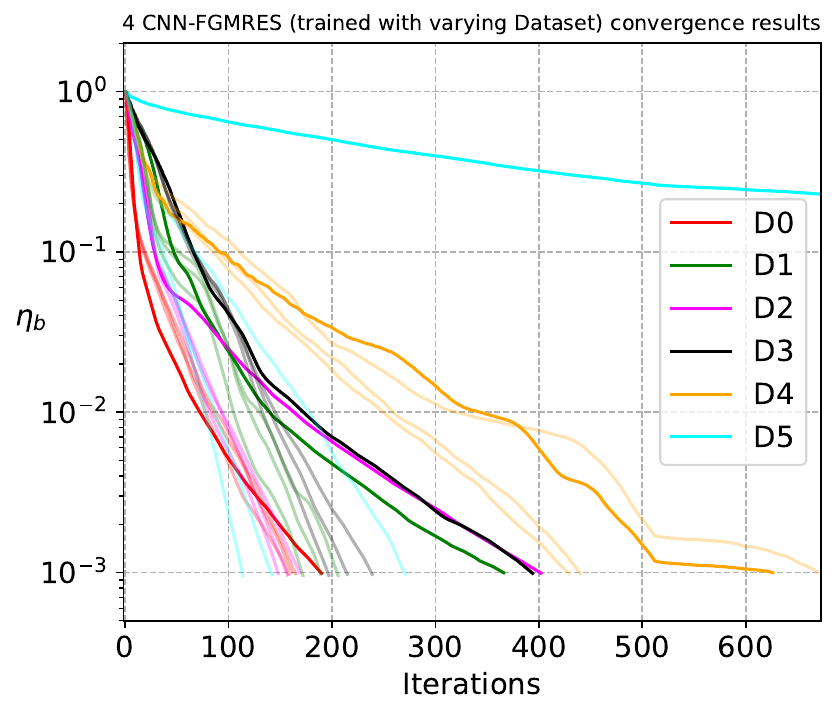}
        \caption{\FGMRES with $\K=4,$ and the CNN was trained with the six different datasets \texttt{D0}-\texttt{D5}} 
    \end{subfigure}\hfill
    \begin{subfigure}{.49\textwidth}
        \includegraphics[width=\textwidth]{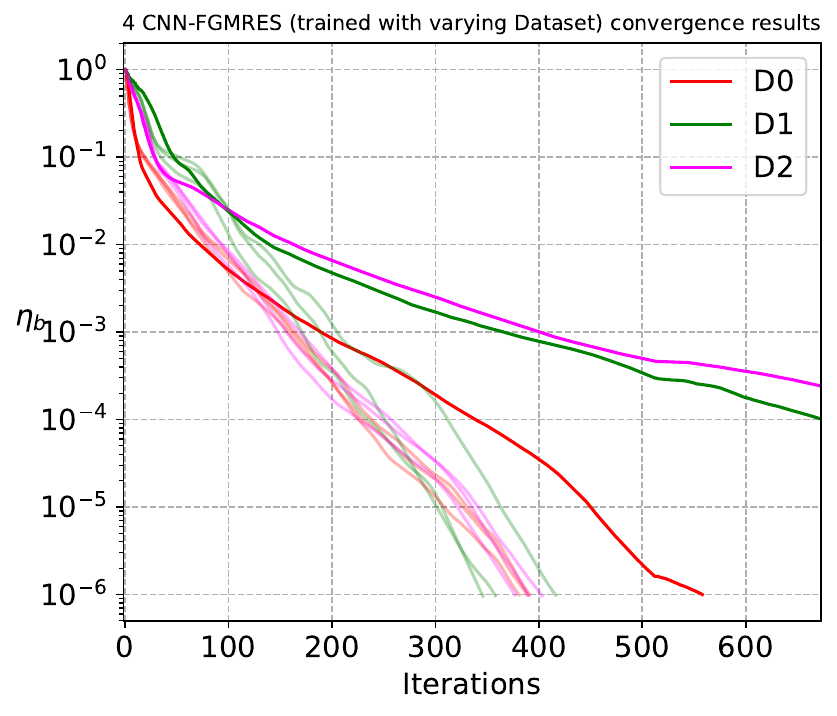}
        \caption{Counterpart of (c) with the three different datasets \texttt{D0}-\texttt{D2} and $\eta_b=10^{-6}$} 
    \end{subfigure}
    \caption{Convergence in terms of $\eta_b=10^{-3}$ for (a)-(c), and $\eta_b=10^{-6}$ for (d) of  \FGMRES and GMRES with $\maxit=5120$ with restrat equals 10, together with the learned optimizer proposed in~\cite{Stanziola2021_HelmholtzEquationSolver} with $\maxit=200, 5120,$ on four (i.e., $\K=4$) practical head CT scans problems. 
    In (a)-(b), the presented \FGMRES results are for the CNN-preconditioning trained with dataset \texttt{D0}. In (c), it shows corresponding results from the CNN-preconditioning trained with the six different datasets \texttt{D0}-\texttt{D5} described in Table~\ref{tab:cases_mixed_dataset}. The (d) is the counterpart results of (c) with datasets \texttt{D0}-\texttt{D2} and $\eta_b=10^{-6}.$}
    \label{fig:nn-fgmres_gen_sos}
\end{figure}
%******************************************************

Table~\ref{tab_results_skull} reports the numerical results of GMRES, the learned optimizer and \FGMRES for the four practical transcranial ultrasound examples, measured in terms of the iteration counts $\its$ and CPU elapsed time $\ET$ under different prescribed accuracy requirements.
Corresponding convergence histories are presented in Figure~\ref{fig:nn-fgmres_gen_sos}.
%s
The conclusions drawn from the previous experiments in terms of comparing \FGMRES to GMRES remain valid for these challenging large-scale test cases. In particular, the trained neural operator preconditioners consistently reduce both the number of iterations and the computational time compared with unpreconditioned GMRES.  
For instance, \FGMRES with the neural operator preconditioner trained on dataset \texttt{D0} requires at most 191~$\its$ with around 4~minutes for the most challenging example (i.e., $\bc$ from \texttt{CT000112.dcm} presented in Figure~\ref{fig:nn-fgmres_gen_sos}~(a)), and approximately 12~minutes in total to solve all four examples with accuracy $\eta_b=10^{-3}.$ 
In contrast, GMRES only successfully solves the simplest case (i.e., $\bc$ from \texttt{CT000255.dcm} presented in Figure~\ref{fig:nn-fgmres_gen_sos}~(b)), requiring 4546~$\its$, and fails to reach the prescribed accuracy for other three examples within the maximum iteration limit $\maxit.$\footnote{To further illustrate the difficulty of these problems, we note that, even when allowing additional restarts, GMRES requires 14,849~$\its$ and approximately 13.16~hours to solve the most challenging example shown in Figure~\ref{fig:nn-fgmres_gen_sos}~(a).}
On the other hand, we observed that the learned optimizer proposed in~\cite{Stanziola2021_HelmholtzEquationSolver} only successfully converges for one example shown in Figure~\ref{fig:nn-fgmres_gen_sos}~(a)-(b), which is the same test case (i.e., $\bc$ from \texttt{CT000112.dcm}) considered in~\cite[Section~3.3]{Stanziola2021_HelmholtzEquationSolver}, reaching a stagnated attainable accuracy of approximately $\eta_b=10^{-4}$ after 3084 Richardson iterations (noted by a green $*$ in the green curve in Figure~\ref{fig:nn-fgmres_gen_sos}~(b)). It fails to converge for the remaining three examples (see Figure~\ref{fig:nn-fgmres_gen_sos}~(a)-(b) or Figure~\ref{fig:learnt_optimizer_all_results} in Appendix~\ref{subsec:nn-fgmres_three_wave_evl}). In contrast, the proposed \FGMRES successfully solves all four examples with substantially fewer iterations and achieves arbitrary prescribed accuracy.

%******************************************************
\begin{figure}[!ht]
    \centering
    \begin{subfigure}{.33\textwidth}
        \includegraphics[width=\textwidth]{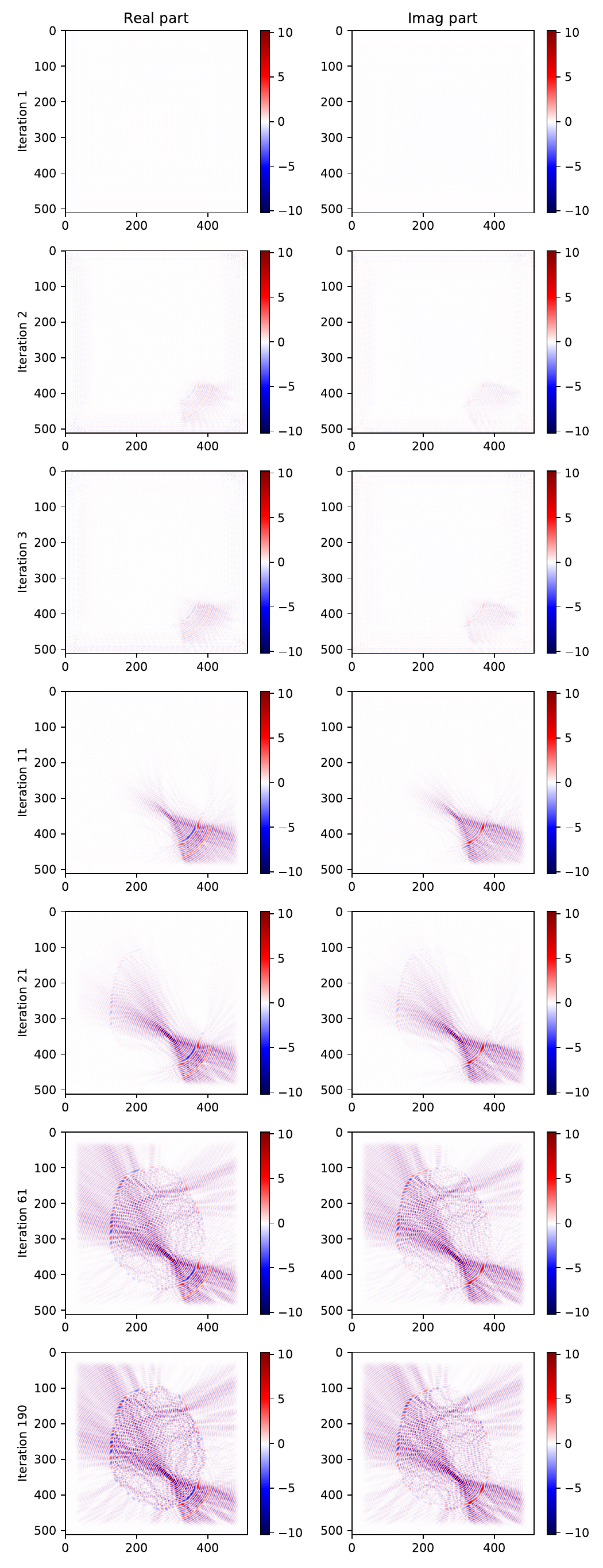}
        \caption{\FGMRES with dataset~\texttt{D0}}~\label{method_a} 
    \end{subfigure}\hfill
    \begin{subfigure}{.33\textwidth}
        \includegraphics[width=\textwidth]{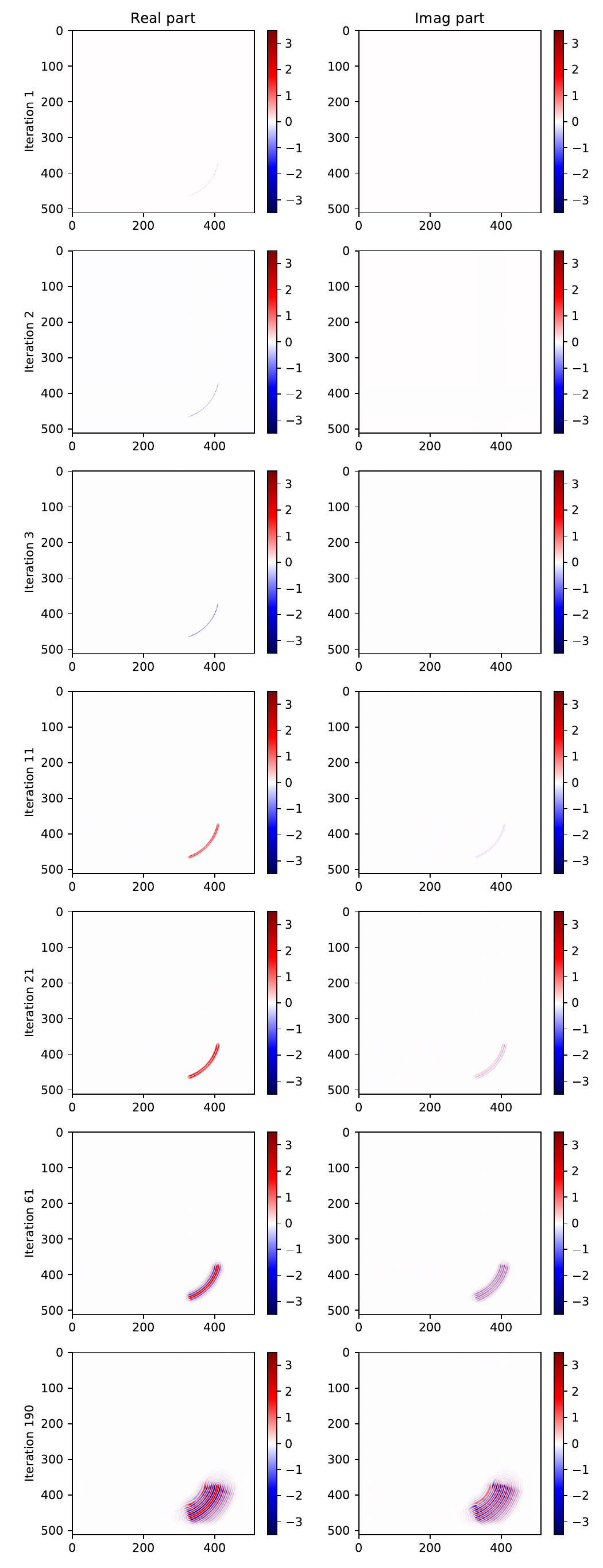}
        \caption{GMRES}~\label{method_b}
    \end{subfigure}\hfill
    \begin{subfigure}{.33\textwidth}
        \includegraphics[width=\textwidth]{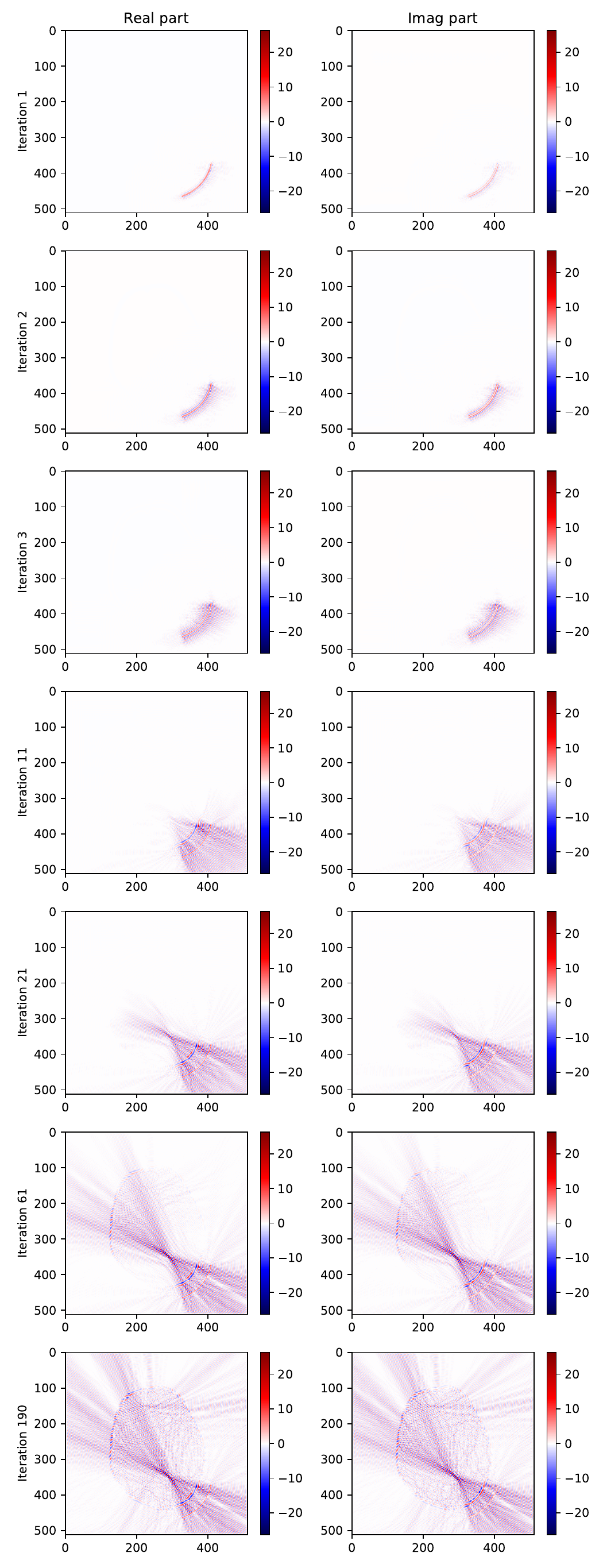} 
        \caption{Learned optimizer~\cite{Stanziola2021_HelmholtzEquationSolver}}~\label{method_c}
    \end{subfigure}
    \caption{Evolution of the complex wavefield predicted by the learned optimizer~\cite{Stanziola2021_HelmholtzEquationSolver}, \FGMRES and GMRES at iterations 1, 2, 3, 11, 21, 61, 190 for the transcranial ultrasound example based on the head CT scan $\bc$ from \texttt{CT000112.dcm} defined 2D domain $512 \times 512.$} %
    \label{fig:nn-fgmres_gmres}
\end{figure}
%******************************************************

The evolution of the complex wavefield predicted by GMRES, the learned optimizer and \FGMRES trained with dataset \texttt{D0} for the most challenging example with $\bc$ from \texttt{CT000112.dcm} is presented in Figure~\ref{fig:nn-fgmres_gmres}. The \FGMRES iterates rapidly recover the main wavefield features: the focal region and the corresponding focal pressure are already well established after approximately 21 iterations, while most of the fine-scale wavefield structure is captured within about 61 iterations. In contrast, GMRES does not exhibit comparable convergence within the same number of iterations. In fact, approximately 4000 iterations are required before the overall wavefield structure becomes discernible.
Furthermore, the learned optimizer requires more than 200 iterations before recovering most of the wavefield structure for the same transcranial ultrasound example (see Figure~\ref{fig:nn-fgmres_gmres}~(c) or~\cite[Figure~13]{Stanziola2021_HelmholtzEquationSolver}), and it fails to converge other three examples.\footnote{Refer to Figure~\ref{fig:nn-fgmres_three_wave_evl} in Appendix~\ref{subsec:nn-fgmres_three_wave_evl} for the evolution of the complex wavefields predicted by \FGMRES trained on dataset \texttt{D0} for the other three head CT scan examples.} 
This qualitative comparison further highlights the effectiveness and robustness of the proposed neural operator preconditioning strategy in accelerating the convergence of Krylov subspace methods for the Helmholtz equations arising in transcranial ultrasound simulations.
Figure~\ref{fig:nn-fgmres} provides visualizations of the resulting final complex wavefields together with the corresponding physical parameter configurations.

%******************************************************
\begin{figure}[!ht]
    \centering
    \begin{subfigure}{.49\textwidth}
        \includegraphics[width=\textwidth]{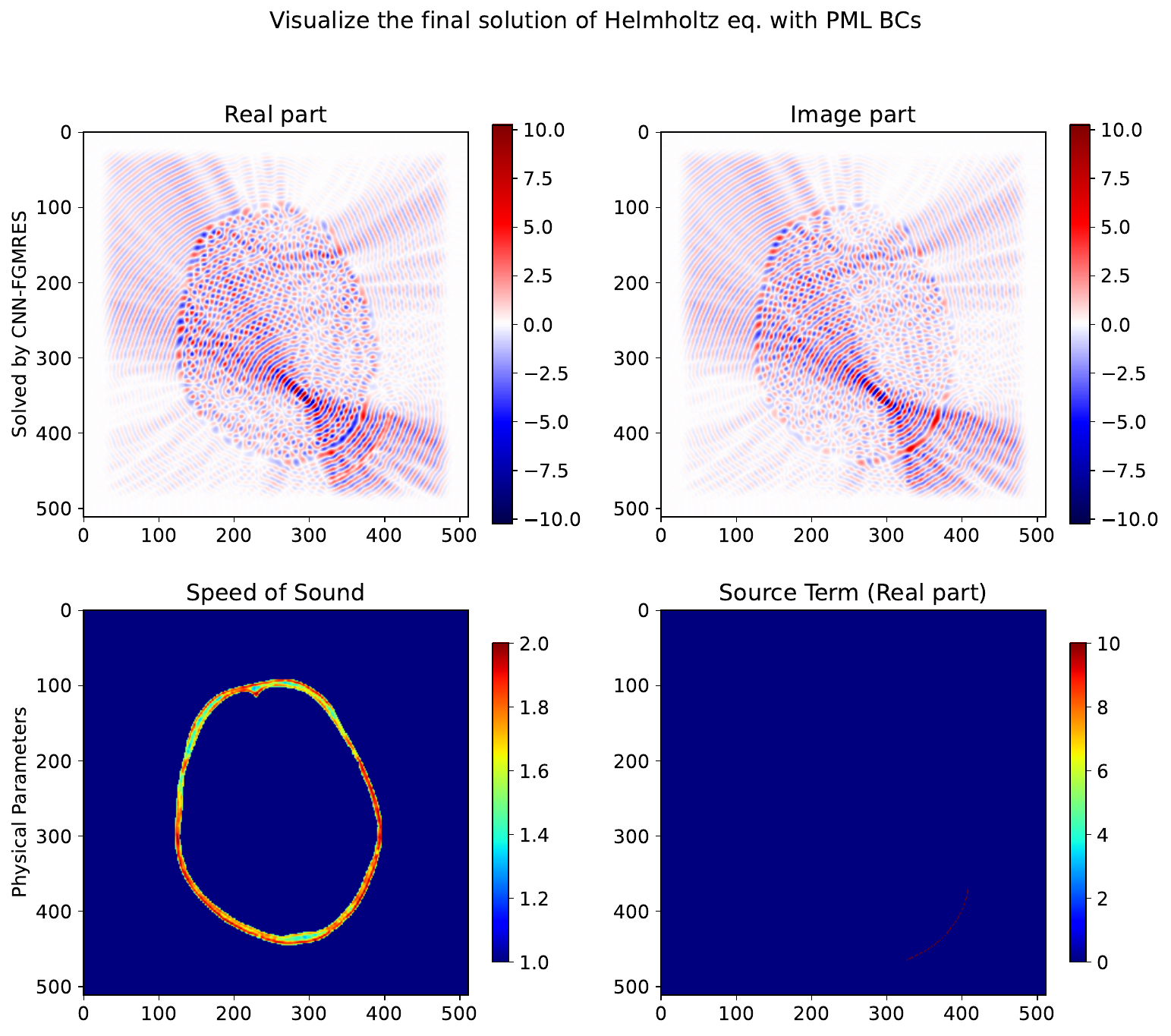}
        \caption{$\bc$ is from \texttt{CT000112.dcm}} 
    \end{subfigure}\hfill
    \begin{subfigure}{.49\textwidth}
        \includegraphics[width=\textwidth]{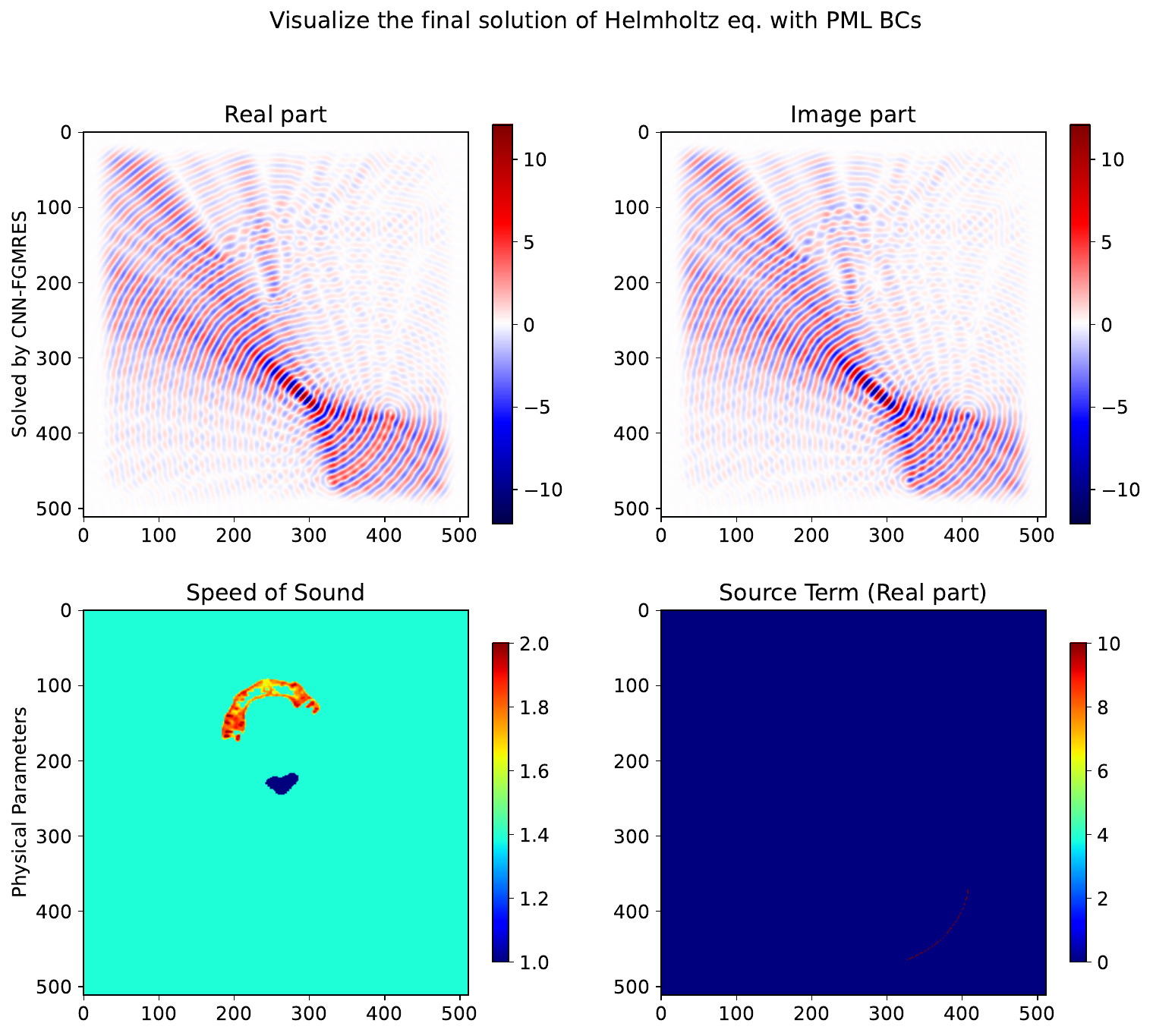}
        \caption{$\bc$ is from \texttt{CT000255.dcm}}
    \end{subfigure}
    \medskip
    \begin{subfigure}{.49\textwidth}
        \includegraphics[width=\textwidth]{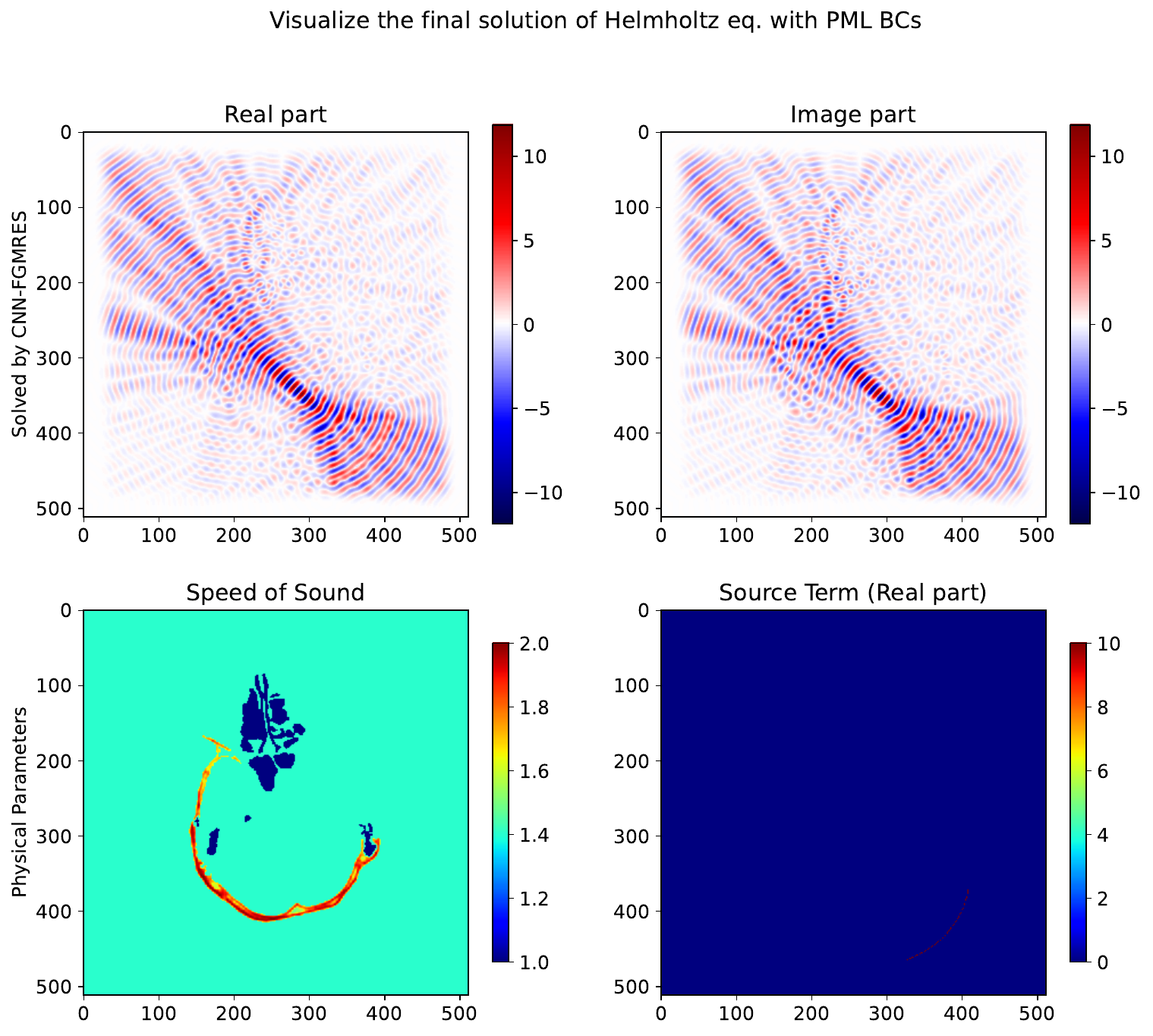}
        \caption{$\bc$ is from \texttt{CT000200.dcm}}
    \end{subfigure}\hfill
    \begin{subfigure}{.49\textwidth}
        \includegraphics[width=\textwidth]{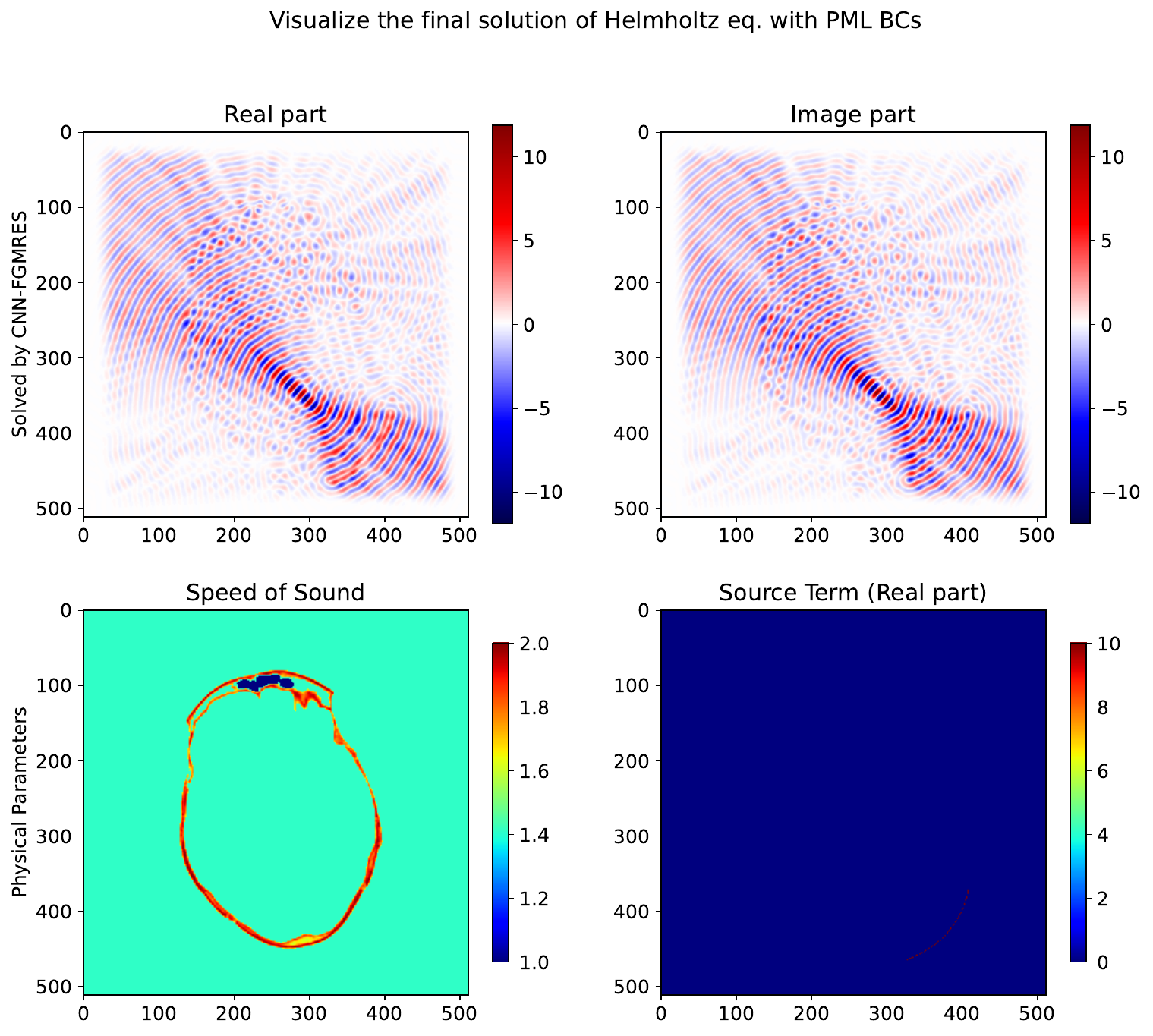}
        \caption{$\bc$ is from \texttt{CT000100.dcm}}
        % \label{fig:nn-fgmres-zoom}
    \end{subfigure}
    \caption{Visualize four head CT scans on 2D domain $512 \times 512$ solved by \FGMRES with dataset~\texttt{D0}, and $\eta_b=10^{-6}.$}
    \label{fig:nn-fgmres}
\end{figure}
%******************************************************

On the other hand, several additional observations can be made.
A key finding is that the performance differences among the six trained neural operator preconditioners become significantly more pronounced in these out-of-distribution testing scenarios. For example, comparing the models trained on datasets \texttt{D0} and \texttt{D1} reveals that increasing the diversity of the source field configurations can substantially improve efficiency, especially when a tighter accuracy requirement of $10^{-6}$ is imposed. Likewise, comparing the results obtained with datasets \texttt{D0} and \texttt{D2} highlights the influence of the velocity distributions used during training. The corresponding differences become increasingly evident as the testing problems depart further from the training configurations.
Furthermore, when comparing dataset \texttt{D1} with datasets \texttt{D4} and \texttt{D5}, which rely exclusively on either idealized skull configurations or fully random velocity fields, respectively, the model trained on \texttt{D1} consistently achieves better performance. This suggests that combining structured and heterogeneous training samples can be more beneficial than relying solely on highly specialized or completely randomized datasets.
Overall, these results indicate that training on specialized datasets can improve performance on closely related target problems, whereas training on more diverse datasets tends to enhance robustness and generalization across a broader range of scenarios. 
Among the six models considered, the preconditioner trained on dataset \texttt{D0} exhibits the best overall performance for these practical transcranial ultrasound simulations. A possible explanation is that \texttt{D0} achieves a favorable balance between specialized idealized skulls and randomized velocity fields while simultaneously incorporating the greatest diversity of source field configurations.
The convergence histories of \FGMRES with different network preconditoners and convergence threshold $\etab$ are presented in Figure~\ref{fig:nn-fgmres_gen_sos}~(c)-(d).

Up to this point, two important insights can be summarized from the perspective of generalization to practical simulations.
First, the composition of the training dataset plays a critical role in determining the effectiveness of learned models, particularly when they are applied to challenging out-of-distribution problems. Second, for scientific machine learning approaches targeting practical applications, training datasets should carefully balance application-specific structures and randomized samples to achieve both specialization and robust generalization.

We specially note that, beyond the primary application in transcranial ultrasound, the same trained neural operator preconditioners can also be applied to accelerate simulations in other settings. We refer the reader to Appendix~\ref{subsec:other_apps} for their use in simplified models of spherical wave propagation and seismology.

\section{Concluding remarks}~\label{sec:conclusion}
In this study, we investigated the use of neural networks trained on mixed datasets as nonlinear preconditioners for solving parametric Helmholtz equations with Krylov subspace methods. Six mixed training datasets were considered to construct neural operator preconditioners, with a primary focus on large-scale transcranial ultrasound simulations.
A systematic investigation of the training datasets demonstrates that dataset composition plays a critical role in determining the effectiveness and generalization capability of the learned preconditioners. In particular, a suitable balance between application-specific structured samples and randomized samples leads to improved performance on both in-distribution and out-of-distribution problems.
The proposed hybrid approach combines the strong generalization capability of neural operator preconditioners with the robustness of classical iterative methods, enabling accelerated solutions of large-scale Helmholtz problems to arbitrary prescribed tolerances. Numerical experiments demonstrate that the proposed neural operator preconditioners consistently accelerate Krylov subspace methods for both idealized and practical head~CT scan examples, while remaining effective for other Helmholtz-based applications.

Overall, this work highlights the importance of dataset design in scientific machine learning and provides practical guidelines for constructing training datasets for neural operator preconditioners targeting real-world large-scale PDE simulations. More broadly, it demonstrates the potential of integrating learned approximations of inverse PDE operators into established numerical linear algebra algorithms through a matrix-free framework for efficient large-scale simulations. Rather than replacing classical iterative solvers, neural operator preconditioners trained on mixed datasets serve as complementary components that provide strong generalization capability while retaining the flexibility and robustness of Krylov subspace methods.

\section*{Acknowledgements}
The author would like to thank Victor Michel-Dansac and Killian Lutz from Inria MACARON team and Universit\'e de Strasbourg for the valuable discussions and exchanges regarding the physical parameters settings and numerical experiments considered in Section~\ref{sec:raw_nn_gen_b}.

The author would also like to express her sincere gratitude to her Ph.D. supervisors, Luc Giraud (Inria Concace team, France) and Paul Mycek (Concace, Cerfacs, France). Research on transcranial ultrasound simulations was initiated during her doctoral studies under their supervision and continued to be developed after graduation. Their encouragement, support, and confidence in the potential of this research direction provided lasting motivation for pursuing the present study. This work is therefore, in many respects, a continuation of ideas and research interests that originated during the author's Ph.D.

% \section*{Declaration of generative AI and AI-assisted technologies in the manuscript preparation process}
% During the preparation of this work the author(s) used ChatGPT in order to improve language and readability. After using this tool/service, the author(s) reviewed and edited the content as needed and take(s) full responsibility for the content of the publication.

\printbibliography

%%%%%%%%%%%%%%
\newpage
\appendix

% \newpage
% \clearpage
\section{Further results on the four transcranial ultrasound examples}~\label{subsec:nn-fgmres_three_wave_evl}
%%***********************************************************

%******************************************************
\begin{figure}[!ht]
    \centering
    \begin{subfigure}{.285\textwidth}
        \includegraphics[width=\textwidth]{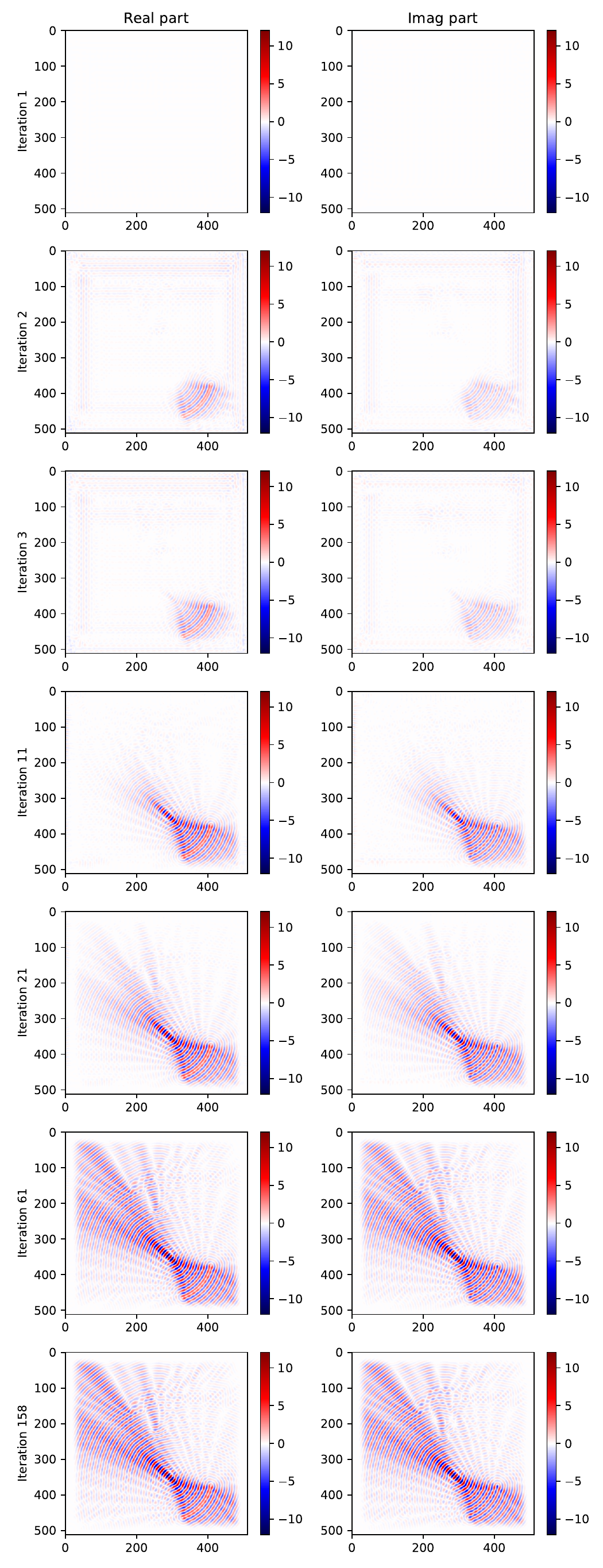}
        \caption{$\bc$ is from \texttt{CT000255.dcm}} 
    \end{subfigure}\hfill
    \begin{subfigure}{.285\textwidth}
        \includegraphics[width=\textwidth]{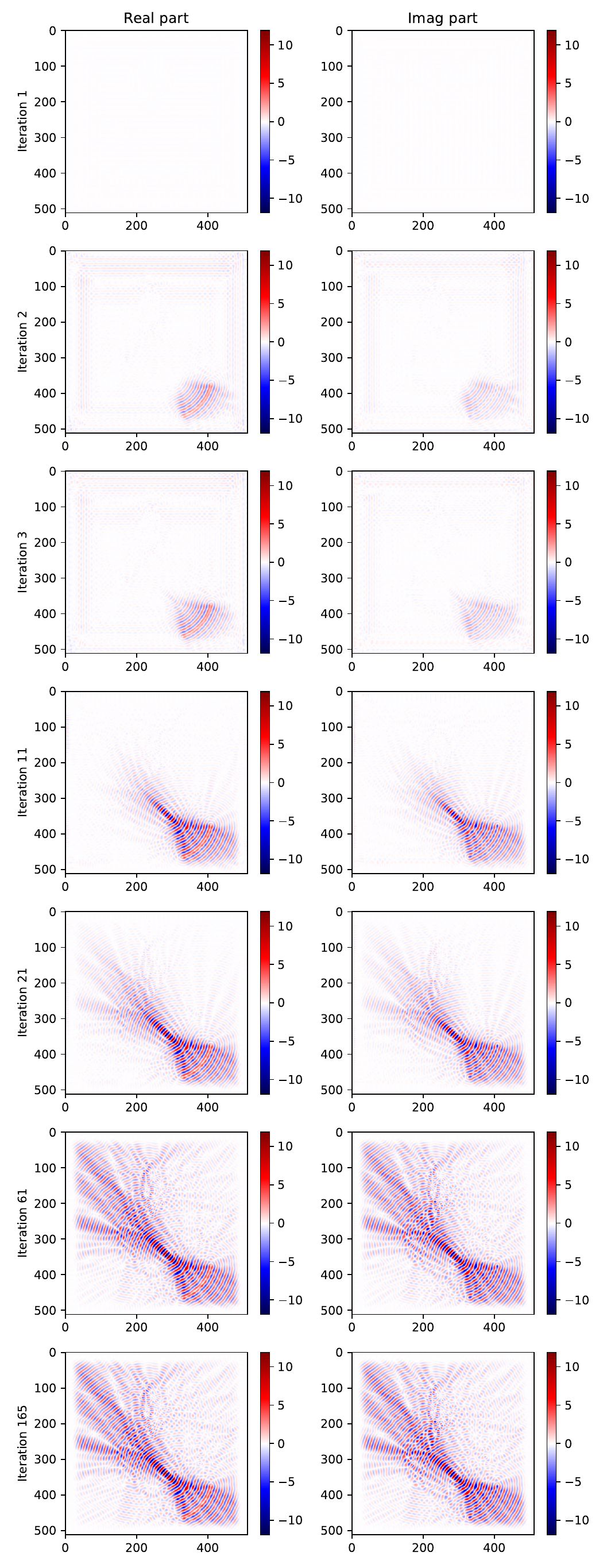} 
        \caption{$\bc$ is from \texttt{CT000200.dcm}}
    \end{subfigure}\hfill
    \begin{subfigure}{.285\textwidth}
        \includegraphics[width=\textwidth]{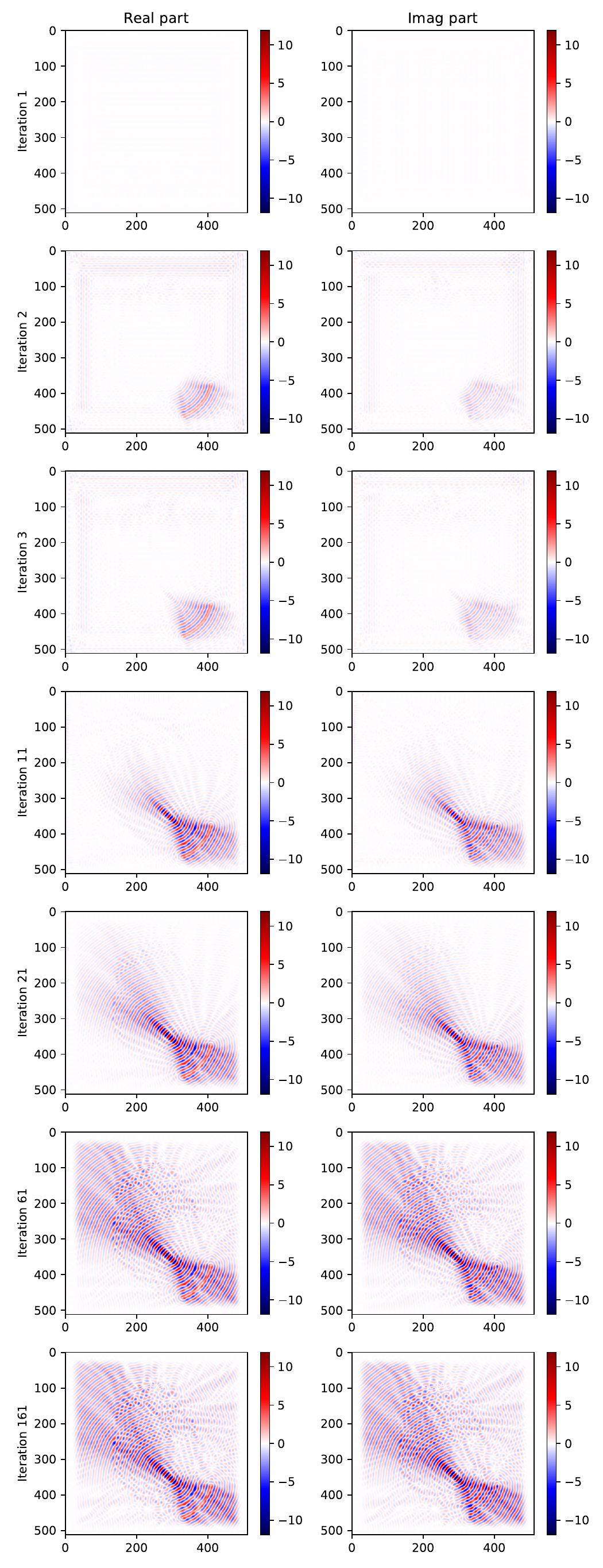} 
        \caption{$\bc$ is from \texttt{CT000100.dcm}}
    \end{subfigure}
    \caption{Counterpart of Figure~\ref{fig:nn-fgmres_gmres}: Evolution of the complex wavefield predicted by \FGMRES trained with dataset~\texttt{D0} at iterations 1, 2, 3, 11, 21, 61, and 158 for (a) or 165 for (b) 161 for (c) for the other three transcranial ultrasound examples defined 2D domain $512 \times 512.$} %
    \label{fig:nn-fgmres_three_wave_evl}
\end{figure}
%******************************************************

%******************************************************
\begin{figure}[!ht]%
    \centering%
    \begin{subfigure}{.35\textwidth}
        \includegraphics[width=\textwidth]{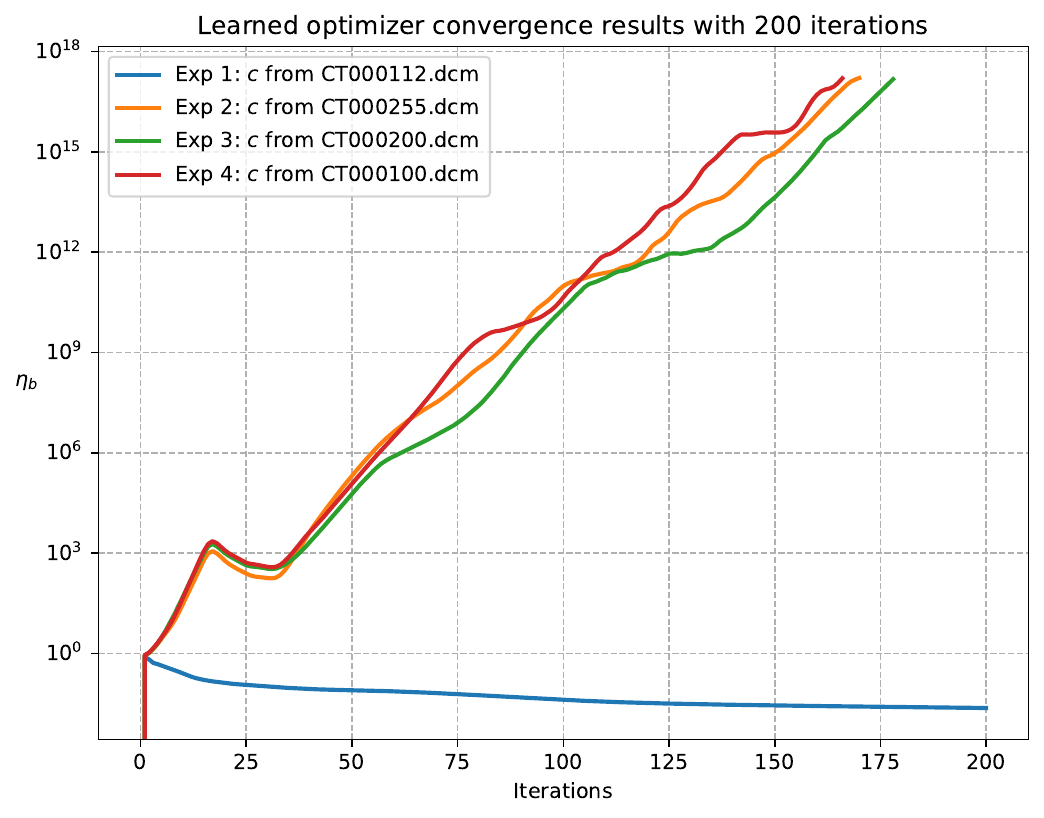}
        \caption{$\maxit=200$ Richardson iterations}
    \end{subfigure}\hfill
    \begin{subfigure}{.35\textwidth}%
    \includegraphics[width=\textwidth]{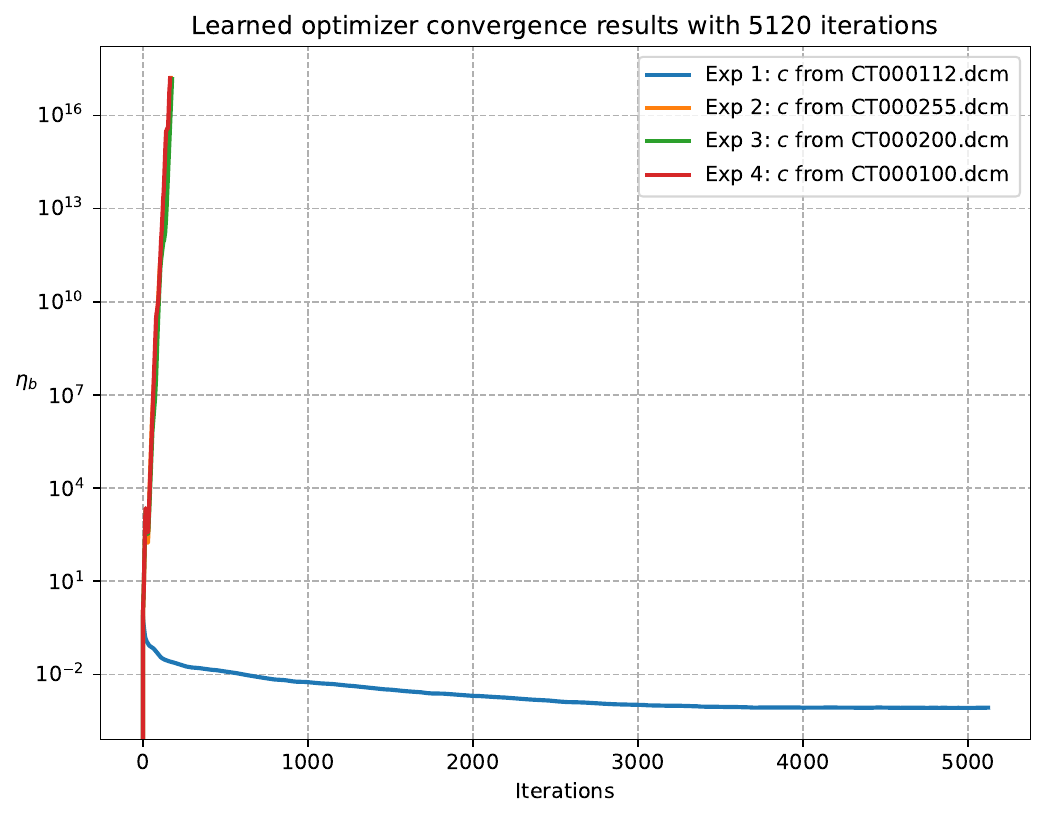}%
    \caption{$\maxit=5120$ Richardson iterations}%
    \end{subfigure}
    \caption{Convergence histories of the learned optimizer~\cite{Stanziola2021_HelmholtzEquationSolver} on the head CT scan examples in Section~\ref{subsec:skull_exp}.}
    \label{fig:learnt_optimizer_all_results}
\end{figure}
%******************************************************

\newpage
% \clearpage
\section{Simplified models of spherical wave propagation and seismology}~\label{subsec:other_apps}
%%***********************************************************
This section explores further applications of the same trained CNN-based preconditioner for solving Helmholtz equations arising from other applications, including the simplified models of spherical wave propagation with single or multiple source terms, as well as from idealized seismological settings involving wave propagation in heterogeneous media.
To represent these scenarios, we define the source term $\bb$ following the $\texttt{dirac1\_smooth}$ and $m$~\texttt{dirac} configurations, where $m \ge 1 \in \mathbb{N}$ is an integer. Specifically, $m$ entries of $\bb$ are set to $\pm 1$ at randomly selected grid locations, while all remaining entries are zero.
To model spherical wave propagation in a homogeneous medium, we first consider the \texttt{$\gamma$-constant} case, in which the wave speed $\bc$ is spatially constant. All entries of $\bc$ are set to a value $\gamma$ drawn uniformly at random from the interval $[1,2].$
To mimic more realistic seismological configurations, we construct heterogeneous wave speed fields $\bc$ with layered structures as shown in the publicly available OpenFWI dataset~\cite{NEURIPS2022_27d3ef26}. In each case, the layers are well aligned, and all elements within a given layer share a constant value sampled uniformly from $[\alpha, \beta]$ over the domain. Unless otherwise specified, we take $(\alpha, \beta)=(1, 2)$. The following configurations are considered:
\begin{enumerate}
	\item \texttt{4layered}: $\bc$ consists of four aligned layers, each with a constant value;
	\item \texttt{4faulted}: $\bc$ consists of four aligned but discontinuous (faulted) layers, with shifted segments.
\end{enumerate}
We refer the reader to the lower-left panel of each subplot of Figure~\ref{fig:viso_nn-fgmresgen_sos_dirac2} for visualizations of these homogeneous and heterogeneous wave speed fields with layered structures.
We note that these configurations are intended as simplified yet representative test cases for assessing the robustness of the proposed method under both homogeneous and heterogeneous wave propagation regimes.

\begin{figure}[!ht]%
    \centering%
    \begin{subfigure}{.21\textwidth}
        \includegraphics[width=\textwidth]{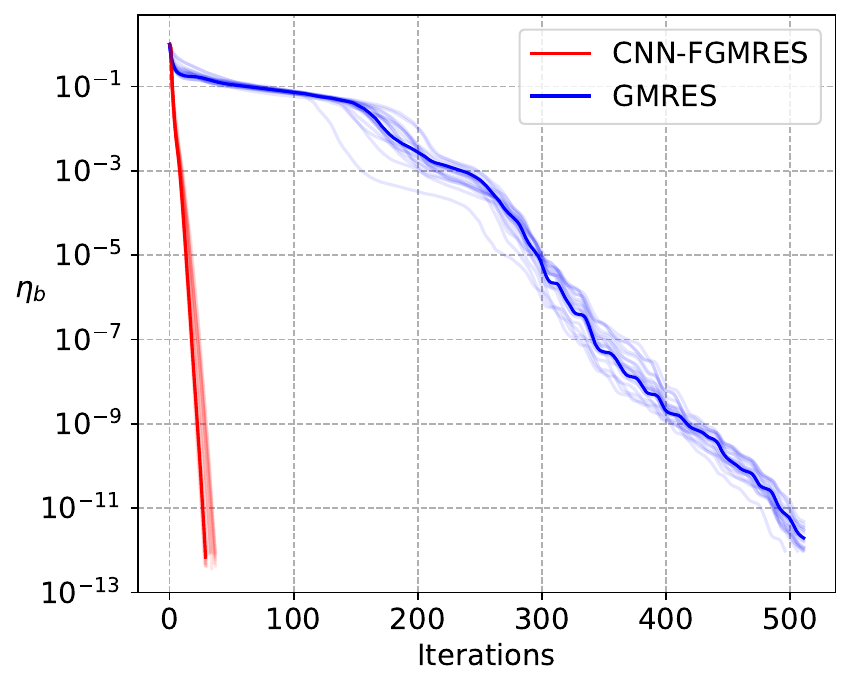}
        \caption{\scriptsize $\bc  \simiid \texttt{$\gamma$-constant},$ $\bb  \simiid ~\texttt{dirac},$ dataset \texttt{D3}}
    \end{subfigure}\hfill
    \begin{subfigure}{.21\textwidth}%
        \includegraphics[width=\textwidth]{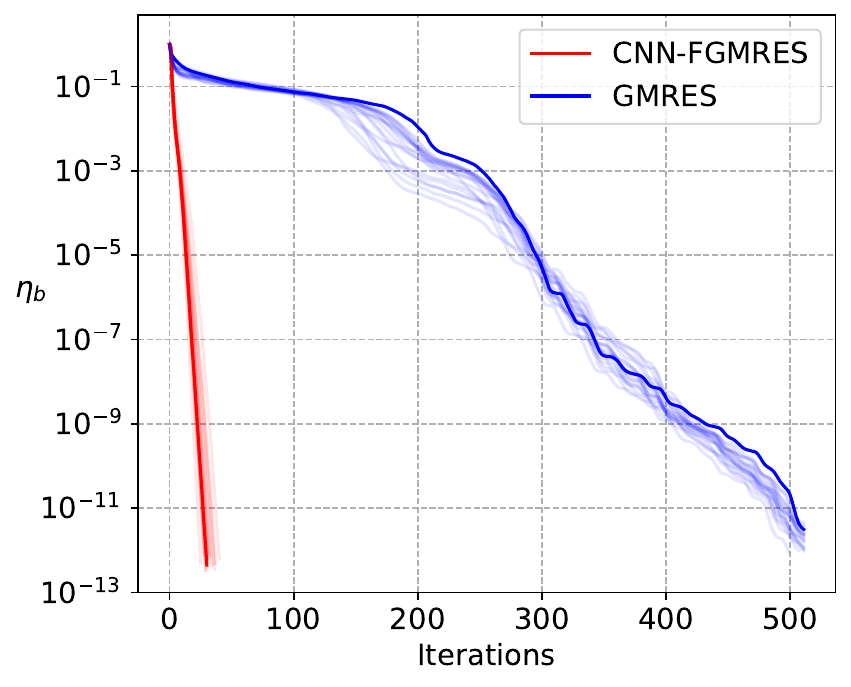}
        \caption{\scriptsize $\bc  \simiid \texttt{$\gamma$-constant},$ $\bb  \simiid 2~\texttt{dirac},$ dataset \texttt{D3}}%
    \end{subfigure}\hfill
    \begin{subfigure}{.21\textwidth}
        \includegraphics[width=\textwidth]{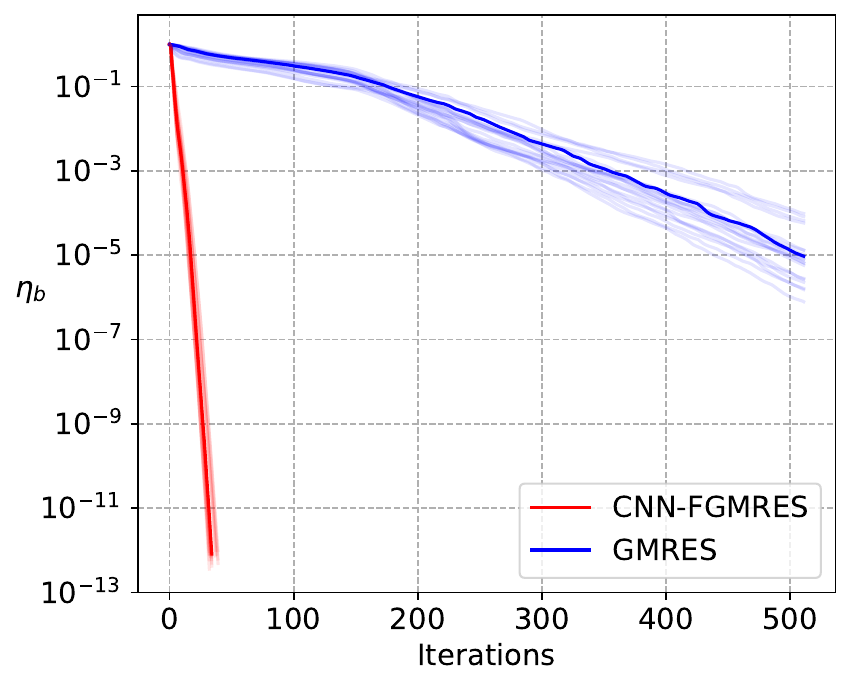}
        \caption{\scriptsize $\bc  \simiid \texttt{4layered},$ $\bb  \simiid \texttt{dirac1\_smooth},$ dataset \texttt{D3}}
    \end{subfigure}\hfill
    \begin{subfigure}{.21\textwidth}%
        \includegraphics[width=\textwidth]{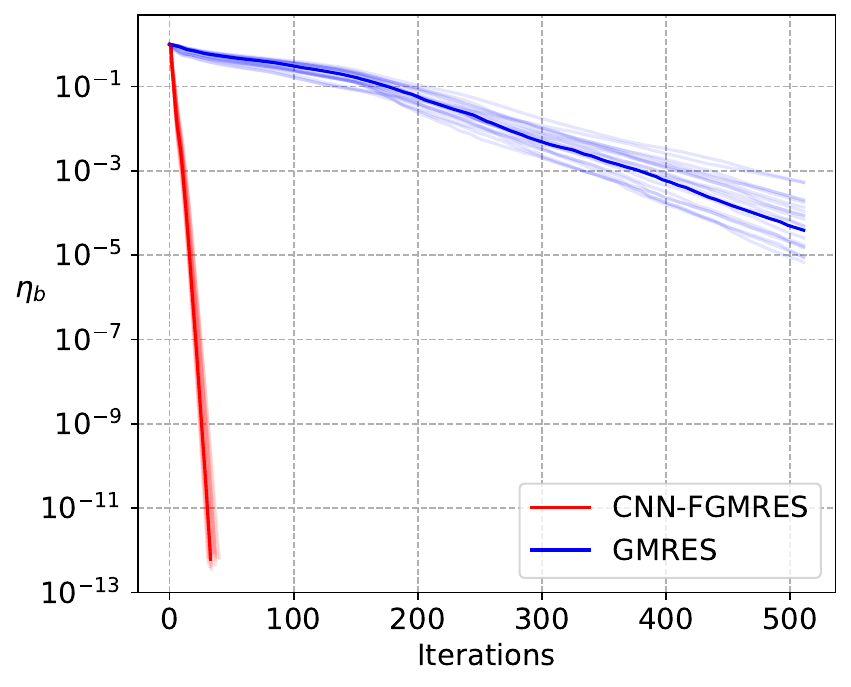}
        \caption{\scriptsize $\bc  \simiid \texttt{4faulted},$ $\bb  \simiid \texttt{dirac1\_smooth},$ dataset \texttt{D3}}%
    \end{subfigure}
    \medskip
    \begin{subfigure}{.21\textwidth}
        \includegraphics[width=\textwidth]{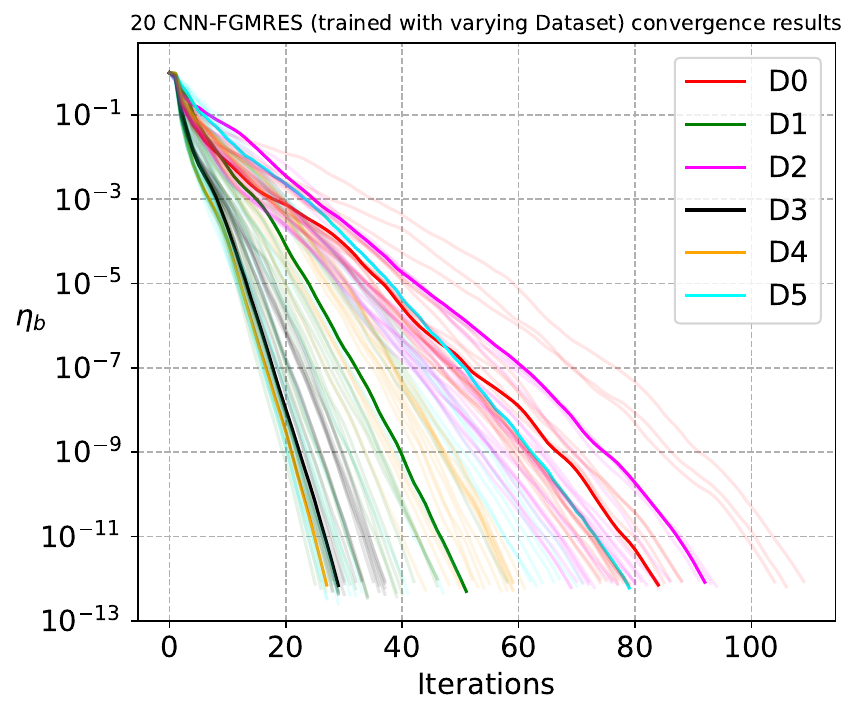}
        \caption{\scriptsize Counterpart of (a) with CNN from datasets \texttt{D0}-\texttt{D5}}%
    \end{subfigure}\hfill
    \begin{subfigure}{.21\textwidth}%
        \includegraphics[width=\textwidth]{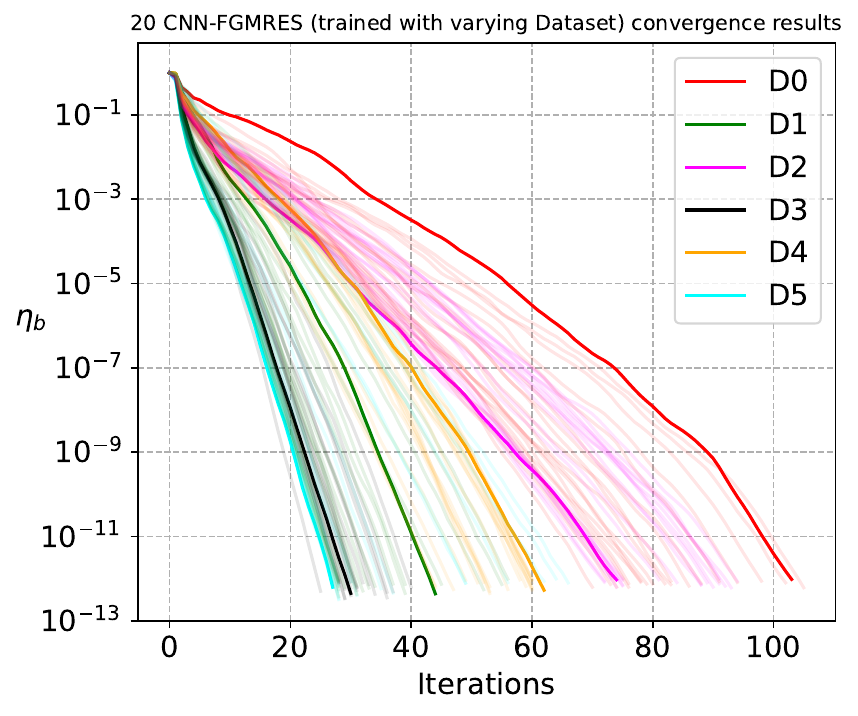}
        \caption{\scriptsize Counterpart of (b) with CNN from datasets \texttt{D0}-\texttt{D5}}%
    \end{subfigure}\hfill
    \begin{subfigure}{.21\textwidth}
        \includegraphics[width=\textwidth]{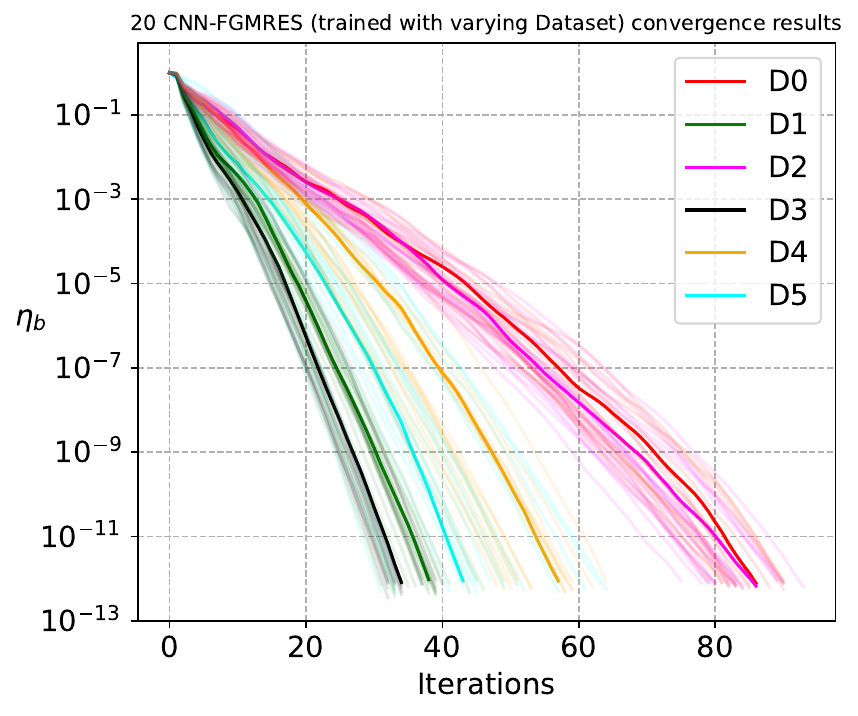}
        \caption{\scriptsize Counterpart of (c) with CNN from datasets \texttt{D0}-\texttt{D5}}%
    \end{subfigure}\hfill
    \begin{subfigure}{.21\textwidth}%
        \includegraphics[width=\textwidth]{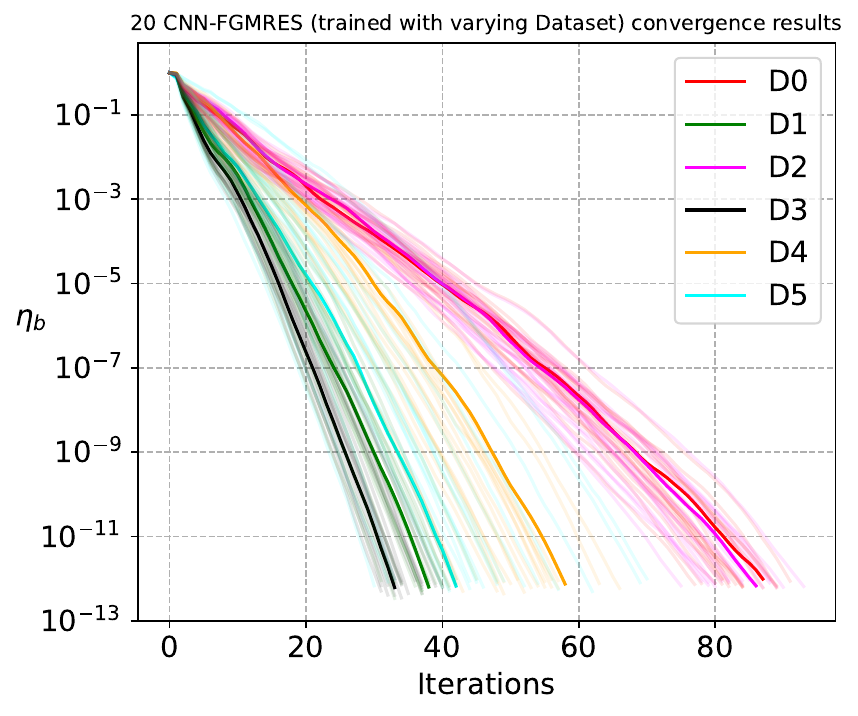}
        \caption{\scriptsize Counterpart of (d) with CNN from datasets \texttt{D0}-\texttt{D5}}%
    \end{subfigure}
    \caption{\small Convergence in terms of $\eta_b=10^{-12}$ of \FGMRES and GMRES with $\maxit=512$ without restart for $\K=20$ simplified models of spherical waves and seismology problems.  
    In (a)-(d), the presented \FGMRES results are for the CNN-preconditioning trained with dataset \texttt{D3}.
    In (e)-(h), it shows corresponding results from the CNN-preconditioning trained with the six different datasets \texttt{D0}-\texttt{D5} described in Table~\ref{tab:cases_mixed_dataset}.}
    \label{fig:nn-fgmres_gen_sos_dirac2}
\end{figure}

\begin{figure}[!ht]
    \centering
    \begin{subfigure}{.48\textwidth}
        \includegraphics[width=\textwidth]{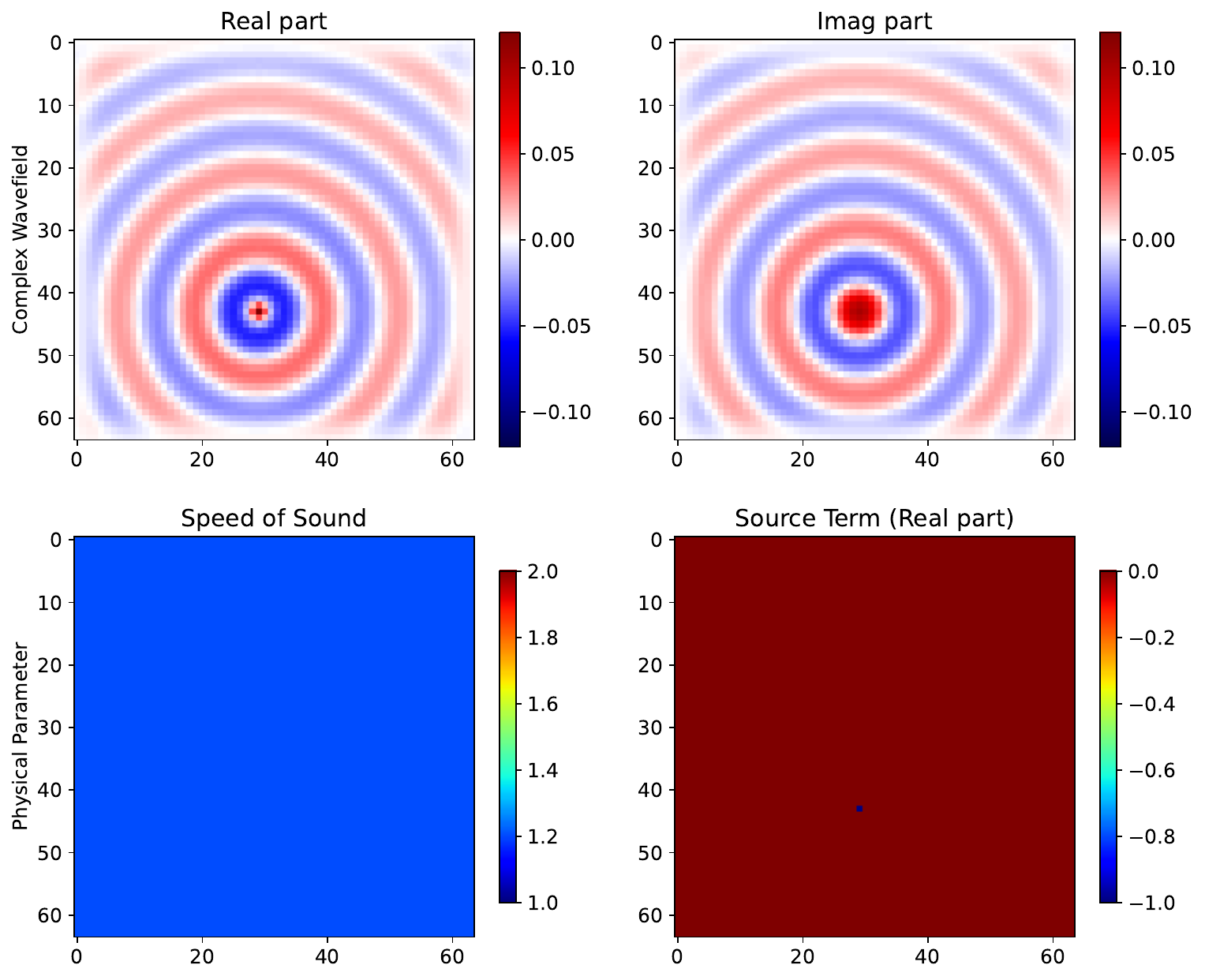}
        \caption{$\bc  \simiid \texttt{$\gamma$-constant}$ with $\bb  \simiid ~\texttt{dirac}$}
    \end{subfigure}\hfill
    \begin{subfigure}{.48\textwidth}%
        \includegraphics[width=\textwidth]{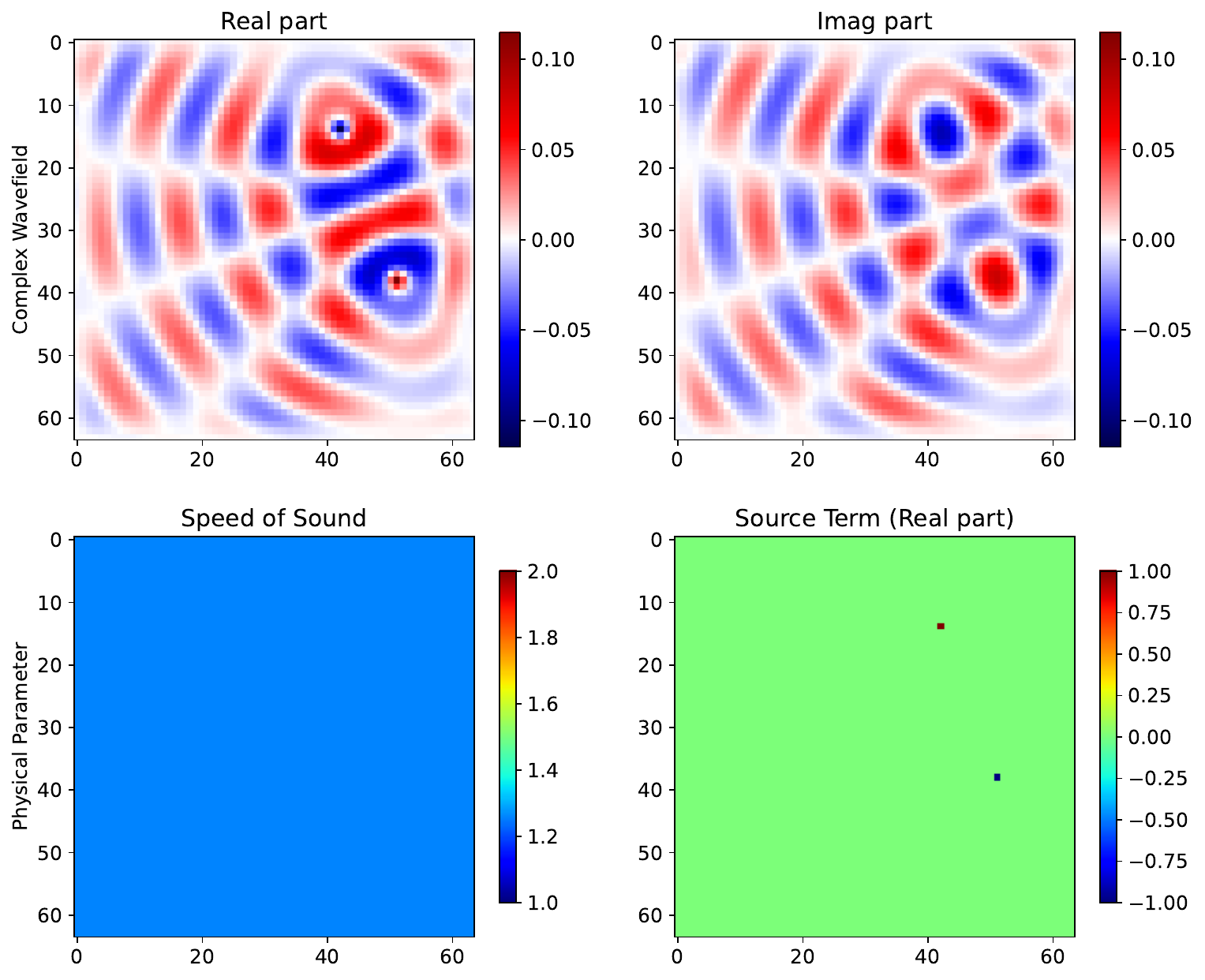} % 0
        \caption{$\bc  \simiid \texttt{$\gamma$-constant}$ with $\bb  \simiid 2~\texttt{dirac}$}%
    \end{subfigure}
    \vskip .1cm
    \medskip
    \begin{subfigure}{.48\textwidth}
        \includegraphics[width=\textwidth]{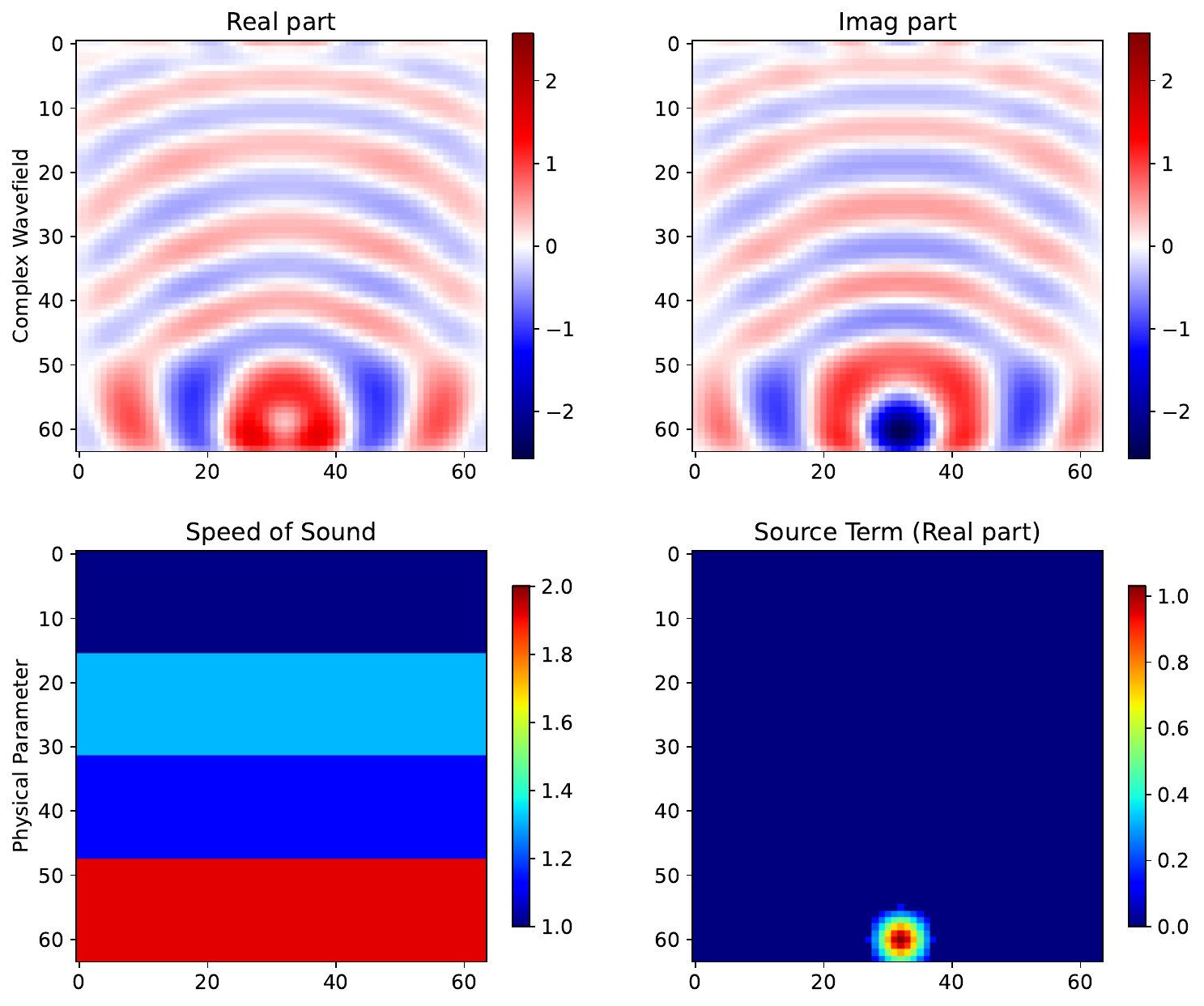} % 9, 8, 3, 1
        \caption{$\bc  \simiid \texttt{4layered}$ with $\bb  \simiid \texttt{dirac1\_smooth}$}
    \end{subfigure}\hfill
    \begin{subfigure}{.48\textwidth}%
        \includegraphics[width=\textwidth]{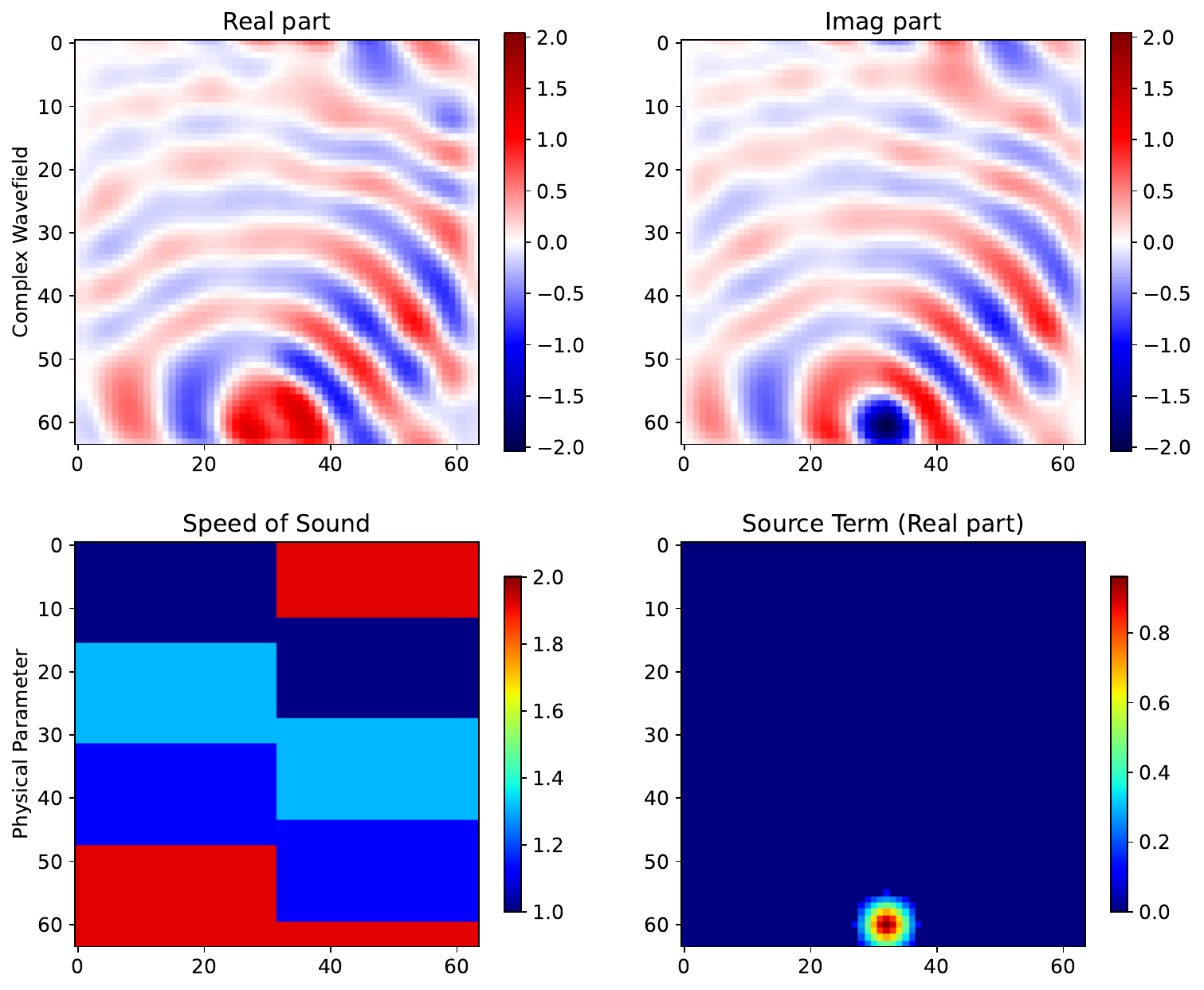}
        \caption{$\bc  \simiid \texttt{4faulted}$ with $\bb  \simiid \texttt{dirac1\_smooth}$}%
    \end{subfigure}
    \vskip .1cm
    \medskip
    \begin{subfigure}{.48\textwidth}
        \includegraphics[width=\textwidth]{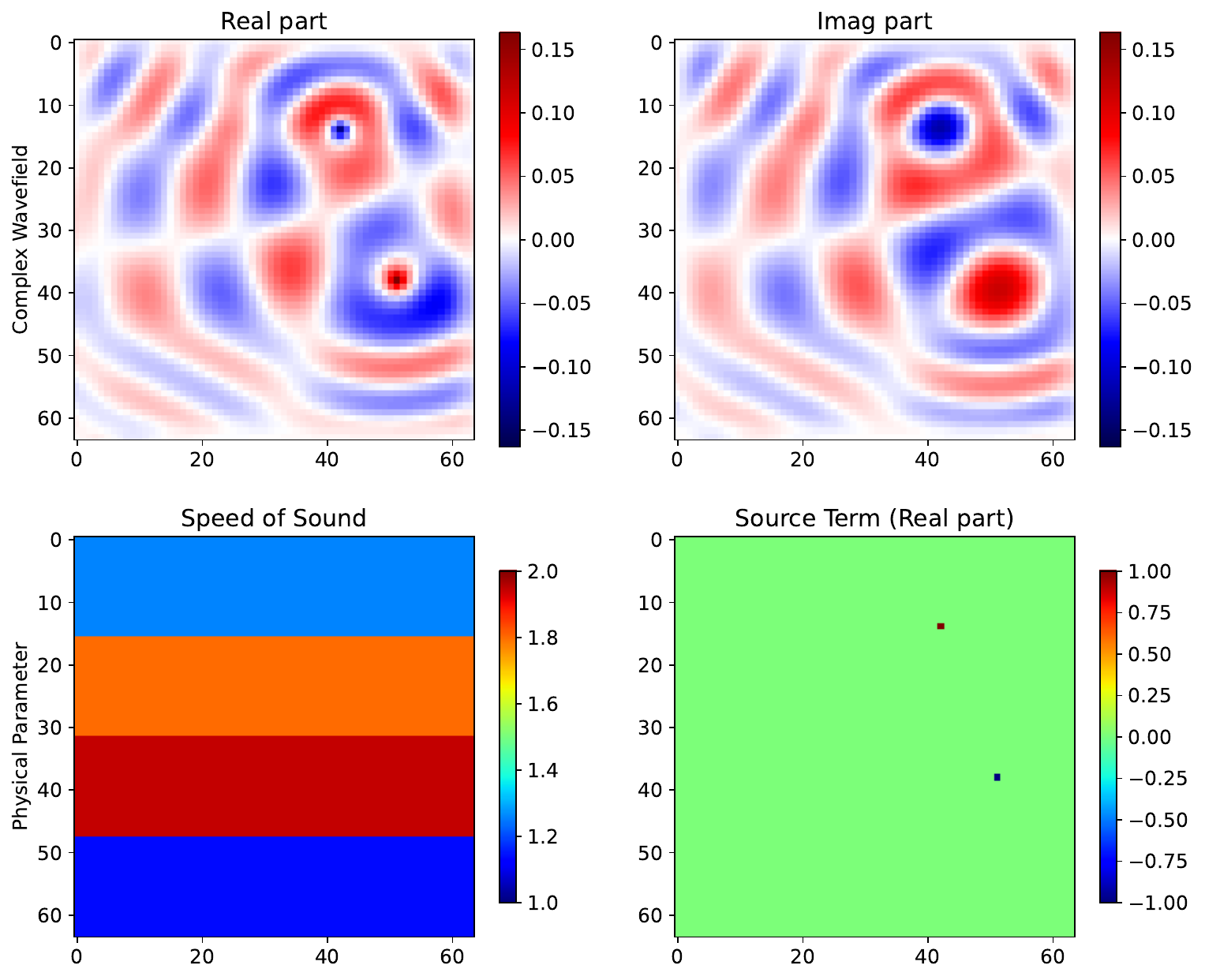}
        \caption{$\bc  \simiid \texttt{4layered}$ with $\bb  \simiid 2~\texttt{dirac}$}
    \end{subfigure}\hfill
    \begin{subfigure}{.48\textwidth}%
        \includegraphics[width=\textwidth]{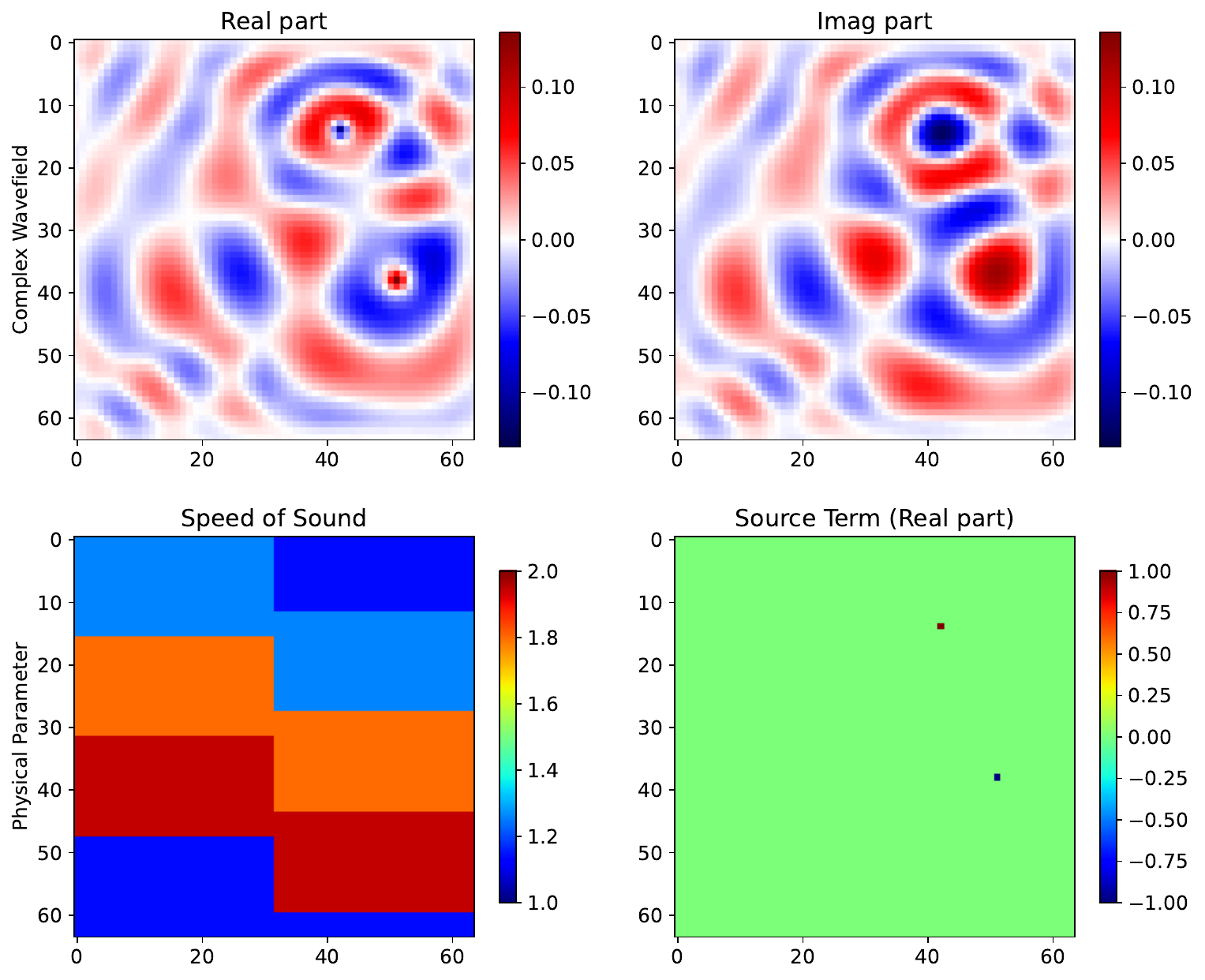}
        \caption{$\bc  \simiid \texttt{4faulted}$ with $\bb  \simiid 2~\texttt{dirac}$}%
    \end{subfigure}
    \caption{\small Visualize spherical waves at (a)-(b) and seismology examples (c)-(f) on 2D domain $64 \times 64$ solved by \FGMRES with $\eta_b=10^{-12}$ and dataset \texttt{D3}.}
    \label{fig:viso_nn-fgmresgen_sos_dirac2}
\end{figure}

The numerical results of GMRES and \FGMRES for the $\K=20$ linear systems with these physical parameter settings are reported in Table~\ref{tab_results_test_gen_sos_dirac2}. Among the six trained models, \FGMRES with network preconditioner trained under dataset~\texttt{D3} exhabites the best performance for these testing cases. 
The corresponding convergence histories for the model trained with dataset~\texttt{D3} are displayed in Figure~\ref{fig:nn-fgmres_gen_sos_dirac2}~(a)-(d), while the convergence histories for all six trained models are presented in Figure~\ref{fig:nn-fgmres_gen_sos_dirac2}~(e)-(h).
The conclusions drawn in the previous numerical sections remain valid in this section.
For completeness, Figure~\ref{fig:viso_nn-fgmresgen_sos_dirac2} provides visualizations of the resulting wavefields under some representative physical parameter configurations.

%******************************************************
%******************************************************

%***Begain Table 9*****************************************
%%***********************************************************
\begin{table}[!htbp]
	\center{
        \scriptsize
		\begin{tabular}
			{@{}llllrr@{}}
			\toprule
            \raisebox{-1.50ex}[0cm][0cm]{\begin{tabular}{c} \# \\ {Speed of sound}\end{tabular}}&
			\raisebox{-1.50ex}[0cm][0cm]{\begin{tabular}{c} \# \\ {Domain/Dataset}\end{tabular}}&
			\raisebox{-1.50ex}[0cm][0cm]{Method}&
			\raisebox{-1.50ex}[0cm][0cm]{$\K$}&
			\raisebox{-1.50ex}[0cm][0cm]{$(\its_{max}, \its_{avg}, \its_{min})$}&
			\raisebox{-1.50ex}[0cm][0cm]{$\ET$}\\
            &
            &
			&
			&
			&\\
			\midrule
            \raisebox{-13.00ex}[0cm][0cm]{\begin{tabular}{c} $\bc  \simiid \texttt{$\gamma$-constant}$ \\ with $\bb  \simiid ~\texttt{dirac}$ \end{tabular}}
			& \raisebox{-0.000ex}[0cm][0cm]{$64 \times 64$/}
            & GMRES & 20 & (512$^{*}_{15}$, 510, 497) &  2192.18s  \\  
            \cmidrule{3-6}
            \raisebox{-10.00ex}[0cm][0cm]{}
			& \raisebox{-0.000ex}[0cm][0cm]{~\qquad/\texttt{D0}}
            & \FGMRES  & 20 & (110, 85, 75) & 33.96s  \\
            \cmidrule{3-6}
			& \raisebox{-0.000ex}[0cm][0cm]{~\qquad/\texttt{D1}}
            & \FGMRES  & 20 & (59, 38, 26) & 8.05s \\ 
            \cmidrule{3-6}
			& \raisebox{-0.000ex}[0cm][0cm]{~\qquad/\texttt{D2}}
            & \FGMRES  & 20 & (95, 80, 70) &  25.30s  \\ 
            \cmidrule{3-6}
			& \raisebox{-0.000ex}[0cm][0cm]{~\qquad/\texttt{D3} $\star$}
            & \FGMRES  & 20 & (38, 33, 29) &  6.86s \\ 
            \cmidrule{3-6}
			&  \raisebox{-0.000ex}[0cm][0cm]{~\qquad/\texttt{D4}}
            & \FGMRES  & 20 & (62, 52, 28) &  11.61s \\ 
            \cmidrule{3-6}
			&  \raisebox{-0.000ex}[0cm][0cm]{~\qquad/\texttt{D5}}
            & \FGMRES  & 20 & (80, 48, 27) &  11.27s \\ 
            \midrule
            \raisebox{-13.00ex}[0cm][0cm]{\begin{tabular}{c} $\bc  \simiid \texttt{$\gamma$-constant}$ \\ with $\bb  \simiid 2~\texttt{dirac}$ \end{tabular}}
			& \raisebox{-0.000ex}[0cm][0cm]{$64 \times 64$/}
            & GMRES & 20 & (512$^{*}_{16}$, 511, 500) &  2130.23s \\  
			% \midrule
            \cmidrule{3-6}
            \raisebox{-10.00ex}[0cm][0cm]{}
			& \raisebox{-0.000ex}[0cm][0cm]{~\qquad/\texttt{D0}}
            & \FGMRES  & 20 & (106, 85, 71) &  23.36s \\ 
            \cmidrule{3-6}
			& \raisebox{-0.000ex}[0cm][0cm]{~\qquad/\texttt{D1}}
            & \FGMRES  & 20 & (57, 37, 29) & 10.02s  \\ 
            \cmidrule{3-6}
			& \raisebox{-0.000ex}[0cm][0cm]{~\qquad/\texttt{D2}}
            & \FGMRES  & 20 & (95, 83, 73) &  22.07s   \\
            \cmidrule{3-6}
			& \raisebox{-0.000ex}[0cm][0cm]{~\qquad/\texttt{D3} $\star$}
            & \FGMRES  & 20 & (41, 32, 26) & 6.87s  \\ 
            \cmidrule{3-6}
			&  \raisebox{-0.000ex}[0cm][0cm]{~\qquad/\texttt{D4}}
            & \FGMRES  & 20 & (63, 56, 45) &  12.81s  \\
            \cmidrule{3-6}
			&  \raisebox{-0.000ex}[0cm][0cm]{~\qquad/\texttt{D5}}
            & \FGMRES  & 20 & (67, 38, 28) &   9.50s \\ 
            \midrule
            \raisebox{-13.00ex}[0cm][0cm]{\begin{tabular}{c} $\bc  \simiid \texttt{4layered}$ \\ with $\bb  \simiid \texttt{dirac1\_smooth}$ \end{tabular}}
			& \raisebox{-0.00ex}[0cm][0cm]{$64 \times 64$/}
            & GMRES & 20 & (512, 512, 512)$^{*}_{20}$  &  520.62s   \\  
            \cmidrule{3-6}
            \raisebox{-10.00ex}[0cm][0cm]{}
			& \raisebox{-0.00ex}[0cm][0cm]{~\qquad/\texttt{D0}}
            & \FGMRES  & 20 &  (91, 85, 80)  & 19.70s  \\ 
            \cmidrule{3-6}
			& \raisebox{-0.000ex}[0cm][0cm]{~\qquad/\texttt{D1}}
            & \FGMRES  & 20 &  (46, 38, 32) &  8.28s  \\ 
            \cmidrule{3-6}
			& \raisebox{-0.000ex}[0cm][0cm]{~\qquad/\texttt{D2}}
            & \FGMRES  & 20 & (94, 84, 76)  &  19.94s  \\ % 
            \cmidrule{3-6}
			& \raisebox{-0.000ex}[0cm][0cm]{~\qquad/\texttt{D3} $\star$}
            & \FGMRES  & 20 &  (40, 35, 33) &  7.06s  \\ % 17086777
            \cmidrule{3-6}
			&  \raisebox{-0.000ex}[0cm][0cm]{~\qquad/\texttt{D4}}
            & \FGMRES  & 20 &  (65, 55, 47) &  11.84s  \\ % 
            \cmidrule{3-6}
			&  \raisebox{-0.000ex}[0cm][0cm]{~\qquad/\texttt{D5}}
            & \FGMRES  & 20 & (65, 48, 32) &  10.18s  \\ % 
            \midrule
            \raisebox{-13.00ex}[0cm][0cm]{\begin{tabular}{c} $\bc  \simiid \texttt{4faulted}$ \\ with $\bb  \simiid \texttt{dirac1\_smooth}$ \end{tabular}}
			& \raisebox{-0.00ex}[0cm][0cm]{$64 \times 64$/}
            & GMRES & 20 & (512, 512, 512)$^{*}_{20}$  & 2379.37s  \\  
            \cmidrule{3-6}
            \raisebox{-10.00ex}[0cm][0cm]{}
			& \raisebox{-0.00ex}[0cm][0cm]{~\qquad/\texttt{D0}}
            & \FGMRES  & 20 & (92, 85, 78) & 19.51s \\ 
            \cmidrule{3-6}
			& \raisebox{-0.000ex}[0cm][0cm]{~\qquad/\texttt{D1}}
            & \FGMRES  & 20 &  (43, 38, 32) & 7.99s  \\ % 
            \cmidrule{3-6}
			& \raisebox{-0.000ex}[0cm][0cm]{~\qquad/\texttt{D2}}
            & \FGMRES  & 20 & (94, 84, 76)  & 19.67s  \\ % 
            \cmidrule{3-6}
			& \raisebox{-0.000ex}[0cm][0cm]{~\qquad/\texttt{D3} $\star$}
            & \FGMRES  & 20 & (41, 35, 32)  &  7.08s  \\ 
            \cmidrule{3-6}
			&  \raisebox{-0.000ex}[0cm][0cm]{~\qquad/\texttt{D4}}
            & \FGMRES  & 20 & (67, 54, 44)  & 11.64s   \\ % 
            \cmidrule{3-6}
			&  \raisebox{-0.000ex}[0cm][0cm]{~\qquad/\texttt{D5}}
            & \FGMRES  & 20 & (71, 48, 31)  & 10.50s \\ % 
            \midrule
            % %% seisemic test cases with two sources
            \raisebox{-13.00ex}[0cm][0cm]{\begin{tabular}{c} $\bc  \simiid \texttt{4layered}$ \\ with $\bb  \simiid 2~\texttt{dirac}$ \end{tabular}}
			& \raisebox{-0.00ex}[0cm][0cm]{$64 \times 64$/}
            & GMRES & 20 & (512, 512, 512)$^{*}_{20}$ & 818.80s  \\  
            \cmidrule{3-6}
            \raisebox{-10.00ex}[0cm][0cm]{}
			& \raisebox{-0.00ex}[0cm][0cm]{~\qquad/\texttt{D0}}
            & \FGMRES  & 20 & (92, 86, 77) & 20.35s \\ 
            \cmidrule{3-6}
			& \raisebox{-0.000ex}[0cm][0cm]{~\qquad/\texttt{D1}}
            & \FGMRES  & 20 & (47, 38, 33) &  8.23s \\ 
            \cmidrule{3-6}
			& \raisebox{-0.000ex}[0cm][0cm]{~\qquad/\texttt{D2}}
            & \FGMRES  & 20 & (93, 82, 74) &  19.32s  \\ 
            \cmidrule{3-6}
			& \raisebox{-0.000ex}[0cm][0cm]{~\qquad/\texttt{D3} $\star$}
            & \FGMRES  & 20 & (42, 35, 32) &  7.39s \\ 
            \cmidrule{3-6}
			&  \raisebox{-0.000ex}[0cm][0cm]{~\qquad/\texttt{D4}}
            & \FGMRES  & 20 & (66, 54, 48) &  15.68s \\ 
            \cmidrule{3-6}
			&  \raisebox{-0.000ex}[0cm][0cm]{~\qquad/\texttt{D5}}
            & \FGMRES  & 20 & (63, 49, 34) & 10.47s  \\ 
            \midrule
            \raisebox{-13.00ex}[0cm][0cm]{\begin{tabular}{c} $\bc  \simiid \texttt{4faulted}$ \\ with $\bb  \simiid 2~\texttt{dirac}$ \end{tabular}}
			& \raisebox{-0.00ex}[0cm][0cm]{$64 \times 64$/}
            & GMRES & 20 & (512, 512, 512)$^{*}_{20}$ & 2133.76s  \\  
            \cmidrule{3-6}
            \raisebox{-10.00ex}[0cm][0cm]{}
			& \raisebox{-0.00ex}[0cm][0cm]{~\qquad/\texttt{D0}}
            & \FGMRES  & 20 & (94, 84, 77)  &  22.37s \\ 
            \cmidrule{3-6}
			& \raisebox{-0.000ex}[0cm][0cm]{~\qquad/\texttt{D1}}
            & \FGMRES  & 20 & (46, 37, 33) & 7.73s \\ 
            \cmidrule{3-6}
			& \raisebox{-0.000ex}[0cm][0cm]{~\qquad/\texttt{D2}}
            & \FGMRES  & 20 & (91, 82, 74) &  21.57s  \\ 
            \cmidrule{3-6}
			& \raisebox{-0.000ex}[0cm][0cm]{~\qquad/\texttt{D3} $\star$}
            & \FGMRES  & 20 & (40, 34, 31) &  7.47s  \\ 
            \cmidrule{3-6}
			&  \raisebox{-0.000ex}[0cm][0cm]{~\qquad/\texttt{D4}}
            & \FGMRES  & 20 & (63, 54, 43) &  11.64s \\ 
            \cmidrule{3-6}
			&  \raisebox{-0.000ex}[0cm][0cm]{~\qquad/\texttt{D5}}
            & \FGMRES  & 20 & (78, 50, 31)  & 10.63s  \\ 
			\bottomrule
		\end{tabular}}
		\caption{{\small  Numerical results shown in Figure~\ref{fig:nn-fgmres_gen_sos_dirac2} 
        in terms of $(\its_{max}, \its_{avg}, \its_{min})$ and the GPUs~$\ET$ on 4~V100~GPUs for \text{$\K=20$} examples with varying wave speed $\bc  \simiid \texttt{$\gamma$-constant},$ and the heterogeneous ones with different types of layered structures, and the random source term \text{$\bb  \simiid m~\text{\texttt{dirac}}$ and $\bb \simiid \texttt{dirac1\_smooth}.$} The $\eta_b = 10^{-12}$ and $\maxit=512$ without restrat.}}\label{tab_results_test_gen_sos_dirac2}
\end{table}
%%***********************************************************

%******************************************************
%******************************************************

\end{document}